\documentclass[reqno,twoside,11pt]{amsart}
\usepackage{cite}
\usepackage{amsmath}
\usepackage{amsfonts}
\usepackage{amssymb}
\usepackage{epsfig}
\usepackage{verbatim}
\usepackage{color}

\usepackage{ulem}

\usepackage{mathpazo}

\DeclareFontFamily{OT1}{cmss}{} \DeclareFontShape{OT1}{cmss}{m}{n} {<5> <6> <7> <8> <9> <10> <11> <12> <13> <14.4> cmss10}{}
\DeclareMathAlphabet{\cmss}{OT1}{cmss}{m}{n}

\DeclareFontFamily{OT1}{fraktura}{}
\DeclareFontShape{OT1}{fraktura}{m}{n} {<5> <6> <7> <8> <9> <10> <11> <12> <13> <14.4> [1.1] eufm10}{}
\DeclareMathAlphabet{\fraktura}{OT1}{fraktura}{m}{n}

\newtheoremstyle{thm}{1.5ex}{1.5ex}{\itshape\rmfamily}{} {\bfseries\rmfamily}{}{2ex}{}

\newtheoremstyle{def}{1.5ex}{1.5ex}{\rmfamily\sl}{} {\bfseries\rmfamily}{}{2ex}{}

\newtheoremstyle{rem}{1.3ex}{1.3ex}{\rmfamily}{} {\bfseries\rmfamily}{}{2ex}{}

\newtheoremstyle{ass}{1.5ex}{1.5ex}{\rmfamily\sl}{} {\bfseries\rmfamily}{}{2ex}{}

\newenvironment{proofsect}[1] {\vskip0.1cm\noindent{\rmfamily\itshape#1.}}{\qed\vspace{0.15cm}}

\theoremstyle{thm}
\newtheorem{theorem}{Theorem}[section]
\newtheorem{lemma}[theorem]{Lemma}

\newtheorem*{Main Theorem}{Main Theorem.}
\newtheorem{corollary}[theorem]{Corollary}

\newtheorem{definition}[theorem]{Definition}

\theoremstyle{rem}
\newtheorem{remark}[theorem]{{Remark}}

\numberwithin{equation}{section}

\renewcommand{\section}{\secdef\sct\sect}
\newcommand{\sct}[2][default]{\refstepcounter{section}
\addcontentsline{toc}{section}
{{\tocsection {}{\thesection}{\!\!\!\!#1\dotfill}}{}}
\vspace{0.7cm}
\centerline{ 
\scshape\arabic{section}.\ #1} \nopagebreak \vspace{0.2cm}}
\newcommand{\sect}[1]{
\vspace{0.4cm} \centerline{\large\scshape\rmfamily #1}
\vspace{0.2cm}}

\renewcommand{\subsection}{\secdef\subsct\sbsect}
\newcommand{\subsct}[2][default]{\refstepcounter{subsection}
\addcontentsline{toc}{subsection}
{{\tocsection{\!\!}{\hspace{1.2em}\thesubsection}{\!\!\!\!#1\dotfill}}{}}
\nopagebreak\vspace{0.45\baselineskip} {\flushleft\bf
\thesection.\arabic{subsection}~\bf #1.~}
\\*[3mm]\noindent
\nopagebreak}
\newcommand{\sbsect}[1]{
\vspace{0.1cm}\noindent
\textbf{#1.~}\vspace{0.1cm}}

\renewcommand{\subsubsection}{%
\secdef \subsubsect\sbsbsect}
\newcommand{\subsubsect}[2][default]{%
\refstepcounter{subsubsection} 
\addcontentsline{toc}{subsubsection}{{\tocsection{\!\!}
{\hspace{3.05em}\thesubsubsection}{\!\!\!\!#1\dotfill}}{}}
\nopagebreak
\vspace{0.15\baselineskip} \nopagebreak {\flushleft\rmfamily
\itshape\arabic{section}.\arabic{subsection}.\arabic{subsubsection}
\ \rmfamily #1\/.}\ }
\newcommand{\sbsbsect}[1]{\vspace{0.1cm}\noindent
\rmfamily \itshape
\arabic{section}.\arabic{subsection}.\arabic{subsubsection} \
\sffamily #1\/.\ }

\renewcommand{\caption}[1]{%
\vglue0.5cm
\refstepcounter{figure}
\begin{center}
\begin{minipage}[c]{0.8\textwidth}\small {\sc Fig.~\thefigure\ }#1\end{minipage}
\end{center}
}

\newcommand{\supp}{\operatorname{supp}}

\newcommand{\textd}{\text{\rm d}\mkern0.5mu}
\newcommand{\texti}{\text{\rm  i}\mkern0.7mu}
\newcommand{\texte}{\text{\rm  e}\mkern0.7mu}

\newcommand{\Var}{\text{\rm Var}}

\renewcommand{\AA}{\mathcal A}

\newcommand{\LL}{\mathcal L}

\newcommand{\NN}{\mathcal N}

\newcommand{\CalS}{\mathcal S}
\newcommand{\TT}{\mathcal T}

\newcommand{\E}{\mathbb E}

\newcommand{\N}{\mathbb N}

\newcommand{\BbbP}{\mathbb P}

\newcommand{\R}{\mathbb R}

\newcommand{\Z}{\mathbb Z}

\newcommand{\twoeqref}[2]{(\ref{#1}--\ref{#2})}

\newcommand{\cc}{{\text{\rm c}}}

\newcommand{\fraka}{\fraktura a}

\newcommand{\frakb}{\fraktura b}

\newcommand{\frake}{\fraktura e}
\newcommand{\fraks}{\fraktura s}

\def\myffrac#1#2 in #3{\raise 2.6pt\hbox{$#3 #1$}\mkern-1.5mu\raise 0.8pt\hbox{$#3/$}\mkern-1.1mu\lower 1.5pt\hbox{$#3 #2$}}
\newcommand{\ffrac}[2]{\mathchoice%
	{\myffrac{#1}{#2} in \scriptstyle}
	{\myffrac{#1}{#2} in \scriptstyle}
	{\myffrac{#1}{#2} in \scriptscriptstyle}
	{\myffrac{#1}{#2} in \scriptscriptstyle}
}

\newcommand{\wh}{\widehat}
\newcommand{\wt}{\widetilde}
\newcommand{\ol}{\overline}

\newcommand{\laweq}{\,\overset{\text{\rm law}}=\,}

\newcommand{\leb}{{\rm Leb}}

\newcommand{\Lawarrow}{{\,\overset{\text{\rm law}}\longrightarrow\,}}

\newcommand\independent{\protect\mathpalette{\protect\independenT}{\perp}}
\def\independenT#1#2{\mathrel{\rlap{$#1#2$}\mkern3mu{#1#2}}}

\newcommand{\cspecial}{\fraktura c}

\newcommand{\myemph}[1]{\textit{#1}}

\newcommand{\dinfty}{\text{\rm d}_\infty}
\begin{document}

\title[Random walk local time\hfill]{Exceptional points of two-dimensional random\\walks at multiples of the cover time}
\author[\hfill Y.~Abe and M.~Biskup]
{Yoshihiro Abe$^{\,1}$ and Marek~Biskup$^{\,2}$}
\thanks{\hglue-4.5mm\fontsize{9.6}{9.6}\selectfont\copyright\,\textrm{2022}\ \ \textrm{Y. Abe and M.~Biskup.
Reproduction, by any means, of the entire
article for non-commercial purposes is permitted without charge.\vspace{2mm}}}
\maketitle

\vspace{-5mm}
\centerline{\textit{$^1$
Department of Mathematics and Informatics, Chiba University, Chiba, Japan}} 
\centerline{\textit{$^2$
Department of Mathematics, UCLA, Los Angeles, California, USA}}


\vskip0.5cm
\begin{quote}
\footnotesize \textbf{Abstract:}
 We study exceptional sets of the local time of the continuous-time simple
random walk in scaled-up (by~$N$) versions $D_N\subseteq  \mathbb Z^2$ of bounded open domains $D\subseteq  \mathbb R^2$. Upon exit from~$D_N$, the walk  lands on a ``boundary vertex'' and then  reenters~$D_N$  through a random boundary edge in the next step. In the parametrization by the local time at the ``boundary vertex'' we prove that, at times corresponding to a~$\theta$-multiple of the cover time of~$D_N$, the sets of suitably defined $\lambda$-thick (i.e., heavily visited) and $\lambda$-thin (i.e., lightly visited) points are, as $N\to\infty$, distributed according to the Liouville Quantum Gravity~$Z^D_\lambda$ with parameter~$\lambda$-times  the critical value.  For $\theta<1$, also the set of avoided vertices  (a.k.a.\ late points)  and the set where the local time is  of  order unity are distributed according to~$Z^D_{\sqrt\theta}$. The local structure of the exceptional sets is described as well, and is that of a pinned Discrete Gaussian Free Field for the thick and thin points and that of random-interlacement occupation-time field for the avoided points. The results demonstrate universality of the Gaussian Free Field for these extremal problems. 
\end{quote}


\section{Introduction}
\vglue-3mm\subsection{Motivation}
\noindent
In a famous paper from 1960, Erd\H os and Taylor~\cite{ET60} studied the most-frequently visited site by the simple random walk on~$\Z^2$ of time-length~$n$. They showed that the time spent at that site is of order~$(\log n)^2$ and conjectured that the time is asymptotically sharp on that scale. This conjecture was proved in 2001 by Dembo, Peres, Rosen and Zeitouni~\cite{DPRZ01} (see also Rosen~\cite{R05}) who in addition described the multifractal structure of the set of \myemph{thick points}; namely, those points where the local time is at least  a given positive multiple of its maximum. The problem has been revisited numerous times; e.g., by Dembo, Peres, Rosen and Zeitouni~\cite{DPRZ06} who studied random walk late points, by Okada~\cite{O16} who studied the most visited site on the inner boundary of the range, or by Jego~\cite{J18} who extended the results of \cite{DPRZ01,R05} to more general random walks.

Over the past two decades, it has become increasingly clear that many questions about the local time can be usefully rephrased as questions about an associated Discrete Gaussian Free Field (DGFF). This connection,  discovered originally in mathematical physics (Symanzik~\cite{Symanzik}, Brydges, Fr\"ohlich and Spencer~\cite{BFS}), is now elegantly expressed via Dynkin-type Isomorphism/Second Ray-Knight theorems (Dynkin~\cite{D83}, Eisenbaum, Kaspi, Marcus, Rosen and Shi~\cite{EKMRS}). 
 Isomoporphism results of this kind drive the analysis of  many important objects; e.g., random interlacements (Sznitman~\cite{S12}, Rodriguez~\cite{R13}, etc), loop-soups (Lawler and Werner~\cite{LW},  Le Jan~\cite{LeJan}, Lupu~\cite{Lupu}, etc) and the cover time (Ding, Lee and Peres~\cite{DLP12}, Ding~\cite{D14}, etc). 

In the present paper we use the Second Ray-Knight theorem of~\cite{EKMRS} to study  the precise statistics of the thick points for the simple random walk run for times  proportional to  the cover time  of an underlying ``planar'' graph. In addition to the thick points, we analyze also  the sets of \myemph{thin points}, which are those where the local time is less than a fraction of its typical value, \myemph{avoided points}, which are those not visited at all, and \myemph{light points}, where the local time is at most a given constant. We show that all these level sets are intimately connected with the corresponding (so called intermediate) level sets of the Discrete Gaussian Free Field studied earlier by O.~Louidor and the second author~\cite{BL4}. In particular, their limiting statistics is captured by the Liouville Quantum Gravity measures introduced and studied by Duplantier and Sheffield~\cite{DS}.

\subsection{Setting for the random walk}
\label{sec1.2}\noindent
In order to take full advantage of the prior work~\cite{BL4} on the DGFF, we will consider a slightly different setting than the earlier references~\cite{ET60,DPRZ01} and Abe~\cite{A15}, who studied the leading order of the number of thick and thin points for random walk on two-dimensional lattice tori. Indeed, our random walk will behave as the simple random walk only inside a large finite subset of~$\Z^2$; when it exits this set it reenters in the next step through a uniformly-chosen boundary~edge.

To describe the dynamics of our random walk, consider first a general finite, unoriented, connected graph $G = (V\cup\{\varrho\},E)$, where~$\varrho$ is a distinguished vertex (not belonging to~$V$). We assume that each edge~$e\in E$ is endowed with a number~$c_e>0$, called the conductance of~$e$. Let~$X$ denote a continuous-time (constant-speed) Markov chain on~$V\cup\{\varrho\}$  that  makes jumps at independent  rate-1  exponential  random  times to a neighbor selected with the help of transition probabilities
\begin{equation}
\cmss P (u, v) := 
\begin{cases} 
\frac{c_{e}}{\pi(u)}, &~\text{if}~e:=(u, v) \in E, 
\\
0, &~\text{otherwise},
\end{cases}
\end{equation}
where $\pi(u)$ is the sum of~$c_e$ for all edges incident with~$u$. We will use~$P^u$ to denote the law of~$X$ with~$P^u(X_0=u)=1$. 

Given a path~$X$ of the above Markov chain, the local time at~$v\in V\cup\{\varrho\}$ at  time~$t$  is then given by 
\begin{equation} 
\label{E:local_time}
\ell_t^V(v) := \frac1{\pi(u)}\int_0^t\textd s\,1_{\{X_s=u\}},\quad t\ge0,
\end{equation}
where the normalization by $\pi(u)$ ensures that the leading-order growth of~$t\mapsto\ell_t^V(v)$ is the same for all vertices. 
We will henceforth work in the time parametrization by the local time at the distinguished vertex~$\varrho$. For this we set $\hat\tau_{\varrho}(t):=\inf\{s\ge0\colon\ell_s^V(\varrho)  >  t\}$ and denote
\begin{equation}
\label{E:LVt}
L_t^V(v):=\ell_{\hat\tau_{\varrho}(t)}^V(v).
\end{equation}
In this parametrization, $t$ is the  expected (and leading-order) value of~$L^V_t(v)$ under~$P^\varrho$, for every~$v\in V\cup\{\varrho\}$. 

Our derivations will make heavy use of the connection of the above Markov chain with an instance of the Discrete Gaussian Free Field (DGFF). Denoting by
\begin{equation}
H_{v} := \inf \bigl\{t \geq 0 \colon X_t = v \bigr\}
\end{equation}
the first hitting time of vertex~$v$, this DGFF is the centered Gaussian process $\{h_v^{V}\colon v \in V\}$ with covariances given by
\begin{equation}
\label{E:cov}
\E\bigl(h_u^{V} h_v^{V} \bigr) = G^{V} (u, v) :=E^u\bigl(\ell_{H_{\varrho}}^V(v)\bigr).
\end{equation} 
Here and henceforth,~$\E$ denotes expectation with respect to the law~$\BbbP$ of~$h^V$. The field naturally extends to~$\varrho$ by~$h^V_\varrho=0$. 

\begin{figure}[t]
\centerline{\includegraphics[height = 2.0in]{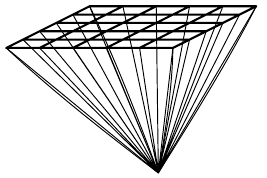}
}
\begin{quote}
\small 
\vglue-0.2cm
\caption{
\label{fig1}
The graph corresponding to~$V$ being the square of $6\times 6$ vertices. Each vertex on the outer perimeter of~$V$ has an edge to the ``boundary vertex''~$\varrho$; the corner vertices that have two edges to~$\varrho$. The ``boundary vertex'' plays the role of the wired boundary condition used often in statistical mechanics. For us this ensures that the associated DGFF vanishes outside~$V$.}
\normalsize
\end{quote}
\end{figure}

\smallskip
Returning back to random walks on~$\Z^2$, in our setting  $V$ stands for  a large finite subset~$V\subseteq  \Z^2$ while~$\varrho$  is  the \myemph{boundary vertex}  obtained by collapsing the set of vertices outside~$V$ to a single point.  The set of edges~$E$ is that between the nearest-neighbor pairs  in~$V$ plus all the edges from~$V$ to~$\Z^2\smallsetminus V$ that now ``end'' in~$\varrho$; see Fig.~\ref{fig1}. The transition rule  of the Markov chain  is that of the simple random walk on the underlying graph; indeed, all conductances take a unit value, $c_e:=1$, at all  the involved  edges including those incident with~$\varrho$. The DGFF associated with this network then corresponds to the ``standard'' DGFF in~$V$ (cf the review by Biskup~\cite{B-notes}) with zero boundary conditions outside~$V$ except that our normalization is slightly different than the one used in~\cite{B-notes} --- indeed, our fields are  half the size of  those in \cite{B-notes}.

For the lattice domains, we will take sequences of subsets of~$\Z^2$ that approximate, in the scaling limit, well-behaved continuum domains. The following definitions are taken from Biskup and Louidor~\cite{BL2}:

\begin{definition} 
An admissible domain 
is a bounded open subset of $\mathbb{R}^2$ 
that consists of a finite number of connected components
and whose boundary is composed of a finite number of connected sets each of which has 
a positive Euclidean diameter.
\end{definition} 

We write $\mathfrak{D}$ to denote the family of all admissible domains and let $d_{\infty} (\cdot, \cdot)$ denote the $\ell^{\infty}$-distance on $\mathbb{R}^2$.

\begin{definition}
An admissible lattice approximation of $D \in \mathfrak{D}$ is a sequence  $\{D_N\}_{N\ge1}$ of subsets of~$\mathbb{Z}^2$  such that the following holds: There is~$N_0\in\N$ such that for all~$N\ge N_0$ we have
\begin{equation}
\label{E:1.8i}
D_N \subseteq   \Bigl\{x \in \mathbb{Z}^2 : 
d_{\infty}\bigl(\ffrac{x}{N}, \mathbb{R}^2 \smallsetminus D\bigr) > \frac{1}{N} \Bigr\}
\end{equation}
and, for any~$\delta>0$ there is~$N_1=N_1(\delta)\in\N$ such that for all~$N\ge N_1$,
\begin{equation}
D_N \supseteq \bigl\{x \in \mathbb{Z}^2 : 
d_{\infty} (\ffrac{x}{N}, \mathbb{R}^2 \smallsetminus D) > \delta \bigr\}.
\end{equation}
\end{definition} 

As shown in \cite{BL2}, these choices ensure that the discrete harmonic measure on~$D_N$ tends, under the scaling of space by~$N$, weakly to the harmonic measure on~$D$. This yields a precise asymptotic expansion of the associated Green functions; see \cite[Chapter~1]{B-notes} for a detailed exposition. In particular, we have $G^{D_N}(x,x)=g\log N+O(1)$ for
\begin{equation}
g:=\frac1{2\pi}
\end{equation}
whenever~$x$ is deep inside~$D_N$. (This is by a factor~$4$ smaller than the corresponding constant in~\cite{B-notes,BL2}  due to a different normalization of our fields.) 

Our random walk will invariably start from the ``boundary vertex''~$\varrho$; throughout we will thus write~$P^\varrho$ for the corresponding law of the Markov chain~$X$. (This law depends on~$N$ but we suppress that notationally.)

\section{Main results}
\noindent
Our aim in this work is to describe the random walk at times that correspond to a $\theta$-multiple of the \myemph{cover time}, for every~$\theta>0$. Recall that the cover time of a graph is the first time that every vertex  of the graph has been visited. Although this is a random quantity, it is quite well concentrated  (provided that the maximal hitting time is of smaller order than the expected cover time; see Aldous~\cite{Aldous}).
 In particular, at the cover time of~$D_N$ the local time at a typical vertex is asymptotic to~$2g(\log N)^2$.  This suggests that we henceforth take~$t$ proportional to $(\log N)^2$ as $N\to\infty$.  

\subsection{Maximum, minimum and exceptional sets}
Let us begin by noting the range of values that the local time takes on~$D_N$:

\begin{theorem}
\label{thm-minmax}
Let~$\{t_N\}_{N\ge1}$ be a positive sequence such that, for some~$\theta>0$,
\begin{equation}
\label{E:1.12}
\lim_{N\to\infty}\frac{t_N}{(\log N)^2}=2g\theta.
\end{equation}
Then for any $D\in\mathfrak D$, any admissible sequence~$\{D_N\}_{N\ge1}$ of lattice approximations of~$D$, the following limits hold in $P^{\varrho}$-probability: 
\begin{equation}
\label{E:max}
\frac1{(\log N)^2}\,\max_{x\in D_N} L^{D_N}_{t_N}(x)\,\,\,\underset{ N\to\infty}\longrightarrow\,\,\,2 g\bigl(\sqrt\theta+1\bigr)^2
\end{equation}
and
\begin{equation}
\label{E:min}
\frac1{(\log N)^2}\,\min_{x\in D_N} L_{t_N}^{D_N}(x)\,\,\,\underset{ N\to\infty}\longrightarrow\,\,\,2 g\bigl[(\sqrt\theta-1)\vee0\,\bigr]^2.
\end{equation}
\end{theorem}

\begin{figure}[t]
\centerline{\includegraphics[width=0.4\textwidth]{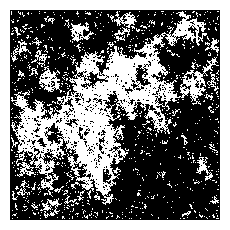}
\includegraphics[width=0.4\textwidth]{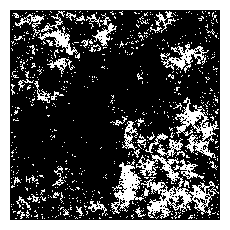}
}
\begin{quote}
\small 
\vglue-0.2cm
\caption{
\label{fig2}
Plots of the $\lambda$-thick (left) and $\lambda$-thin (right) level sets for the same sample of the random walk on a square of side length~$1000$ and parameter choices $\theta:=10$ and $\lambda:=0.1$.}
\normalsize
\end{quote}
\end{figure}

\begin{figure}[t]
\centerline{\includegraphics[width=0.4\textwidth]{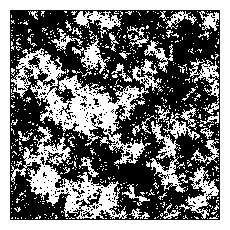}
\includegraphics[width=0.4\textwidth]{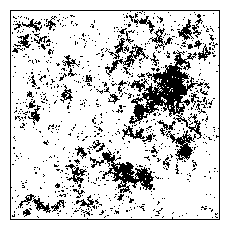}
}
\begin{quote}
\small 
\vglue-0.2cm
\caption{
\label{fig2b}
The sets of avoided points for a sample of the random walk on a square of side-length $N=2000$ observed at times corresponding to $\theta$-multiple of the cover time for $\theta:=0.1$ (left) and $\theta:=0.3$ (right).}
\normalsize
\end{quote}
\end{figure}

These conclusions have previously been obtained by Abe~\cite[Corollary~1.3]{A15} for the  continuous-time  walk on  the~$N\times N$ torus.
As is checked from \eqref{E:min}, the cover time indeed corresponds to $\theta=1$. Noting that the typical value of the local time at a~$\theta$-multiple of the cover time is asymptotic to~$2g\theta(\log N)^2$,  we are naturally led to consider  the set of $\lambda$-\myemph{thick points},
\begin{equation}
\label{E:Lplus}
\TT_N^{+}(\theta,\lambda):=\Bigl\{x\in D_N\colon L_{t_N}^{D_N}(x)\ge 2g(\sqrt\theta+\lambda)^2(\log N)^2\Bigr\},\quad \lambda\in(0,1],
\end{equation}
and $\lambda$-\myemph{thin points},
\begin{equation}
\label{E:Lminus}
\TT_N^{-}(\theta,\lambda):=\Bigl\{x\in D_N\colon L_{t_N}^{D_N}(x)\le 2g(\sqrt\theta-\lambda)^2(\log N)^2\Bigr\},\quad\lambda\in(0,\sqrt\theta\wedge1],
\end{equation}
where the upper bounds on~$\lambda$ reflect on \twoeqref{E:max}{E:min}. 
As a boundary case of $\TT_N^{-}(\theta,\lambda)$, we single out the set of $r$-\myemph{light points},
\begin{equation}
\label{E:2.6i}
\LL_N(\theta,r):=\bigl\{x\in D_N\colon L_{t_N}^{D_N}(x)\le r\bigr\},\quad r\ge0,
\end{equation}
including the special case of the set of \myemph{avoided points},
\begin{equation}
\label{E:2.7i}
\AA_N(\theta):=\bigl\{x\in D_N\colon L_{t_N}^{D_N}(x)=0\bigr\}
\end{equation}
 (Dembo, Peres, Rosen and Zeitouni~\cite{DPRZ06} refer to \eqref{E:2.7i} as the \myemph{late points} but we prefer the above in order to make the distinction between \eqref{E:2.6i} and \eqref{E:2.7i} clear.) 
By \eqref{E:min}, the latter two sets will only be relevant for~$\theta\in(0,1]$.
Our aim is to describe the scaling limit of  all  these sets in the limit as~$N\to\infty$. As shown in Figs.~\ref{fig2} and~\ref{fig2b}, this scaling limit should be a random fractal.

\subsection{Digression on exceptional sets of DGFF}
As noted previously, Biskup and Louidor~\cite{BL4} have addressed similar questions in the context of the DGFF. There the maximum of~$h^{D_N}$ is asymptotic to $2\sqrt g\log N$ and so the set of $\lambda$-thick points is naturally defined as that where the field exceeds $2\lambda\sqrt g\log N$. It was noted that taking a limit of these sets directly does not lead to interesting conclusions as, after scaling space by~$N$, they become increasingly dense in~$D$. A proper way to capture their structure is via the random measure
\begin{equation}
\label{E:etaDGFF}
\eta_N^D:=\frac1{K_N}\sum_{x\in D_N}\delta_{x/N}\otimes\delta_{h^{D_N}_x-a_N},
\end{equation}
where~$\{a_N\}_{N\ge1}$ is a centering sequence with the asymptotic $a_N\sim 2\lambda\sqrt{g}\log N$ and
\begin{equation}
\label{E:1.19e}
K_N:=\frac{N^2}{\sqrt{\log N}}\,\texte^{-\frac{a_N^2}{2g\log N}}.
\end{equation}
In \cite[Theorem~2.1]{BL4} it was then shown that, for each~$\lambda\in(0,1)$ there is~$\cspecial(\lambda)>0$ such that, in the sense of vague convergence of measures on~$\overline D\times(\R\cup\{+\infty\})$,
\begin{equation}
\label{E:1.19}
\eta_N^D\,\,\underset{N\to\infty}\Lawarrow\,\,\cspecial(\lambda)\,Z^D_\lambda(\textd x)\otimes\texte^{-\alpha\lambda h}\textd h,
\end{equation}
where $\alpha:=2/\sqrt g$ and~$Z_\lambda^D$ is a random measure in~$D$ called the Liouville Quantum Gravity (LQG)  at parameter~$\lambda$-times critical.  (While $\eta^D_N$ is defined \myemph{a priori} as a measures on $D\times\R$, we will at times regard it as a measure on $\overline D\times(\R\cup\{+\infty\})$, where $\overline D$ is the closure of~$D$ and the topology on $\R\cup\{+\infty\}$ extends that on~$\R$ so that the intervals of the form $[a,+\infty]$ are compact.) 
The constant~$\cspecial(\lambda)$, given explicitly in terms of~$\lambda$  and the constants in the asymptotic expansion of the potential kernel on~$\Z^2$, allows us take~$Z^D_\lambda$ to be normalized so that, for each Borel set~$A\subseteq   D$,
\begin{equation}
\label{E:1.19a}
\E Z_\lambda^D(A)=\int_A r^{\,D}(x)^{2\lambda^2}\textd x,
\end{equation} 
where~$r^D$ is an explicit function supported on~$D$ that, for~$D$ simply connected, is simply the conformal radius; see~\cite[(2.10)]{BL4}. 

A construction of  the LQG  measures goes back to Kahane's Multiplichative Chaos theory~\cite{Kahane}; they were recently reintroduced  and further studied  by Duplantier and Sheffield~\cite{DS}. Shamov~\cite{Shamov}  neatly  characterized the LQG measures for all~$\lambda\in(0,1)$ by  their expected value  and the behavior under Cameron-Martin shifts of the underlying continuum Gaussian Free Field.

\subsection{Thick and thin points}
\label{sec2.3}\noindent
Inspired by the above developments, we will encode the level sets $\TT_N^{\pm}(\theta,\lambda)$ via the random measures
\begin{equation}
\label{E:zetaND}
\zeta^D_N:=\frac1{W_N}\sum_{x\in D_N}\delta_{x/N}\otimes\delta_{(L_{t_N}^{D_N}(x)-a_N)/\log N}\,,
\end{equation}
where~$\{a_N\}_{N\ge1}$ is a centering sequence and $\{t_N\}_{N\ge1}$ is a sequence of times, both growing proportionally to~$(\log N)^2$, and
\begin{equation}
\label{E:WN}
W_N:=\frac{N^2}{\sqrt{\log N}}\texte^{-\frac{(\sqrt{2t_N}-\sqrt{2a_N})^2}{2g\log N}}.
\end{equation}
The normalization by~$\log N$ in the second delta-mass in \eqref{E:zetaND} indicates that we are tracking variations of the local time of scale~$\log N$.  (As we will see in Section~\ref{sec-2.5},  this is also the order of the variation of the local time between nearest neighbors.) We then get: 

\begin{theorem}[Thick points]
\label{thm-thick}
Suppose $\{t_N\}_{N\ge1}$ and $\{a_N\}_{N\ge1}$ are positive sequences such that, for some $\theta>0$ and some $\lambda\in(0,1)$,
\begin{equation}
\label{E:1.20}
\lim_{N\to\infty}\frac{t_N}{(\log N)^2}=2g\theta\quad\text{\rm and}\quad\lim_{N\to\infty}\frac{a_N}{(\log N)^2}=2g(\sqrt\theta+\lambda)^2.
\end{equation}
For any~$D\in\mathfrak D$, any sequence $\{D_N\}_{N\ge1}$ of admissible approximations of~$D$, and for~$X$ sampled from~$P^{\varrho}$, in the sense of vague convergence of measures on $\overline D\times(\R\cup\{+\infty\})$,
\begin{equation}
\label{E:1.21dis}
\zeta^D_N\,\,\,\underset{N\to\infty}\Lawarrow\,\,\, \frac{\theta^{1/4}}{2\sqrt{g}\,(\sqrt\theta+\lambda)^{3/2}}\,\,\cspecial(\lambda) Z_\lambda^{D}(\textd x)\otimes\texte^{-\alpha(\theta,\lambda) h}\textd h, 
\end{equation}
where $\alpha(\theta,\lambda):=\frac1g\frac{\lambda}{\sqrt\theta+\lambda}$
and $\cspecial(\lambda)$ is as in \eqref{E:1.19}.
\end{theorem}

For the thin points, we similarly obtain:

\begin{theorem}[Thin points]
\label{thm-thin}
Suppose $\{t_N\}_{N\ge1}$ and $\{a_N\}_{N\ge1}$ are positive sequences such that, for some $\theta>0$ and some $\lambda\in(0,1\wedge\sqrt\theta)$,
\begin{equation}
\label{E:1.22}
\lim_{N\to\infty}\frac{t_N}{(\log N)^2}=2g\theta\quad\text{\rm and}\quad\lim_{N\to\infty}\frac{a_N}{(\log N)^2}=2g(\sqrt\theta-\lambda)^2.
\end{equation}
For any~$D\in\mathfrak D$, any sequence $\{D_N\}_{N\ge1}$ of admissible approximations of~$D$, and for~$X$ sampled from~$P^{\varrho}$, in the sense of vague convergence of measures on $\overline D\times(\R\cup\{-\infty\})$, 
\begin{equation}
\label{E:1.23dis}
\zeta^D_N\,\,\,\underset{N\to\infty}\Lawarrow\,\,\, \frac{\theta^{1/4}}{2\sqrt{g}\,(\sqrt\theta-\lambda)^{3/2}}\,\cspecial(\lambda) Z_\lambda^{D}(\textd x)\otimes\texte^{+\tilde\alpha(\theta,\lambda) h}\textd h,
\end{equation}
where~$\tilde\alpha(\theta,\lambda):=\frac1g\frac{\lambda}{\sqrt\theta-\lambda}$ and $\cspecial(\lambda)$ is as in \eqref{E:1.19}.
\end{theorem}

Note that, under \eqref{E:1.20} or \eqref{E:1.22}, the above implies
\begin{equation}
|\TT^\pm_N(\theta,\lambda)|=N^{2(1-\lambda^2)+o(1)},
\end{equation}
where~$o(1)\to0$ in probability.
This conclusion has previously been obtained by the first author in~\cite[Theorem~1.2]{A15}, albeit for random walks on tori and under a different parametrization of the level sets. The present theorems tell us considerably more. Indeed, they imply that points picked at random from $\TT^\pm_N(\theta,\lambda)$ have asymptotically the same statistics as those picked from the set where the DGFF is above the $\lambda$-multiple of its absolute maximum. 

The connection with the DGFF becomes nearly perfect if instead of $\log N$ we normalize the second coordinate of~$\zeta^D_N$ by $\sqrt{2a_N}$. In that parametrization, the resulting measure \myemph{coincides} (up to reversal of the second coordinate for the thin points) with that for the DGFF up to an overall normalization constant. This  demonstrates \myemph{universality} of the Gaussian Free Field in these extremal problems. 

\subsection{Light and avoided points}
The level sets \twoeqref{E:Lplus}{E:Lminus} are naturally nested which  suggests  that, for~$\theta\in(0,1)$, also the sets of $r$-light points~$\LL_N(\theta,r)$ and avoided points $\AA_N(\theta)$  bear a close connection  to an intermediate level set of the DGFF, this time with~$\lambda:=\sqrt\theta$. As the next theorem shows, this is true albeit under a  different normalization:

\begin{theorem}[Light points]
\label{thm-light}
Suppose $\{t_N\}_{N\ge1}$ is a positive sequence such that
\begin{equation}
\label{E:1.29}
\theta:=\frac1{2g}\lim_{N\to\infty}\frac{t_N}{(\log N)^2}\in (0,1).
\end{equation}
For any~$D\in\mathfrak D$, any sequence $\{D_N\}_{N\ge1}$ of admissible approximations of~$D$, and for~$X$ sampled from~$P^{\varrho}$, consider the measure
\begin{equation}
\label{E:varthetaND}
\vartheta^D_N:=\frac1{\wh W_N }\sum_{x\in D_N}\delta_{x/N}\otimes\delta_{L_{t_N}^{D_N}(x)},
\end{equation}
where
\begin{equation}
\label{E:1.31}
\wh W_N :=N^2\texte^{-\frac{t_N}{g\log N}}.
\end{equation}
Then, in the sense of vague convergence of measures on~$\overline D\times[0,\infty)$,
\begin{equation}
\label{E:2.22ii}
\vartheta^D_N\,\,\,\underset{N\to\infty}\Lawarrow\,\,\,  \sqrt{2\pi g}\,\cspecial(\sqrt\theta)\,\,  Z_{\sqrt{\theta}\,}^{D}(\textd x)\otimes\mu(\textd h),
\end{equation}
where $\cspecial(\lambda)$ is as in \eqref{E:1.19} and~$\mu$ is the Borel measure 
\begin{equation}
\label{E:1.33}
\mu(\textd h):=\delta_0(\textd h)+\biggl(\,\sum_{n=0}^\infty\frac1{n!(n+1)!}\Bigl(\frac{\alpha^2\theta}2\Bigr)^{n+1} h^n\biggr)1_{(0,\infty)}(h)\,\textd h.
\end{equation}
\end{theorem}

Note that the density of the continuous part of the measure in \eqref{E:1.33} is uniformly positive on~$[0,\infty)$ and grows exponentially in~$\sqrt h$. Naturally, the atom at zero has the interpretation of the contribution of the avoided points and so we also get:

\begin{theorem}[Avoided points]
\label{thm-avoid}
Suppose $\{t_N\}_{N\ge1}$ is a  positive  sequence such that \eqref{E:1.29} holds. For any $D\in\mathfrak D$, any sequence $\{D_N\}_{N\ge1}$ of admissible approximations of~$D$, and for~$X$ sampled from~$P^{\varrho}$,  consider the measure
\begin{equation}
\label{E:kappaND}
\kappa^D_N:=\frac1{\wh W_N }\sum_{x\in D_N}1_{\{L_{t_N}^{D_N}(x)=0\}}\,\delta_{x/N},
\end{equation}
where~$\wh W_N $ is as in \eqref{E:1.31}. Then, in the sense of vague convergence of measures on~$\overline D$,
\begin{equation}
\label{E:kappa-lim}
\kappa^D_N\,\,\,\underset{N\to\infty}\Lawarrow\,\,\, \sqrt{2\pi g}\,\cspecial(\sqrt\theta)\,\,  Z_{\sqrt{\theta}\,}^{D}(\textd x),
\end{equation}
where $\cspecial(\lambda)$ is again as in \eqref{E:1.19}.
\end{theorem}

We conclude that, at times asymptotic to a $\theta$-multiple of the cover time with~$\theta<1$, the total number of avoided points is~proportional to~$\wh W_N =N^{2(1-\theta)+o(1)}$. Moreover, when normalized by~$\wh W_N $, it tends in law to a constant times the total mass of~$Z^D_{\sqrt\theta}$.

\subsection{Local structure: thick and thin points}
\label{sec-2.5}\noindent
Similarly to the case of the DGFF treated in \cite{BL4}, the convergence of the point measures associated with the exceptional sets can be extended to include information about the local structure of the exceptional sets under consideration. 
For the case of thick and thin points, this structure is captured by the measure on Borel subsets of~$D\times\R\times\R^{\Z^2}$ (under the product topology) defined by
\begin{equation}
\label{E:zeta-local}
\wh \zeta^D_N :=
\frac{1}{W_N} \sum_{x \in D_N} \delta_{x/N} \otimes 
\delta_{(L_{t_N}^{D_N} (x) - a_N)/\log N} 
\otimes \delta_{\{(L_{t_N}^{D_N} (x) - L_{t_N}^{D_N} (x+z))/\log N\colon z \in \mathbb{Z}^2 \}}. 
\end{equation}
In order to express the limit measure, we need to introduce the DGFF $\phi$ on $\Bbb{Z}^2$ pinned to zero at the origin.
This is a centered Gaussian field on $\Bbb{Z}^2$ with law~$\nu^0$  determined by 
\begin{equation}
\label{E:phi-cov}
E_{\nu^0} (\phi_x \phi_y) = \fraka(x) + \fraka(y) - \fraka (x-y),
\end{equation}
where $\fraka \colon \Bbb{Z}^2 \to [0, \infty)$ is the potential kernel,
i.e., the unique function with $\fraka (0) = 0$ which is discrete harmonic 
on $\Bbb{Z}^2 \smallsetminus \{0 \}$ and satisfies $\fraka (x) = g \log |x| + O(1)$
as $|x| \to \infty$. For the thick points, we then get:

\begin{theorem}[Local structure of the thick points]
\label{thm-thick-local}
Under the conditions of Theorem~\ref{thm-thick}
and denoting by $\zeta^D$ the limit measure on the right of (\ref{E:1.21dis}), 
\begin{equation}
\wh \zeta^D_N\,\,\,\underset{N\to\infty}\Lawarrow\,\,\, \zeta^D \otimes \nu_{\theta, \lambda},
\end{equation}
where $\nu_{\theta, \lambda}$ is the law of $2\sqrt{g} (\sqrt{\theta} + \lambda) 
(\phi +  \alpha\lambda \fraka)$ under $\nu^0$.
\end{theorem}

For the thin points, we in turn get:

\begin{theorem}[Local structure of the thin points]
\label{thm-thin-local}
Under the condition of Theorem~\ref{thm-thin} and
denoting by $\zeta^D$ the limit measure on the right of (\ref{E:1.23dis}),
\begin{equation}
\wh \zeta^D_N\,\,\,\underset{N\to\infty}\Lawarrow\,\,\, \zeta^D \otimes \wt \nu_{\theta, \lambda},
\end{equation}
where $\wt \nu_{\theta, \lambda}$ is the law of $2\sqrt{g} (\sqrt{\theta} - \lambda) 
 (\phi -  \alpha \lambda\fraka) $ under $\nu^0$.
\end{theorem} 

As shown in \cite{BL4}, the field $\phi+\lambda\alpha\fraka$ describes the local structure of the DGFF near the points where it takes values (close to) $2\sqrt g\lambda\log N$. As before, the prefactor $2\sqrt g(\sqrt\theta\pm\lambda)$ disappears when instead of~$\log N$ we normalize the third coordinate of~$\wh\zeta^D_N$ by$~\sqrt{2a_N}$. The above results thus extend the universality of the DGFF to the local structure as well. 

\subsection{Local structure: avoided points}
The local structure of the local time near the avoided points will be radically different. Indeed, in the vicinity of an avoided point, the local time will remain of order unity and so a proper way to extend the measure~$\kappa_N^D$ is
\begin{equation}
\wh \kappa_N^D :=
\frac{1}{\wh W_N}
\sum_{x \in D_N}
1_{\{L_{t_N}^{D_N} (x) = 0 \}} \delta_{x/N} \otimes \delta_{\{L_{t_N}^{D_N} (x+z) \colon z \in \Bbb{Z}^2 \}},
\end{equation}
which is now a Borel measure on $D\times[0,\infty)^{\Z^2}$. Moreover, near an avoided point~$x$, the walk itself should behave as if conditioned not to hit~$x$. This naturally suggests that its trajectories will look like two-dimensional \myemph{random interlacements} introduced recently by Comets, Popov and Vachkovskaia~\cite{CPV} and Rodriguez~\cite{R18}, building on earlier work of Sznitman~\cite{S12} and Teixeira~\cite{T09} in transient dimensions. In order to state our limit theorem, we need to review some of the main conclusions from \cite{CPV,R18}.

\begin{figure}[t]
\centerline{\includegraphics[width = 3.2in]{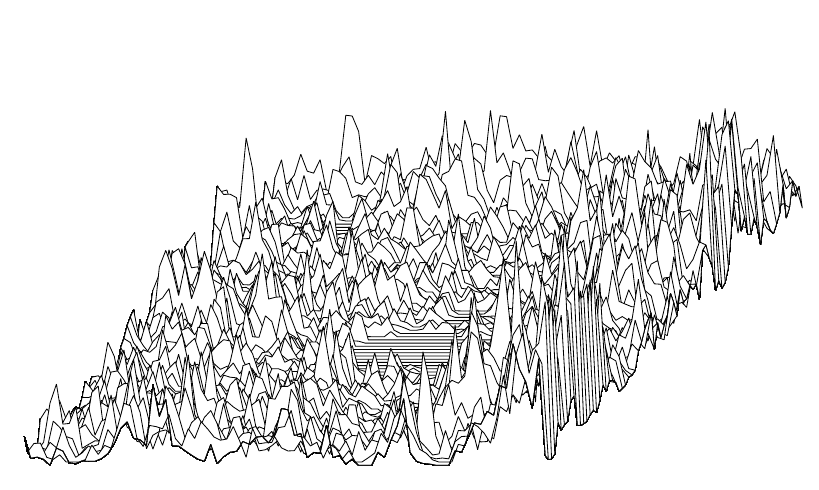}
\hglue-1.7cm\includegraphics[width = 3.2in]{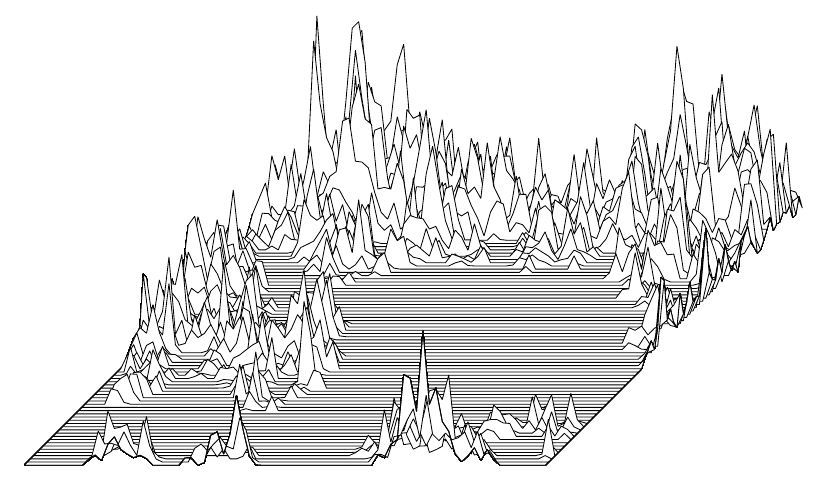}
}
\begin{quote}
\small 
\vglue-0.2cm
\caption{
\label{fig3b}
Samples of the occupation-time field near two randomly-selected avoided points of a random walk run for $0.2$-multiple of the cover time in a square of side-length $N=2000$. Only the square of side-length $81$ centered at the chosen avoided point is depicted.}
\normalsize
\end{quote}
\end{figure}

First we need some notation. Let~$W$ be the set of all doubly-infinite
transient random-walk trajectories on~$\Z^2$; namely, piece-wise constant right-continuous maps $X\colon\R\to\Z^2$ that make only jumps between nearest neighbors and spend only finite time (measured by the Lebesgue measure) in every finite subset of~$\Z^2$. We endow~$W$ with the $\sigma$-field~$\mathcal{W}$ generated by finite-dimensional coordinate projections, $\mathcal{W}:=\sigma(X_t\colon t\in\R)$.
For $A \subseteq   \mathbb{Z}^2$ finite, we write $W_A$ for the subset of~$W$ of the trajectories that visit~$A$.

Next we will put a measure $Q_A^{0, \Z^2}$ on~$W_A$ as follows. Let $\mathfrak h_A$ denote the harmonic measure of~$A$ from infinity (i.e., the distribution of the first entry point to~$A$ by a random walk started at infinity).  Assume~$0\in A$ and let  $\widehat{P}^x$ denote the law of a constant-speed continuous-time random walk on $\mathbb{Z}^2 \smallsetminus \{0\}$ started at~$x$ with conductance $\fraka (y) \fraka (z)$ at nearest-neighbor edges~$(y,z)$ in~$\Z^2$.  By Doob's $h$-transform argument, $\widehat{P}^x$ is the law of the simple random walk on~$\Z^2$ started from~$x$ and conditioned to avoid~$0$. 
For all cylindrical events $E^+, E^- \in \sigma(X_t\colon t\ge0)$ and any  $x \in \Z^2$, we then set
\begin{multline}
\qquad
Q_A^{0, \Z^2} \Bigl((X_{-t})_{t \geq 0} \in E^-,\,X_0 = x,\,(X_t)_{t \geq 0} \in E^+\Bigr)
\\
:= 4\,\fraka (x) \mathfrak h_A (x) \widehat{P}^x(E^+) \widehat{P}^x (E^-\,|\,H_A = \infty).
\qquad
\end{multline}
Note that, since cylindrical events are unable to distinguish left and right path continuity, writing $(X_{-t})_{t \geq 0} \in E^-$ is meaningful. The transience of~$\widehat{P}^x$ implies $\widehat{P}^x(H_A=\infty)>0$  whenever $\mathfrak h_A(x)>0$  and so the conditioning on the right is non-singular. 

The measure $Q_A^{0, \Z^2}$ represents the (un-normalized) law of doubly-infinite trajectories of the simple random walk that hit~$A$ (recall that $\mathfrak{h}_A(x)=0$ unless $x\in A$) but avoid~$0$ for all times. As the main results of \cite{CPV,R18} show, the normalization is chosen such that these measures are consistent, albeit only after factoring out time shifts. To state this precisely, we need some more notation. Regarding two trajectories $w,w'\in W$ as equivalent if they are time shifts of each other --- i.e., if there is $t\in\R$ such that $w(s)=w'(s+t)$ for all~$s\in\R$ --- we use~$W^{\star}$ to denote the quotient space of~$W$ induced by this equivalence relation. Writing $\Pi_\star \colon W \to W^{\star}$ for the canonical projection, the induced $\sigma$-field on~$W^{\star}$ is given by $\mathcal{W}^\star:=\{E\subseteq   W^\star\colon \Pi_\star^{-1}(E)\in \mathcal W\}$. Note that $W_A^{\star} := \Pi_\star (W_A)\in\mathcal{W}^\star$. 

Theorems~3.3 and~4.2 of~\cite{R18} (building on~\cite[Theorem~2.1]{T09}, see also \cite[page 133]{CPV}) then ensure the existence of a (unique) measure on~$W^{\star}$ such that for any finite $A \subseteq   \mathbb{Z}^2$ and any $E\in\mathcal W^\star$,
\begin{equation}
\nu^{0, \Z^2}(E \cap W_A^\star) = Q_A^{0, \Z^2}\circ\Pi_\star^{-1}(E\cap W_A^\star).
\end{equation}
Since $Q_A^{0, \Z^2}$ is a finite measure and the set of finite $A\subseteq  \Z^2$ is countable, $\nu^{0,\Z^2}$ is $\sigma$-finite. We may thus consider a Poisson point process on~$W^\star\times[0,\infty)$ with intensity $\nu^{0, \Z^2}\otimes\leb$. Given a sample~$\omega$ from this process, which we may write as $\omega = \sum_{i\in\N} \delta_{(w_i^{\star}, u_i)}$, and any $u\in[0,\infty)$, we define the \myemph{occupation time field at level~$u$} by
\begin{equation}
\label{E:occ-time-field}
L_u(x) := \sum_{i \in \N} 1_{\{u_i \leq u\}}\,\frac{1}{4} \int_\R \textd t\, 1_{\{w_i (t) = x \}},
\quad x \in \mathbb{Z}^2,
\end{equation}
where $w_i \in W$ is any representative of the class of trajectories marked by~$w_i^\star$; i.e., $\Pi_\star (w_i) = w_i^{\star}$. (The integral does not depend on the choice of the representative.) We are now ready to state the convergence of the measures $\wh \kappa_N^D$.

\begin{theorem}[Local structure of the avoided points]
\label{thm-avoid-local}
Under the conditions of Theorem~\ref{thm-avoid} and for~$\kappa^D$ denoting the measure on the right of \eqref{E:kappa-lim},
\begin{equation}
\wh \kappa^D_N\,\,\,\underset{N\to\infty}\Lawarrow\,\,\, 
\kappa^D \otimes \nu_{\theta}^{\text{\rm RI}},
\end{equation}
where $\nu_{\theta}^{\text{\rm RI}}$ is the law of the occupation time field 
$(L_u(x))_{x \in \Bbb{Z}^2}$ at $u:=\pi\theta$.
\end{theorem}

We expect a similar result to hold for the light points as well but with the random interlacements replaced by a suitably modified version that allows the walks to hit the origin but only accumulating a given (order unity) amount of local time there. Samples of the occupation time field near an avoided point are shown in~Fig.~\ref{fig3b}.

\section{Main ideas, extensions and outline}
\label{sec2}\nopagebreak\noindent
Let us proceed by a brief overview of the main ideas of the proof and then a list of possible extensions and refinements. We also outline the remainder of this paper.

\subsection{Main ideas}
\label{sec2.1}\noindent
As already noted, key for all developments in this paper is the connection of the local time $L^V_t$ and the associated DGFF~$h^V$.
Our initial take on this connection was through the fact that the DGFF represents the fluctuations of~$L^V_t$ at large times via
\begin{equation}
\label{E:GFFlimit}
\frac{L_t^V(\cdot)-t}{\sqrt{2t}}\,\,\,\underset{t\to\infty}\Lawarrow\,\,\, h^V
\end{equation}
 which is proved by decomposing the local time in individual excursions and applying the Central Limit Theorem. 
(The observation \eqref{E:GFFlimit} also guided  the parametrization in the earlier work on this problem, e.g.,~\cite{A15}.) However, as noted at the end of Subsection~\ref{sec2.3}, for the thick and thin points, the effective~$t$ in the correspondence \eqref{E:GFFlimit} of the local time with the DGFF  turns out to be~$a_N$,  rather than~$t_N$,  due to conditioning on large local time. In particular,   approximating the local time  fluctuations  by the DGFF becomes accurate only \myemph{beyond} the times of the order of the cover time. 

We thus base our proofs on a deeper version of the connection, known under the name  \myemph{Second Ray-Knight Theorem} after Ray~\cite{Ray} and Knight~\cite{Knight}  or \myemph{Dynkin isomorphism} after Dynkin~\cite{D83},  although the statement we use is due to Eisenbaum, Kaspi, Marcus, Rosen and Shi~\cite{EKMRS} (with an interesting new proof by Sabot and Tarres~\cite{Sabot-Tarres}): 

\begin{theorem}[Dynkin isomorphism]
\label{thm-Dynkin}
 Consider the random walk on~$V\cup\{\varrho\}$ as detailed in Section~\ref{sec1.2}. 
For each~$t>0$ there exists a coupling of $L^V_t$ (sampled under~$P^\varrho$) and two copies of the DGFF $h^V$ and~$\wt h^V$ such that
\begin{equation}
\label{E:Dynkin1}
h^V\text{\rm\ and } L^V_t \text{\rm\ are independent}
\end{equation}
and
\begin{equation}
\label{E:Dynkin2}
L_t^V(u)+\frac12(h_u^V)^2 = \frac12\bigl(\wt h_u^V+\sqrt{2t}\bigr)^2, \quad u\in V.
\end{equation}
\end{theorem}

\noindent
This is  usually stated  as a distributional identity; the coupling  version is then a result of abstract-nonsense theorems in probability (see Zhai~\cite[Section~5.4]{Zhai}). 

 Our proofs are based on the following natural idea: 
If we could simply disregard the DGFF on the left-hand side  of \eqref{E:Dynkin2},  the relation would tie the level set corresponding to~$L^{D_N}_{t_N}\approx a_N$ to the level sets of the DGFF where
\begin{equation}
\text{either}\quad \wt h^{D_N}\approx\sqrt{2a_N}-\sqrt{2t_N}\quad\text{or}\quad \wt h^{D_N}\approx-\sqrt{2a_N}-\sqrt{2t_N}.
\end{equation}
For~$a_N\to\infty$, the second level set lies further away from the mean of~$\wt h^{D_N}$ than the first and its contribution can therefore be disregarded. (This is true for the thick and thin points; for the light and avoided points both levels play a similar role). One could then simply hope to plug to the existing result \eqref{E:1.19}. 

Unfortunately, since $\Var(h^{D_N}_x)$ is of order~$\log N$, the square of the DGFF on the left of \eqref{E:Dynkin2} is typically of the size of the anticipated fluctuations of~$L^{D_N}_{t_N}$ and so it definitely affects the limiting behavior of the whole quantity. The main technical challenge of the present paper is thus to understand the contribution of this term precisely.  A key observation that makes this possible is that even for~$x\in D_N$ where $L^{D_N}_{t_N}(x)+\frac12(h^{D_N}_x)^2$ takes exceptional values, the DGFF  $h^{D_N}_x$ remains typical (and $L^{D_N}_{t_N}(x)$ is thus dominant). This requires proving fairly sharp single-site tail estimates for the local time and combining them with the corresponding tail bounds for the DGFF. 

Once that is done, we include the field~$h^{D_N}$, properly scaled,  as a third ``coordinate'' of the point process and study weak subsequential limits of these. For instance, for the thick and thin points this concerns the measure
\begin{equation}
\label{E:2.5i}
\frac1{W_N}\sum_{x\in D_N}\delta_{x/N}\otimes\delta_{(L^{D_N}_{t_N}(x)-a_N)/\log N}\otimes\delta_{h^{D_N}_x/\sqrt{\log N}}.
\end{equation}
Here the key is to show that the DGFF part acts, in the limit, as an explicit \myemph{deterministic} measure.  For instance, for the thick and thin points this means that if~$\zeta^D_N$ converges to some~$\zeta^D$ along a subsequence of~$N$'s, the measure in \eqref{E:2.5i} converges to $\zeta^D\otimes\mathfrak g$ where~$\mathfrak g$ is the  normal law~$\NN(0,g)$;  see Lemma~\ref{lemma-add-field}.

Denoting by~$\ell$ the second variable and by~$h$ the third variable in \eqref{E:2.5i}, the Dynkin isomorphism now tells us that the ``law'' of $\ell+\frac{h^2}2$ under \myemph{any} weak subsequential limit of the measures in \eqref{E:2.5i} is the same as the limit ``law'' of the DGFF  centered at  $\sqrt{2a_N}-\sqrt{2t_N}$  (for the thick points) which we know from \eqref{E:1.19}. This produces a convolution-type identity for subsequential limits of the local-time point process. Some technical work then shows that this identity has a unique solution which can be identified explicitly in all cases of interest.

 We note that an important benefit of our reliance on the Dynkin isomorphism is that our arguments  --- and, in particular, the proof of convergence of the measures in \eqref{E:2.5i} --- avoid the need to work with the second moments of the local time. Unlike the first moments, these are harder to control explicitly and that particularly so under additional truncation that would be required to cover the whole regime of interesting behavior.

Our control of the local structure of the exceptional points  also   relies on isomorphism theorems.  For the thick and thin points, we  combine the Dynkin isomorphism with  Theorem~2.1 of~\cite{BL4} that  captures  the local structure of intermediate level sets of the DGFF. For the avoided points, we  instead invoke  the \myemph{Pinned Isomorphism Theorem} of Rodriguez~\cite[Theorem~5.5]{R18} that links the random-interlacement occupation-time field $(L_u(x))_{x\in\Z^2}$ introduced in \eqref{E:occ-time-field} to the pinned DGFF~$\phi$ defined via \eqref{E:phi-cov} as follows:

\begin{theorem}[Pinned Isomorphism Theorem]
\label{thm-pinned-iso}
Let~$u>0$ and suppose that $(L_u(x))_{x\in\Z^2}$ with law $\nu^{\text{\rm RI}}_{u/\pi}$ is independent of~$\{\phi_x\colon x\in\Z^2\}$ with law~$\nu^0$. Then
\begin{equation}
\label{E:3.6uiw}
L_u+\frac12\phi^2\,\laweq\frac12\bigl(\phi+ 2 \sqrt{2u}\,\fraka\bigr)^2,
\end{equation}
where~$\fraka$ is the potential kernel. (The extra factor of~2 compared to \cite[Theorem~5.5]{R18} is due to different normalizations of the local time, the pinned field and the potential kernel.) 
\end{theorem}

It is exactly the generalization of this theorem that blocks us from extending control of the local structure to the light points. Indeed, we expect that, for the light points, the associated process is still that of random interlacements but with the local time at the origin fixed to a given  positive  number. Developing the theory of this process explicitly goes beyond the scope of the present paper.

\subsection{Extensions and refinements}
 We see a number of possible ways  the existing conclusions  may  be refined so let us discuss these in some more detail.

\begin{figure}[t]
\centerline{\includegraphics[height = 2.4in]{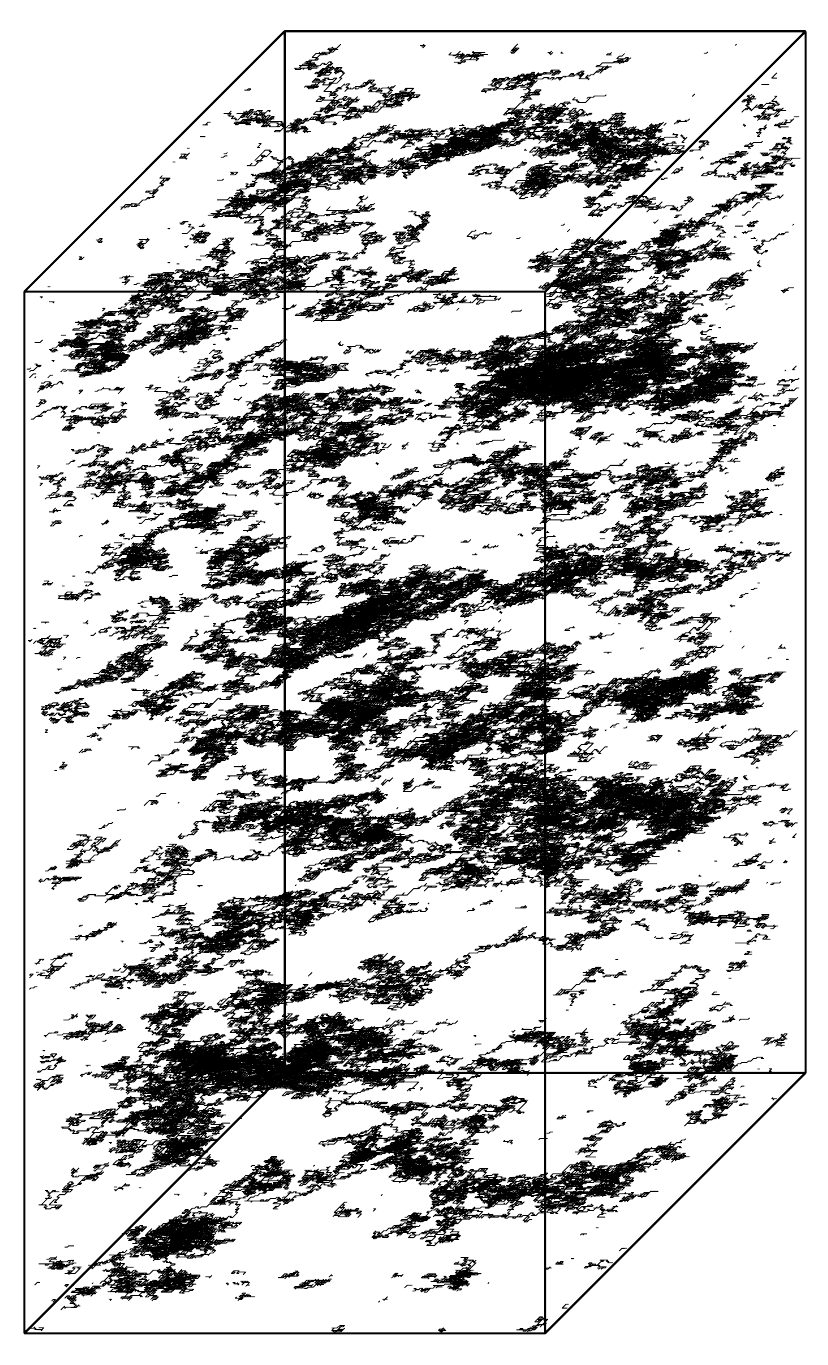}
\includegraphics[width = 2.8in]{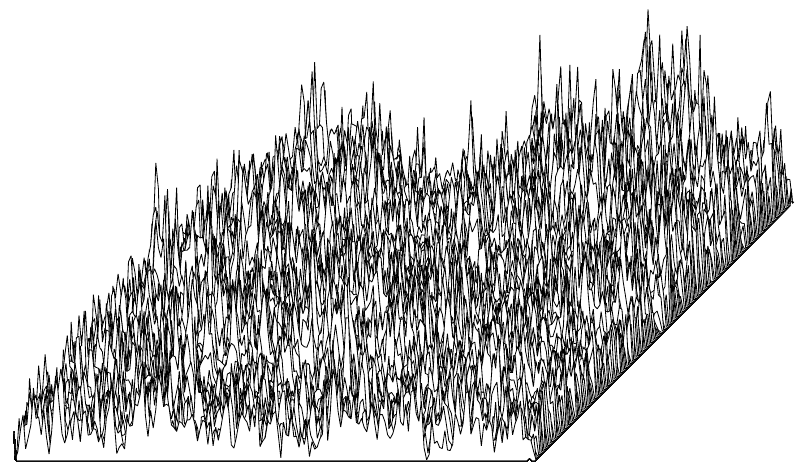}
}
\begin{quote}
\small 
\vglue-0.2cm
\caption{
\label{fig3}
Left: Plot of the trajectory of the random walk on a $200\times200$ square run for $0.3$-multiple of the cover time. The time runs in the vertical direction. Right: The corresponding local time profile. Note that while short excursions near the boundary are numerous, most contribution to the local time profile comes from the excursions that reach ``deep'' into the domain.}
\normalsize
\end{quote}
\end{figure}

\smallskip\noindent
\textsl{Other ``boundary'' conditions: } Perhaps the most significant  deficiency of our setting is the somewhat unnatural mechanism by which the walk returns back to~$D_N$ after each exit.  Contrary to  the intuition one might have,  this does not lead to the local time exploding near the boundary; see Fig.~\ref{fig3} or the fact that~$Z^D_\lambda$ puts no mass on~$\partial D$. The main reason for using the specific setting worked out here is that it allows us to seamlessly plug in the existing results from~\cite{BL4} on the ``intermediate'' level sets of the~DGFF. The natural alternatives are
\settowidth{\leftmargini}{(11111)}
\begin{enumerate}
\item[(1)] running the walk on an $N\times N$ torus, or
\item[(2)]
running the walk as a simple random walk on all of~$\Z^2$ but only  recording  the local time spent inside~$D_N$.
\end{enumerate}
 Both of these require  developing the level-set analysis of a DGFF  on a finite graph pinned at one vertex. 

\smallskip\noindent
\textsl{Time parametrization: }
Another feature for which our setting may  be considered somewhat unnatural  is the parametrization  of the walk  by the time spent at the ``boundary vertex.'' A reasonable question is then what happens when we instead use the parametrization by the actual time of the walk (continuous-time parametrization), or even by the number of discrete steps that the walk has taken (discrete-time parametrization). The main problem here is the lack of a direct connection with the underlying DGFF; instead, one has to rely on approximations.

Preliminary calculations have so far shown that, at least approximately, the local time in the continuous-time parametrization is still connected with the DGFF as in \eqref{E:Dynkin2} but now with the field $\wt h^{D_N}$ reduced by its arithmetic mean over~$D_N$. This implies that, for both continuous and discrete-time analogues of the measures $\zeta^D_N$, $\vartheta^D_N$ and~$\kappa^D_N$, their $N\to\infty$ limits still take the product form as in \eqref{E:1.21dis}, \eqref{E:1.23dis}, \eqref{E:2.22ii} and \eqref{E:kappa-lim}, respectively, albeit now with~$Z^D_\lambda$ replaced by a suitable substitute reflecting on the reduction of the CGFF by its arithmetic mean.  \textit{Update in revision}: These statements have now been established rigorously in Abe, Biskup and Lee~\cite{ABL}.

\smallskip\noindent
\textsl{Critical cases: }
Another natural  extension to consider concerns  various borderline  parameter regimes left out   in the present paper;  namely,  $\lambda:=1$ for the $\lambda$-thick points and~$\lambda:=\sqrt\theta\wedge1$ for the $\lambda$-thin points as well as~$\theta:=1$ for the avoided points. In analogy with the corresponding question for the DGFF (Biskup and Louidor~\cite{BL1,BL2,BL3}), we expect that the corresponding measures  will require a different scaling --- essentially, boosting by an additional factor of~$\log N$ --- and the limit spatial behavior will  be governed by the 
\myemph{critical} LQG measure $Z^D_1$. For the simple random walk on a homogeneous tree of depth~$n$, this program has already been carried out by the first author~(Abe~\cite{A18}). A breakthrough result along these lines describing the limit law of the cover time on homogenous trees has recently been posted by Cortines, Louidor and Saglietti~\cite{CLS}  and by Dembo, Rosen and Zeitouni \cite{DRZ19}.  \textit{Update in revision}: The limit law of the maximum cover time has recently been established by Biskup and Louidor~\cite{BL21}. .

\smallskip\noindent
\textsl{Brownian local time: } 
Yet another  potentially  interesting extension concerns the corresponding problem for the Brownian local time. This requires working with the $\epsilon$-cover time defined as the first time when every disc of radius~$\epsilon>0$ inside~$D$ has been visited; the limit behavior is then studied as~$\epsilon\downarrow0$. We actually expect that, with proper definitions, very similar conclusions will hold here as well although we presently do not see other way to prove them than by approximations via random walks.

Jego~\cite{J18b} recently posted a preprint that proves the existence of a scaling limit for the process associated, similarly to our~$\zeta^D_N$ from \eqref{E:zetaND}, with the local-time thick points of the Brownian path killed upon \myemph{first exit} from~$D$. As it turns out, the limit measure still factors into a product of a random spatial part, defined  via limits of exponentials of the root of the local time, and an exponential measure. However, although the spatial part of the measure obeys the expectation identity of the kind \eqref{E:1.19a}, it is certainly not one of the LQG measures $Z^D_\lambda$ above, due to the limited time horizon of the Brownian path.  In \cite{J19}, Jego characterized the limit measure directly by a list of natural properties.

\subsection{Outline}
The rest of this  paper  is organized as follows. In the next section (Section~\ref{sec3b}) we derive tail estimates for the local time that will come handy later in the proofs. These are used to prove tightness of the corresponding point measures. Section~\ref{sec4} then gives the proof of convergence for the measure associated with  the  $\lambda$-thick points following the outline from Section~\ref{sec2.1}. This proof is then used as a blue print for the corresponding proofs for the $\lambda$-thin points (Section~\ref{sec5}) and the light and avoided points (Section~\ref{sec6}).
The results on  the  local structure are proved at the very end (Section~\ref{sec8}).

\section{Tail estimates and tightness}
\label{sec3b}\noindent\nopagebreak
We are now ready to commence the proofs of our results.  All of our derivations will pertain to the continuous-time Markov chain started, and with the local time parametrized by the time spent, at the ``boundary vertex.'' 
Let us pick a domain~$D\in\mathfrak D$ and a sequence~$\{D_N\}_{N\ge1}$ of admissible approximations of~$D$ and consider these fixed throughout the rest of this paper.  Recall the notation $\zeta^D_N$, $\vartheta^D_N$ and~$\kappa^D_N$ for the measures in \eqref{E:zetaND}, \eqref{E:varthetaND} and \eqref{E:kappaND}, respectively.

\subsection{Upper tails}
We begin  with  estimates on the tails of the random variable~$L^{D_N}_{t_N}(x)$ which then readily imply tightness of the random measures of interest. We first derive these estimates in the general setting of a random walk on a graph with  a  distinguished vertex~$\varrho$ and only then specialize to~$N$-dependent domains in the plane. We begin with the upper tail:

\begin{lemma}[Local time upper tail]
\label{lemma-upper}
Consider the random walk on~$V\cup\{\varrho\}$ as detailed in Section~\ref{sec1.2}.
For all~$a,t>0$ and all~$b\in\R$ such that $a+b>t$, and all~$x\in V$,
\begin{equation}
\label{E:2.2a}
P^\varrho\bigl(L_t^V(x)\ge a+b\bigr)\le \frac{\sqrt{G^V(x,x)}}{\sqrt{2(a+b)}-\sqrt{2t}}\,
\,\texte^{-\frac{(\sqrt{2a}-\sqrt{2t})^2}{2G^V(x,x)}}\,\texte^{-b\frac{\sqrt{2a}-\sqrt{2t}}{G^V(x,x)\sqrt{2a}}}\,.
\end{equation}
\end{lemma}

\begin{proofsect}{Proof}
We will conveniently use estimates developed in earlier work on this problem.  Denoting by $(Y_s)_{s \geq 0}$ the $0$-dimensional Bessel process and writing $P_Y^a$ for its law with $P_Y^a(Y_0=a)=1$, Lemma~3.1(e) of Belius, Rosen and Zeitouni~\cite{BRZ} shows 
\begin{equation}
\label{E:Bessel}
L_t^V(x)~\text{under }P^\varrho 
\,\,\laweq\,\,
\frac{1}{2} \left(Y_{G^V (x, x)} \right)^2~\text{under } P_Y^{\sqrt{2t}}.
\end{equation}
 (Strictly speaking, the derivations in~\cite{BRZ} are restricted to random walks on linear graphs. To make them applicable to our setting, we invoke a ``network reduction'' argument that effectively replaces the underlying graph by a single edge connecting $\varrho$ to~$x$. The reduction preserves both~$G^V(x,x)$ and the law of~$L_t^V(x)$ under~$P^\varrho$.) 

Let $P_B^r$ be a law under which
$(B_s)_{s \geq 0}$ is a standard Brownian motion on $\mathbb{R}$ starting at~$r$. The process~$Y$ is absolutely continuous with respect to~$B$ up to the first time it hits zero; after that~$Y$ vanishes identically. The Radon-Nikodym derivative takes the explicit form (see, for example, \cite[(2.13)]{BRZ})
\begin{equation} 
\label{E:Bessel-to-Brown}
\frac{\textd P_Y^r}{\textd P_B^r}\Bigr|_{\mathcal{F}_{ H_0\wedge t}}
= \sqrt{\frac{r}{B_t}} \exp \left\{- \frac{3}{8} \int_0^t \textd s\,\frac{1}{B_s^2} \right\},\quad\text{on }\{H_0 > t \},
\end{equation}
where $\mathcal{F}_t$ is the $\sigma$-field generated by the process up to time~$t$ and~$H_a$ is the first time the process hits level~$a$.

The identification \eqref{E:Bessel} along with the assumptions~$a+b>0$ translates the event $\{L_t^V(x)\ge a+b\}$ to~$\{Y_t\ge\sqrt{2(a+b)}\}$ intersected by~$\{H_0>t\}$. For~$r:=\sqrt{2t}$, the assumption~$a+b>t$ implies that   the quantity in  \eqref{E:Bessel-to-Brown} is less than one everywhere on the event of interest. Hence,
\begin{equation}
\label{E:4.4}
\begin{aligned}
P^\varrho\bigl(L_t^V(x)\ge a+b\bigr)&\le P_B^{\sqrt{2t}}\Bigl(B_{G^V(x,x)}\ge \sqrt{2(a+b)}\,\Bigr)
\\
&= P_B^0\Bigl(B_{G^V(x,x)}\ge\sqrt{2(a+b)}-\sqrt{2t}\Bigr).
\end{aligned}
\end{equation}
 In order to get \eqref{E:2.2a} from this, we  invoke the Gaussian estimate $P(\NN(0,\sigma^2)\ge x)\le\sigma x^{-1}\texte^{-\frac{x^2}{2\sigma^2}}$ valid for all~$x>0$ along with the calculation
\begin{equation}
\begin{aligned}
\Bigl(\sqrt{2(a+b)}-\sqrt{2t}\Bigr)^2
&=2(a+b)+2t-2\sqrt{2a}\sqrt{2t}\Bigl(1+\frac{b}{a}\Bigr)^{1/2}
\\
&\ge 2(a+b)+2t-2\sqrt{2a}\sqrt{2t}\Bigl(1+\frac{b}{2a}\Bigr)
\\
&=\bigl(\sqrt{2a}-\sqrt{2t}\bigr)^2+2b\frac{\sqrt{2a}-\sqrt{2t}}{\sqrt{2a}},
\end{aligned}
\end{equation}
where in the middle line we used that $(1+x)^{1/2}\le 1+x/2$ holds for all~$x>-1$.
\end{proofsect}

From this we readily obtain:

\begin{corollary}[Tightness for the thick points]
\label{cor-tightness-upper}
Suppose that $t_N$ and~$a_N$ are such that the limits in \eqref{E:1.20} exist for some~$\theta>0$ and some $\lambda\in(0,1)$. For each~$b\in\R$, there is~$c_1(b)\in(0,\infty)$ such that for all~$A\subseteq  \R^2$ closed,
\begin{equation}
\label{E:2.1}
\limsup_{N\to\infty}\,
E^\varrho\bigl[\,\zeta^D_N\bigl(A\times[b,\infty)\bigr)\bigr]\le c_1(b)\,\leb(A\cap D).
\end{equation}
\end{corollary}

\begin{proofsect}{Proof}
It suffices to prove the bound for all $b<0$ with~$|b|$ sufficiently large. 
Pick~$x\in D_N$. If~$G^{D_N}(x,x)\ge \frac g{b^2}\log N$, then Lemma~\ref{lemma-upper} with~$a:=a_N$, $t:=t_N$ and~$b$ replaced by~$b\log N$ and  the uniform  bound~$G^{D_N}(x,x)\le g\log N+c$ give
\begin{equation}
\label{E:2.2}
P^\varrho\bigl(L^{D_N}_{t_N}(x)\ge a_N+b\log N\bigr)
\le \frac{\tilde c}{\sqrt{\log N}}\,\texte^{-\frac{(\sqrt{2a_N}-\sqrt{2t_N})^2}{2g\log N}}\,\texte^{\beta |b|^3},
\end{equation}
for some constants~$\tilde c<\infty$ and~$\beta>0$ independent of~$b$ and~$N$, once~$N$ is sufficiently large. This is  of  order~$W_N/N^2$. If, on the other hand, $G^{D_N}(x,x)\le \frac g{b^2}\log N$, then we use that~$G^{D_N}(x,x)\ge\frac14$ in the second exponential on the right of \eqref{E:2.2a} to get
\begin{equation}
\label{E:2.2b}
P^\varrho\bigl(L^{D_N}_{t_N}(x)\ge a_N+b\log N\bigr)
\le \frac{\tilde c'}{\sqrt{\log N}}\,\texte^{-b^2\frac{(\sqrt{2a_N}-\sqrt{2t_N})^2}{2g\log N}}\,\texte^{\beta' |b|\log N},
\end{equation}
where again~$\tilde c'<\infty$ and~$\beta'>0$ do not depend on~$b$ or~$N$ once~$N$ is sufficiently large. Since the first exponent in \eqref{E:2.2b} is  of  order~$\log N$, for~$|b|$ large enough, this is again at most order~$W_N/N^2$. Now write $A_\epsilon:=\{x\in\R^2\colon \dinfty(x,A)<\epsilon\}$ and note that, in light of \eqref{E:1.8i}, we have 
\begin{equation}
\label{E:3.9ui}
\#\bigl\{x\in D_N\colon x/N\in A\bigr\}\le N^2\leb(A_{1/N}\cap D).
\end{equation}
Summing the relevant bound from \twoeqref{E:2.2}{E:2.2b} over~$x\in D_N$ with~$x/N\in A$, the claim follows by noting that, since~$A$ is closed, $\leb(A_{1/N}\cap D)\to\leb(A\cap D)$ as~$N\to\infty$. 
\end{proofsect}

\subsection{Lower tails}
For the lower tail we similarly get:

\begin{lemma}[Local time lower tail]
\label{lemma-lower}
Consider the random walk on~$V\cup\{\varrho\}$ as in Section~\ref{sec1.2}.
For all~$a,t>0$ and all~$b'<b$ such that $a+b'>0$ and $a+b<t$, and all~$x\in V$,
\begin{equation}
\label{E:4.7a2}
P^\varrho\bigl(L_t^V(x)-a\in[b',b]\bigr)\le \Bigl(\frac{t}{a+b'}\Bigr)^{1/4}\!\!\!
\frac{\sqrt{G^V(x,x)}}{\sqrt{2t}-\sqrt{2(a+b)}}\,
\,\texte^{-\frac{(\sqrt{2t}-\sqrt{2a})^2}{2G^V(x,x)}}\,\texte^{+b\frac{\sqrt{2t}-\sqrt{2a}}{G^V(x,x)\sqrt{2a}}}\,.
\end{equation}
\end{lemma}

\begin{proofsect}{Proof}
We use again the passage \twoeqref{E:Bessel}{E:Bessel-to-Brown} via the Bessel process and Brownian motion except that here we can no longer bound the prefactor in \eqref{E:Bessel-to-Brown} by one. Instead, we get the root of the ratio of the roots of~$2t$ and~$2(a+b')$. Therefore, \eqref{E:4.4} is replaced by
\begin{equation}
P^\varrho\bigl(L_t^V(x)-a\in[b',b]\bigr)\le \Bigl(\frac{t}{a+b'}\Bigr)^{1/4}
P_B^0\Bigl(B_{G^V(x,x)}\le\sqrt{2(a+b)}-\sqrt{2t}\Bigr).
\end{equation}
Noting that the difference in the probability on the right is negative, the rest of the calculation is exactly as before.
\end{proofsect}

 Postponing the tightness of the thin points to the end of this subsection, we first deal with estimates for the light and avoided points: 

\begin{lemma}[Vanishing local time]
\label{lemma-vanish}
For each~$t>0$ and each~$x\in V$,
\begin{equation}
\label{E:L-vanish}
P^\varrho\bigl(L_t^V(x)=0\bigr) = \texte^{-\frac{t}{G^V(x,x)}}.
\end{equation}
In fact, for every~$b\ge0$, we have
\begin{equation}
\label{E:L-vanish2}
P^\varrho\bigl(L_t^V(x)\le b\bigr)  \le\texte^{-\frac{t}{G^V(x,x)}\exp\{-\frac{b}{G^V(x,x)}\}} \le \texte^{-\frac{t}{G^V(x,x)}+b\frac{t}{G^V(x,x)^2}}.
\end{equation}
\end{lemma}

\begin{proofsect}{Proof}
Here we proceed by a direct argument  based on excursion decomposition  (see, however, Remark~\ref{rem-Bessel-works}).  Writing $\hat H_u$ for the first time to return to~$u$ after the walk left~$u$, consider the following independent random variables: 
\settowidth{\leftmargini}{(1111)}
\begin{enumerate}
\item[(1)] $N:=$ Poisson$(t/G^V(x,x))$,
\item[(2)] $\{Z_n\colon n\ge1\}:=$  i.i.d.\  Geometric with parameter~$p:=P^x(H_\varrho<\hat H_x)$,
\item[(3)] $\{T_{k,j}\colon k,j\ge1\}:=$  i.i.d.\  Exponentials with mean one.
\end{enumerate}
We then claim
\begin{equation}
\label{E:3.13}
\pi(x)L_t^V(x) \,\laweq\,\sum_{k=1}^N \sum_{j=1}^{Z_k}T_{k,j}.
\end{equation}
To see this, note that thanks to the parametrization by the local time at~$\varrho$, the value~$L_t^V(x)$ is accumulated through a Poisson$(\pi(\varrho)t)$ number of independent excursions that start and end at~$\varrho$. Each excursion that actually visits~$x$, which happens with probability $P^\varrho(H_x<\hat H_\varrho)$, contributes a Geometric$(p)$-number of independent exponential random variables to the total time the walk spends at~$x$. By Poisson thinning, the number of excursions that visit~$x$ is Poisson with parameter~$\pi(\varrho)P^\varrho(H_x<\hat H_\varrho)t$. We claim that this equals $t/G^V(x,x)$. Indeed, since the walk is constant speed, reversibility gives
\begin{equation}
\pi(\varrho)P^\varrho(H_x<\hat H_\varrho) =\pi(x)P^x(H_\varrho<\hat H_x).
\end{equation}
As was just noted, under~$P^x$ the quantity $\pi(x)\ell_{H_\varrho}(x)$ is the sum of Geometric($p$) independent exponentials of mean one. From \eqref{E:cov} we then get $\pi(x)G^V(x,x)=1/p$.

With \eqref{E:3.13} in hand, to get \eqref{E:L-vanish} we just observe that, modulo null sets, the sum in~\eqref{E:3.13} vanishes only if~$N=0$. To get \eqref{E:L-vanish2} we note that, for~$L_t^V(x)\le b$ we must have $\sum_{j=1}^{Z_k}T_{k,j}\le b\pi(x)$ for each $k=1,\dots,N$. The probability that the sum of~$Z_k$ independent exponentials is less than~$b\pi(x)$ equals $1-\texte^{-bp\pi(x)}$, and that this happens for all $k=1,\dots,N$ thus has probability at most 
\begin{equation}
\sum_{n=0}^\infty \frac{(t/G^V(x,x))^n}{n!}\bigl[ 1-\texte^{-bp\pi(x)}\bigr]^n\texte^{-\frac{t}{G^V(x,x)}}
=\texte^{-\frac{t}{G^V(x,x)} \texte^{-bp\pi(x)}}.
\end{equation}
The claim again follows from $1/p=\pi(x)G^V(x,x)$  and the bound $\texte^{-x}\ge1-x$. 
\end{proofsect}

\begin{remark}
\label{rem-Bessel-works}
We note that a proof based on the connection with the $0$-dimensional Bessel process is also possible. Indeed, by~Belius, Rosen and Zeitouni~\cite[(2.8)]{BRZ}, given $x > 0$ the law of $(Y_s)^2$ under~$P^x_Y$ is given by
\begin{equation}
\texte^{- \frac{x}{2s}} \delta_0 (\textd y) + 1_{(0, \infty)} (y) \,\frac{1}{2s} \,\sqrt{\frac{x}{y}} \,I_1 \Bigl(\frac{\sqrt{xy}}{s}\Bigr)
\texte^{- \frac{x+y}{2s}} \textd y,
\end{equation}
where $I_1 (z) := \sum_{k=0}^{\infty} \frac{(z/2)^{2k+1}}{k! (k+1)!}$.
The identity \eqref{E:L-vanish} follows immediately from \eqref{E:Bessel} and
\begin{equation}
I_1 \Bigl(\frac{\sqrt{2ts}}{G^V (x,x)} \Bigr) \leq \frac{\sqrt{2ts}}{2G^V(x,x)} \,\texte^{\frac{ts}{2 G^V(x,x)^2}}
\end{equation}
then implies the inequality in \eqref{E:L-vanish2} as well.
\end{remark}

From Lemma~\ref{lemma-vanish} we get:

\begin{corollary}[Tightness for the light and avoided points]
\label{cor-light}
Suppose $t_N$ is such that \eqref{E:1.29} holds with some~$\theta\in(0,1)$. 
For each~$b>0$ there is a constant~$c_2(b)\in(0,\infty)$ such that for each $A\subseteq  \R^2$ closed,
\begin{equation}
\label{E:vartheta-lim}
\limsup_{N\to\infty}\,E^\varrho\bigl[\,\vartheta^D_N(A\times[0,b])\bigr]\le  c_2(b)\,\leb(A\cap D).
\end{equation}
In particular,
\begin{equation}
\limsup_{N\to\infty}\,E^\varrho\bigl[\,\kappa^D_N(A)\bigr]\le  c_2(b)\,\leb(A\cap D).
\end{equation}
\end{corollary}

\begin{proofsect}{Proof}
 It suffices to prove just \eqref{E:vartheta-lim} and that for~$b>0$ sufficiently large. Denote~$\tilde c:=\sup_{N\ge1}t_N/(\log N)^{2}$. We then claim
\begin{equation}
\label{E:4.21ie}
P^\varrho\bigl(L^{D_N}_{t_N}(x)\le b \bigr)\le \texte^{-bt_N(\log N)^{-1}}+\texte^{-\frac{t_N}{G^{D_N}(x,x)}+\tilde c b^3\texte^{8b}}.
\end{equation}
Indeed, the first term arises for~$x$ with $G^{D_N}(x,x)\le b^{-1}\texte^{-4b}\log N$ by the first inequality in \eqref{E:L-vanish2} along with~$G^{D_N}(x,x)\ge\frac14$. The second term controls the remaining~$x$; we invoke the second inequality in \eqref{E:L-vanish2} along with $bt_N/G^{D_N}(x,x)^2\le \tilde cb^3\texte^{8b}$. 

For~$b$ sufficiently large, the first term on the right of \eqref{E:4.21ie} is~$o(\wh W_N/N^2)$ independently of~$x\in D_N$. The second term is in turn $O( \wh W_N/N^2)$, with the implicit constant depending on~$b$, by the fact that that~$G^{D_N}(x,x)\le g\log N+c$, uniformly in~$x\in D_N$. The sum over such~$x\in D_N$ with~$x/N\in A$ is now handled via \eqref{E:3.9ui}. 
\end{proofsect}

\subsection{Some corollaries}
Combining the conclusions of Lemmas~\ref{lemma-lower} and~\ref{lemma-vanish}, we can now derive the easier halves of Theorem~\ref{thm-minmax}:

\begin{lemma}
\label{lemma-minmax}
Suppose~$\theta>0$ is related to $t_N$ as in \eqref{E:1.12}. Then for each~$\epsilon>0$, the bounds
\begin{equation}
\label{E:max-upper}
\frac{1}{(\log N)^2}\max_{x\in D_N} L^{D_N}_{t_N}(x)\le 2 g\bigl(\sqrt\theta+1\bigr)^2+\epsilon
\end{equation}
and
\begin{equation}
\label{E:min-lower}
\frac{1}{(\log N)^2}\min_{x\in D_N} L^{D_N}_{t_N}(x)\ge 2 g\bigl[(\sqrt\theta-1)\vee0\bigr]^2-\epsilon
\end{equation}
hold with $P^\varrho$-probability tending to one as~$N\to\infty$.
\end{lemma}

\begin{proofsect}{Proof}
For the maximum, pick~$\epsilon>0$ and abbreviate $a_N:=2 g\bigl(\sqrt\theta+1+\epsilon\bigr)^2(\log N)^2$. Then use \eqref{E:2.2a} with~$b:=0$ and $a:=a_N$ to bound the probability that $L^{D_N}_{t_N}(x)\ge a_N$ by order $N^{-2(1+\epsilon)+o(1)}$ uniformly in~$x\in D_N$. The union bound then gives \eqref{E:max-upper}.

For the minimum, it suffices to deal with the case~$\theta>1$. We pick~$\epsilon>0$ such that $\sqrt\theta>1+\epsilon$. Abbreviate  $a_N:=2 g\bigl(\sqrt\theta-1-\epsilon\bigr)^2(\log N)^2$ and  apply Lemma~\ref{lemma-lower} to get, for any~$b>0$, 
\begin{equation}
\label{E:3.18}
\begin{aligned}
P^\varrho\bigl(L^{D_N}_{t_N}&(x)\le \,a_N\bigr)
\\
&=P^\varrho\bigl(L^{D_N}_{t_N}(x)\le b \bigr)+P^\varrho\bigl(b<L^{D_N}_{t_N}(x)\le a_N\bigr)\\
&\le  P^\varrho\bigl(\,L^{D_N}_{t_N}(x)\le b \bigr) 
+\left(\frac{t_N}{b}\right)^{1/4}\frac{\sqrt{G^{D_N}(x,x)}}{\sqrt{2t_N}-\sqrt{2a_N}}\,\texte^{-\frac{(\sqrt{2t_N}-\sqrt{2a_N})^2}{2G^{D_N}(x,x)}}\,.
\end{aligned}
\end{equation}
 The proof of Corollary~\ref{cor-light} bounds the first probability by $N^{-2\theta+o(1)}$, with~$o(1)\to0$ uniformly in~$x\in D_N$. (As the quantity is non-decreasing in~$b$, the requirement that~$b$ be sufficiently large is achieved trivially.) Hence, even after summing over~$x\in D_N$, the contribution of this term is negligible.  
 
For the second term  on the right of \eqref{E:3.18} we note that, invoking the uniform upper bound $G^{D_N}(x,x)\le g\log N+c$, the above choice of $a_N$ yields 
\begin{equation}
\frac{(\sqrt{2t_N}-\sqrt{2a_N})^2}{2G^{D_N}(x,x)} \ge 2\bigl(1+\epsilon+o(1)\bigr)^{2}\log N
\end{equation}
 uniformly in~$x\in D_N$. 
As the prefactors  produce  only polylogarithmic terms in~$N$, also the second term on the right of \eqref{E:3.18} is~$o(N^{-2})$ as~$N\to\infty$. 
\end{proofsect}

A similar argument will allow us to deal with the tightness of the thin points:

\begin{corollary}[Tightness for the thin points]
\label{cor-tightness-lower}
Suppose that $t_N$ and~$a_N$ are such that the limits in \eqref{E:1.22} exist for some~$\theta>0$ and some $\lambda\in(0,\sqrt\theta\wedge1)$. For all~$b\in\R$ there is~$c_3(b)\in(0,\infty)$ such that for all~$A\subseteq  \R^2$ closed,
\begin{equation}
\limsup_{N\to\infty}\,
E^\varrho\bigl[\,\zeta^D_N\bigl(A\times(-\infty,b]\bigr)\bigr]\le c_3(b)\,\leb(A\cap D).
\end{equation}
\end{corollary}

\begin{proofsect}{Proof}
We proceed as in the proof of Lemma~\ref{lemma-minmax}. Let~$a_N\sim 2g(\sqrt\theta-\lambda)^2(\log N)^2$ be as given, pick $\epsilon\in(0,\sqrt\theta-\lambda)$ an abbreviate~$\hat a_N:= 2g\epsilon^2(\log N)^2$. Then for any $b' > 0$,
\begin{multline}
\label{E:3.18ie}
\quad P^\varrho\bigl(L^{D_N}_{t_N}(x)\le \,a_N + b\log N\bigr)
=P^\varrho\bigl(\,L^{D_N}_{t_N}(x)\le b' \bigr)
\\+P^\varrho\bigl(b'<L^{D_N}_{t_N}(x)\le \hat a_N\bigr)
+P^\varrho\bigl(\hat a_N<L^{D_N}_{t_N}(x)\le a_N + b\log N\bigr).
\quad
\end{multline}
Exactly as in \eqref{E:3.18}, the first term on the right is estimated to be $N^{-2\theta+o(1)}=o(W_N/N^2)$ uniformly in $x\in D_N$, where we used that $W_N=N^{2-2\lambda^2+o(1)}$ and $\lambda<\sqrt\theta$. The second term is bounded as in \eqref{E:3.18} by $N^{-2(\sqrt\theta-\epsilon)^2+o(1)}=o(W_N/N^2)$ by our choice of~$\epsilon$. 

 For the last term we invoke Lemma~\ref{lemma-lower} with $a+b'$ and $a+b$ set to $\hat a_N$ and $a_N+b\log N$, respectively. This allows for~$b$ in \eqref{E:4.7a2} to be negative which permits bounding the last factor on the right by one while keeping the prefactors in \eqref{E:4.7a2} bounded by a constant that depends only on~$\epsilon$, uniformly in~$x\in D_N$. Hence, the last term in \eqref{E:3.18ie} is  $O(W_N/N^2)$ uniformly in~$x\in D_N$. The observation \eqref{E:3.9ui} then helps us  deal with the sum over~$x\in D_N$ subject to~$x/N\in A$.
\end{proofsect}

\begin{remark}
\label{rem-restrict}
The reason for using the expressions $\leb(A\cap D)$ to control the first moments of the measures of interest is that this will later allow us to restrict attention to~$A\subseteq   D$ open with~$\overline A\subseteq   D$ in the arguments to follow. Indeed, taking~$\{A_n\}_{n\ge1}$ open with~$A_n\uparrow D$, as $n\to\infty$ the expected measure of the complement $D\smallsetminus A_n$ tends to zero by the fact that $\leb(D\smallsetminus A_n)\to0$.
\end{remark}

\section{Thick points}
\label{sec4}\noindent\nopagebreak
We are now ready to move to the proof of the stated convergence for the point measure associated with~$\lambda$-thick points. Throughout we will assume that~$a_N$ and~$t_N$ satisfy \eqref{E:1.20} with some~$\theta>0$ and some~$\lambda\in(0,1)$. Introduce the auxiliary centering sequence
\begin{equation}
\label{E:4.1}
\wh a_N:=\sqrt{2a_N}-\sqrt{2t_N}
\end{equation}
and note that $\wh a_N\sim 2\lambda\sqrt g\log N$ as~$N\to\infty$.  The arguments below make frequent use of the coupling of~$L_{t_N}^{D_N}$ and an independent DGFF $h^{D_N}$ to another DGFF~$\wt h^{D_N}$ via the Dynkin isomorphism (Theorem~\ref{thm-Dynkin}). We will use these notations throughout and write~$\wh\eta_N^D$ to denote the DGFF process associated with~$\wt h^{D_N}$ and the centering sequence~$\wh a_N$. A key point to note is that~$W_N$ then coincides with normalizing constant from \eqref{E:1.19e}. 

\subsection{Tightness considerations}
The proof of Theorem~\ref{thm-thick} naturally divides into two parts. In the first part we dominate~$\zeta^D_N$ using~$\wh\eta_N^D$ and control the effect of adding~$h^{D_N}$ to the local time~$L^{D_N}_{t_N}$. The second part is then a derivation, and a solution, of a convolution-type identity linking the weak-limits of~$\zeta^D_N$ to those of~$\wh\eta_N^D$. Our tightness considerations start by the following domination lemma: 

\begin{lemma}[Domination by DGFF process]
\label{lemma-domination}
For any~$b\in\R$ and any measurable~$A\subseteq   D$,
\begin{equation}
\label{E:4.7a}
\zeta^D_N\bigl(A\times[b,\infty)\bigr)\,\,\overset{\text{\rm law}}\le\,\,\wh\eta_N^D\Bigl(A\times\bigl[\tfrac1{2\sqrt g}\tfrac{b}{\sqrt\theta+\lambda},\infty\bigr)\Bigr)
+o(1)
\end{equation}
where~$o(1)\to0$ in probability as $N\to\infty$. Similarly, for any measurable~$A\subseteq   D\times D$ and any~$b\in\R$,
\begin{equation}
\label{E:4.7b}
\zeta^D_N\otimes\zeta^D_N\bigl(A\times[b,\infty)^2\bigr)\,\,\overset{\text{\rm law}}\le\,\,\wh\eta_N^D\otimes\wh\eta_N^D\Bigl(A\times\bigl[\tfrac1{2\sqrt g}\tfrac{b}{\sqrt\theta+\lambda},\infty\bigr)^2\Bigr)
+o(1).
\end{equation}
\end{lemma}

\begin{proofsect}{Proof}
Let us start by \eqref{E:4.7a}. The Dynkin isomorphism shows
\begin{equation}
\label{E:4.4u}
L^{D_N}_{t_N}\le L^{D_N}_{t_N}+\frac12 (h^{D_N})^2= \frac12\bigl(\wt h^{D_N}+\sqrt{2t_N}\bigr)^2.
\end{equation}
 For expression on the left of \eqref{E:4.7a} we then get 
\begin{equation}
\label{E:4.10}
\zeta^D_N\bigl(A\times[b,\infty)\bigr)
\le\frac1{W_N}\sum_{\begin{subarray}{c}
x\in D_N\\x/N\in A
\end{subarray}}
1_{\{|\wt h^{D_N}_x+\sqrt{2t_N}|\ge\sqrt{2a_N+2b\log N}\}}.
\end{equation}
Pick any~$b'<b\tfrac1{2\sqrt g}\tfrac{1}{\sqrt\theta+\lambda}$.  Once~$N$ is sufficiently large, the asymptotic formulas for~$a_N$ and~$t_N$  give $\sqrt{2a_N+2b\log N}\ge \sqrt{2a_N}+b'$ and so 
\begin{equation}
\label{E:4.6e}
1_{\{|\wt h^{D_N}_x+\sqrt{2t_N}|\ge\sqrt{2a_N+2b\log N}\}}
\le 1_{\{\wt h^{D_N}_x\ge\wh a_N+b'\}}+
1_{\{\wt h^{D_N}_x\le-\sqrt{2a_N}-\sqrt{2t_N}- b'\}}\,.
\end{equation}
 Writing~$\overline\eta^D_N$ for the process associated with the field $-\wt h^{D_N}$ and the centering sequence $\sqrt{2a_N}+\sqrt{2t_N}$, and~$\bar K_N$ for the associated normalization from \eqref{E:1.19e}, we thus have
\begin{equation}
\zeta^D_N\bigl(A\times[b,\infty)\bigr)
\le \wh\eta^D_N\bigl(A\times[b',\infty)\bigr)+\frac{\bar K_N}{W_N}\,\overline\eta^D_N\bigl(A\times[b',\infty)\bigr).
\end{equation}
Noting that $\{\overline\eta^D_N\colon N\ge1\}$ is tight on~$\overline D\times(\R\cup\{+\infty\})$ and~$\bar K_N=o(W_N)$, the second term is~$o(1)$ in probability as~$N\to\infty$. To get \eqref{E:4.7a} we now take~$b'$ to~$b\tfrac1{2\sqrt g}\tfrac{1}{\sqrt\theta+\lambda}$ and invoke the continuity of the limit measure in \eqref{E:1.19} in the second variable.

The proof of \eqref{E:4.7b} is completely analogous. Indeed, the same reasoning implies, for any Borel~$A\subseteq   D\times D$,
\begin{multline}
\quad\zeta^D_N\otimes\zeta^D_N\bigl(A\times[b,\infty)^2\bigr)
\le \wh\eta^D_N\otimes \wh\eta^D_N\bigl(A\times[b',\infty)^2\bigr)+
\frac{\bar K_N}{W_N}\,\wh\eta^D_N\otimes\overline\eta^D_N\bigl(A\times[b',\infty)^2\bigr)
\\
+\frac{\bar K_N}{W_N}\,\overline\eta^D_N\otimes\wh\eta^D_N\bigl(A\times[b',\infty)^2\bigr)
+\Bigl(\frac{\bar K_N}{W_N}\Bigr)^2\,\overline\eta^D_N\otimes\overline\eta^D_N\bigl(A\times[b',\infty)^2\bigr).
\quad
\end{multline}
Replacing~$A$ by $D\times D$ in the last three terms shows, via~$\bar K_N=o(W_N)$, that these three terms are again all~$o(1)$ in probability as~$N\to\infty$. A continuity argument in the second variable then proves \eqref{E:4.7b} as well.
\end{proofsect}

Note that Lemma~\ref{lemma-domination} provides an independent proof of the tightness of the measures~$\zeta^D_N$. Based on the proof one might think that~$\zeta^D_N$ is asymptotically close to $\wh\eta_N^D$, but this is false: Although \eqref{E:4.6e} is asymptotically sharp, the inequalities in \twoeqref{E:4.4u}{E:4.10} are not. To account for this  fact,  we have to carefully examine the effect of adding the half of the DGFF-squared to the local time. In particular, we have to ensure that the DGFF remains typical even at the points where the local time \myemph{combined with} half of its square is large. This rather important step is the content of:

\begin{lemma}
\label{lemma-no-conspire}
Let $0<\beta<\frac1{2g}\frac{\sqrt\theta}{\sqrt\theta+\lambda}$. Then for each~$b\in\R$ there is~$
 c_4(b)\in(0,\infty)$ such that for all $M\ge0$, all sufficiently large~$N$ and all~$x\in D_N$,
\begin{equation}
\label{E:4.9uiu}
P^\varrho\otimes\BbbP\biggl(L^{D_N}_{t_N}(x)+\frac{(h^{D_N}_x)^2}2\ge a_N+b\log N,\,\frac{|h^{D_N}_x|}{\sqrt{\log N}}\ge M\biggr)\le  c_4(b)\frac{W_N}{N^2}\texte^{-\beta M^2}.
\end{equation}
\end{lemma}

\begin{proofsect}{Proof}
 Since the $b\log N$-correction can be absorbed into a re-definition of~$a_N$, which thanks to the assumed asymptotic behavior of~$a_N$ and~$t_N$ only changes~$W_N$ by a multiplicative constant, we may assume  for simplicity  that~$b=0$.  Assume also that~$M$ is an integer and pick~$\delta$ with
\begin{equation}
\label{E:4.2a}
0<\delta<2\sqrt\theta\lambda.
\end{equation}
Partitioning the event  in \eqref{E:4.9uiu}  according to which interval of the form $[n,n+1)$, with~$n\in\N$ subject to $n\ge M^2$, the ratio $(h^{D_N}_x)^2/\log N$ lies in, the probability  in \eqref{E:4.9uiu} is bounded~by  
\begin{multline}
\label{E:4.2}
\BbbP\Bigl(|h^{D_N}_x|\ge 2\sqrt g\sqrt{\lambda^2+\delta}\,\log N\Bigr)
\\
+\sum_{M^2\le n\le4 g(\lambda^2+\delta)\log N}\!\!\!\!\BbbP\bigl((h^{D_N}_x)^2\ge n\log N\bigr)P^\varrho\Bigl(L^{D_N}_{t_N}(x)\ge a_N-\frac12(n+1)\log N\Bigr).
\end{multline}
A standard Gaussian bound estimates the first probability by a constant times 
$N^{- 2(\lambda^2+\delta)}$ which is $o(W_N/N^2)$ as~$N\to\infty$. Concerning the terms in the sum, here we first note that for all~$n$ under the summation symbol,
\begin{equation}
\begin{aligned}
a_N-\frac12(n+1)\log N
&\ge 2g\Bigl[(\sqrt\theta+\lambda)^2-(\lambda^2+\delta)+o(1)\Bigr](\log N)^2
\\
&=t_N+2g\bigl(2\sqrt\theta\lambda-\delta+o(1)\bigr)(\log N)^2.
\end{aligned}
\end{equation}
Hence, under \eqref{E:4.2a}, Lemma~\ref{lemma-upper} can be applied. Using $G^{D_N}(x,x)\le g\log N+c$, the term corresponding to integer~$n$ in the sum is thus bounded by
\begin{equation}
\label{E:4.4ui}
\tilde c\frac{W_N}{N^2}\exp\biggl\{\frac1{2g}\Bigl[(n+1)\frac{\sqrt{2a_N}-\sqrt{2t_N}}{\sqrt{2a_N}}-n\Bigr]\biggr\},
\end{equation}
where~$\tilde c$ is a constant that depends on~$\theta$, $\lambda$ and our choice of~$\delta$ but not on~$N$ or~$x$ or~$n$.
Since the assumptions on~$a_N$ and~$t_N$ give
\begin{equation}
\frac{\sqrt{2a_N}-\sqrt{2t_N}}{\sqrt{2a_N}}\,\underset{N\to\infty}\longrightarrow\,\frac{\lambda}{\sqrt\theta+\lambda}<1-2g\beta
\end{equation}
as soon as~$N$ is sufficiently large, the quantity in \eqref{E:4.4ui} is summable on~$n$ and the sum in \eqref{E:4.2} is thus dominated by the term with~$n=M^2$. The claim follows.
\end{proofsect}

\subsection{Convolution identity}
We now move to the second part which consists of  the  derivation of, and a solution to, a convolution identity that links weak (subsequential) limits of $\zeta^D_N$ to those of~$\wh\eta^D_N$.  A key input here is the observation  that, at the scale of its typical fluctuations, the field~$h^{D_N}$ that we add to~$L^{D_N}_{t_N}$ in the Dynkin isomorphism acts like white noise:

\begin{lemma}
\label{lemma-add-field}
Suppose $\{N_k\}$ is a subsequence along which~$\zeta^D_N$ converges in law to~$\zeta^D$. Then
\begin{equation}
\label{E:4.7}
\frac1{W_N}\sum_{x\in D_N}\delta_{x/N}\otimes\delta_{(L^{D_N}_{t_N}(x)-a_N)/\log N}\otimes\delta_{h_x^{D_N}/\sqrt{\log N}}\,\,\underset{\begin{subarray}{c}
N=N_k\\k\to\infty
\end{subarray}}{\overset{\text{\rm law}}\longrightarrow}\,\,
\zeta^D\otimes\mathfrak g,
\end{equation}
where~$\mathfrak g$ is the law of~$\NN(0,g)$.
\end{lemma}

\begin{proofsect}{Proof}
Denote by~$\zeta^{D,\text{ext}}_N$ the measure on the left of \eqref{E:4.7}.  We need to show that the integral of any~$f\in C_\cc(D\times\R\times\R)$ with respect to~$\zeta^{D,\text{ext}}_N$ converges in law to that with respect to~$\zeta^D\otimes\mathfrak g$. The restrictions on~$f$ imply that there is a compact set~$A\subseteq   D$ and a number~$b>0$  such that
\begin{equation}
\label{E:4.8}
\bigl|f(x,\ell,h) \bigr|\le \Vert f\Vert_\infty 1_A(x)1_{[-b,\infty)}(\ell)1_{[-b,b]}(h).
\end{equation}
The argument is based on a \myemph{conditional} second moment calculation and domination by the DGFF process from Lemma~\ref{lemma-domination}.

Abbreviate $L(x):=(L^{D_N}_{t_N}(x)-a_N)/\log N$ and $h_x:=h^{D_N}_x/\sqrt{\log N}$.  Writing $\text{Var}_{\BbbP}$, resp., $\text{Cov}_{\BbbP}$ for the conditional variance, resp., covariance  given the local time, we have
\begin{equation}
\label{E:4.11}
\text{Var}_{\BbbP}\bigl(\langle\zeta^{D,\text{ext}}_N,f\rangle\bigr)
=\frac1{W_N^2}\sum_{x,y\in D_N}\text{Cov}_{\BbbP}\Bigl(f\bigl(\ffrac xN, L(x),h_x\bigr),f\bigl(\ffrac yN, L(y),h_y\bigr)\Bigr).
\end{equation}
Pick~$\epsilon>0$ and split the sum according to whether~$|x-y|\ge\epsilon N$ or not.  Focusing first on the former case,  we use the Gibbs-Markov decomposition to write~$h^{D_N}$ using the value~$h^{D_N}_x$ and an independent DGFF in~$D_N\smallsetminus\{x\}$ as
\begin{equation}
\label{E:5.12i}
h^{D_N}\laweq h^{D_N}_x\frakb_{D_N,x}(\cdot)+\wh h^{D_N\smallsetminus\{x\}},\quad h^{D_N}_x\independent \wh h^{D_N\smallsetminus\{x\}},
\end{equation}
where~$\frakb_{D_N,x}\colon\Z^2\to[0,1]$ is the unique function that is discrete harmonic on~$D_N\smallsetminus\{x\}$, vanishes outside~$D_N$ and equals one at~$x$. A key point, proved with the help of monotonicity of~$D\mapsto\frakb_{D,x}(y)$ with respect to  the  set inclusion, is
\begin{equation}
\label{E:5.13i}
\max_{\begin{subarray}{c}
x,y\in D_N\\|x-y|\ge\epsilon N
\end{subarray}}
\frakb_{D_N,x}(y)\le\frac{c(\epsilon)}{\log N},
\end{equation}
where~$c(\epsilon)\in(0,\infty)$ is independent of~$N$. 

 Write~$R_f(\delta)$ is the maximal oscillation of~$f$ in the third variable on intervals of size~$\delta$. In light of \eqref{E:4.8} we then get
\begin{multline}
\label{E:5.20ii}
\quad
E_\BbbP\Bigl(f(\dots,h_x)f(\dots,h_y)\Bigr)
\\
\le \Vert f\Vert_\infty R_f\Bigl(\frac{b c(\epsilon)}{\sqrt{\log N}}\Bigr)+E_\BbbP\bigl(f(\dots,h_x)\bigr)E_\BbbP\bigl(f(\dots,\wh h_y)\bigr),
\quad
\end{multline}
where~$\wh h$ abbreviates the field $\wh h^{D_N\smallsetminus\{x\}}$ and the dots stand for the remaining arguments of~$f$ that are not affected by the expectation with respect to~$\BbbP$.
As to the expectation on the right, for any~$M>b$ we similarly obtain
\begin{equation}
\label{E:5.21ii}
\Bigl|E_\BbbP\bigl(f(\dots,\wh h_y)\bigr)-E_\BbbP\bigl(f(\dots,h_y)\bigr)\Bigr|
\le R_f\Bigl(\frac{M c(\epsilon)}{\sqrt{\log N}}\Bigr)+2\texte^{-\tilde cM^2}\Vert f\Vert_\infty
\end{equation}
by splitting the expectations depending on the containment in $\{|h_x|\le M\sqrt{\log N}\}$ or not and estimating each term separately. The (positive) constant~$\tilde c$ can be taken as close to~$(2g)^{-1}$ as desired by taking~$M$ sufficiently large. 

Putting \twoeqref{E:5.20ii}{E:5.21ii} together and invoking \eqref{E:4.8}, the contribution of the pairs $(x,y)$ with $|x-y|\ge\epsilon N$ to  \eqref{E:4.11} is thus at most 
\begin{equation}
\label{E:5.20uio}
2\Vert f\Vert_\infty \,\biggl(R_f\Bigl(\frac{M c(\epsilon)}{\sqrt{\log N}}\Bigr)+\texte^{-\tilde cM^2}\Vert f\Vert_\infty\biggr)\,\zeta^D_N\otimes\zeta^D_N\bigl(D\times D\times[-b,\infty)^2\bigr),
\end{equation}
Writing the product-measure term on the right of \eqref{E:5.20uio} as the square of $\zeta^D_N(D\times[-b,\infty))$ we note that this term is stochastically bounded in the limit as~$N\to\infty$ by Corollary~\ref{cor-tightness-upper} (or by the domination argument from Lemma~\ref{lemma-domination}). Since $R_f(\delta)\to0$ as~$\delta\downarrow0$ by the uniform continuity of~$f$, taking~$N\to\infty$ followed by~$M\to\infty$ shows that the sum \eqref{E:4.11} restricted to $|x-y|\ge\epsilon N$ vanishes in $P^\varrho$-probability as~$N\to\infty$ for every~$\epsilon>0$.

Moving to the part of the sum  in \eqref{E:4.11} corresponding to~$|x-y|\le\epsilon N$, using  \eqref{E:4.8} this is bounded by $\Vert f\Vert_\infty^2$ times
\begin{equation}
\label{E:4.12}
\zeta^D_N\otimes\zeta^D_N\Bigl(\bigl\{(x,y)\colon |x-y|\le\epsilon\bigr\}\times[-b,\infty)^2\Bigr)
\end{equation}
which by Lemma~\ref{lemma-domination} is stochastically bounded by
\begin{equation}
\wh\eta_N^D\otimes\wh\eta_N^D\Bigl(\bigl\{(x,y)\colon |x-y|\le\epsilon\bigr\}\times
[-\tfrac{1}{2\sqrt{g}} \tfrac{b}{\sqrt{\theta}+\lambda},\infty)^2\Bigr)+o(1).
\end{equation}
As $\{(x,y)\colon |x-y|\le\epsilon\}$ is closed and $\wh a_N\sim2\lambda\sqrt g\log N$ as~$N\to\infty$, \eqref{E:1.19} and the Portmanteau Theorem show that this expression is, in the limit~$N\to\infty$, stochastically dominated by a  $b$-dependent constant times
\begin{equation}
\label{E:4.22uo}
Z^D_\lambda\otimes Z^D_\lambda\bigl(\{(x,y)\colon |x-y|\le\epsilon\}\bigr).
\end{equation}
This tends to zero as $\epsilon\downarrow0$ a.s.\ due to the fact that~$Z^D_\lambda$ has no point masses a.s.

We conclude that $\text{Var}_{\BbbP}(\langle\zeta^{D,\text{ext}}_N,f\rangle)$
tends to zero in $P^\varrho$-probability. This implies
\begin{equation}
\label{E:4.20}
\langle\zeta^{D,\text{ext}}_N,f\rangle-\E\bigl(\langle\zeta^{D,\text{ext}}_N,f\rangle\bigr)\,\,\underset{N\to\infty}\longrightarrow\,\,0,\quad\text{in $P^\varrho\otimes\BbbP$-probability}.
\end{equation}
To infer the desired claim, abbreviate
\begin{equation}
f_{\mathfrak g}(x,\ell):=\int\mathfrak g(\textd h) f(x,\ell,h)
\end{equation}
and note that, since $A$ in \eqref{E:4.8} is compact,~$h^{D_N}_x/\sqrt{\log N}$ tends in law to~$\NN(0,g)$ uniformly for all $x\in\{y\in D_N\colon y/N\in A\}$.  The continuity of~$f$ along with \eqref{E:4.8} yield 
\begin{equation}
\label{E:4.22}
\E\bigl(\langle\zeta^{D,\text{ext}}_N,f\rangle\bigr)-\langle\zeta^D_N,f_{\mathfrak g}\rangle\,\,\underset{N\to\infty}\longrightarrow\,\,0,\quad\text{in $P^\varrho$-probability}.
\end{equation}
Combining \eqref{E:4.20} and \eqref{E:4.22} we then get \eqref{E:4.7}.
\end{proofsect}

As a consequence of the above lemmas, we now get:

\begin{lemma}
\label{lemma-conv}
Recall that~$\mathfrak g$ is the law of~$\NN(0,g)$. Given~$f\in C_\cc(D\times\R)$ with~$f\ge0$, let
\begin{equation}
\label{E:4.24}
f^{\ast\mathfrak g}(x,\ell):=\int\mathfrak g(\textd h)f\Bigl(x,\tfrac1{2\sqrt g(\sqrt\theta+\lambda)}\bigl(\ell+\tfrac{h^2}2\bigr)\Bigr).
\end{equation}
Then for every subsequential weak limit~$\zeta^D$ of~$\zeta^D_N$, simultaneously for all~$f$ as above,
\begin{equation}
\label{E:4.25a}
\langle\zeta^D,f^{\ast\mathfrak g}\rangle\laweq \cspecial(\lambda)\int Z^D_\lambda(\textd x)\otimes\texte^{-\alpha\lambda h}\textd h \,\,f(x,h),
\end{equation}
where, we recall,~$\alpha:=2/\sqrt g$ and~$\cspecial(\lambda)$ is as in \eqref{E:1.19}.
\end{lemma}

\begin{proofsect}{Proof}
Pick~$f$ as above.  Suppressing, for the duration of this proof, the index~$D_N$ on the fields and the local time,  let the DGFF~$\wt h$ in~$D_N$ be related to the local time~$L_{t_N}$ and an independent DGFF~$h$ in~$D_N$ via the Dynkin isomorphism. Recalling \eqref{E:4.1}, for large enough~$N\ge1$ we then have
\begin{equation}
\label{E:4.26}
\begin{aligned}
\langle\wh\eta_N^D, f\rangle&=\frac1{W_N}\sum_{x\in D_N}f\bigl(\ffrac xN,\wt h_x-\wh a_N\bigr)
\\
&=\frac1{W_N}\sum_{x\in D_N}f\Bigl(\ffrac xN,\sqrt{2L_{t_N}(x)+h_x^2}-\sqrt{2a_N}\Bigr)
 1_{\{\wt h_x=\sqrt{2L_{t_N}(x)+h_x^2}-  \sqrt{2t_N}  \}}
\\
&=\frac1{W_N}\sum_{x\in D_N}f\Bigl(\ffrac xN,\sqrt{2L_{t_N}(x)+h_x^2}-\sqrt{2a_N}\Bigr)
\\
&\quad-\frac1{W_N}\sum_{x\in D_N}f\Bigl(\ffrac xN,\sqrt{2L_{t_N}(x)+h_x^2}-\sqrt{2a_N}\Bigr)
1_{\{\wt h_x=-\sqrt{2L_{t_N}(x)+h_x^2}-  \sqrt{2t_N}  \}},
\end{aligned}
\end{equation}
 where we noted that only the positive sign 
in $\wt h_x=\pm\sqrt{2L_{t_N}(x)+h_x^2}-  \sqrt{2t_N} $ 
can contribute in the second line once~$N$ is large due to~$f$ having a compact support and the fact that~$\wh a_N\to\infty$ implied by~$\lambda>0$. 

We start by treating the second term on the extreme right of \eqref{E:4.26} which we note is bounded in absolute value by
\begin{equation}
\label{E:5.29uoi}
\Vert f\Vert_\infty\frac1{W_N}\sum_{x\in D_N}1_{\{\wt h_x\le-\sqrt{2a_N}\}}.
\end{equation}
The result of \cite{BL4}, or even just a simple first-moment estimate, shows that the sum is at most $N^{2[1-(\sqrt\theta+\lambda)^2]+o(1)}$ with high probability. As~$W_N=N^{2(1-\lambda^2)+o(1)}$ and~$\theta>0$, the expression in \eqref{E:5.29uoi} tends to zero in probability as~$N\to\infty$.

We thus need to extract the limit of the first term on the right of \eqref{E:4.26}.  For this we need to first truncate~$h_x$ to values of order~$\sqrt{\log N}$. Let   $\chi\colon[0,\infty)\to[0,1]$ be  non-increasing, continuous with~$\chi(x)=1$ for~$0\le x\le1$ and $\chi(x)=0$ for~$x\ge2$. Then the first term on the right of \eqref{E:4.26} can be written as
\begin{equation}
\label{E:4.27}
\frac1{W_N}\sum_{x\in D_N}f\Bigl(\ffrac xN,\sqrt{2a_N+2[\,L_{t_N}(x)-a_N]+h_x^2}-\sqrt{2a_N}\Bigr)\chi\Bigl(\frac{|h_x|}{M\sqrt{\log N}}\Bigr)
\end{equation}
 plus a quantity bounded, in absolute value, by
\begin{equation}
\label{E:5.32oi}
\Vert f\Vert_\infty\frac1{W_N}\sum_{x\in D_N}1_{\{\sqrt{2L_{t_N}(x)+h_x^2}\ge\sqrt{2a_N}-b\}}1_{\{|h_x|\ge M\sqrt{\log N}\}},
\end{equation}
where~$b>0$ is such that $\supp(f)\subseteq D\times[-b,b]$. Lemma~\ref{lemma-no-conspire} shows that the $L^1$-norm of \eqref{E:5.32oi} under $P^\varrho\otimes\BbbP$ is of order~$\texte^{-\beta M^2}$, uniformly in~$N\ge1$, and so we just need to focus on taking the~$N\to\infty$ limit of \eqref{E:4.27}.

The truncation ensures that, for~$x$ to contribute to the sum  in \eqref{E:4.27},  \myemph{both}~$h_x^2$ and $L_{t_N}(x)-a_N$ must be at most order~$\log N$. Expanding the square root and using the uniform continuity of~$f$ along with the tightness of~$\zeta^D_N$ to replace $a_N$ by its asymptotic expression then recasts \eqref{E:4.27} as
\begin{equation}
\label{E:4.28}
\frac1{W_N}\sum_{x\in D_N}f_{\text{ext}}\Bigl(\ffrac xN,\tfrac{L_{t_N}(x)-a_N}{\log N}, \tfrac{h_x}{\sqrt{\log N}}\Bigr)\chi\Bigl(\frac{|h_x|}{M\sqrt{\log N}}\Bigr)+o(1),
\end{equation}
where
\begin{equation}
f_{\text{ext}}(x,\ell,h):=f\Bigl(x,\tfrac1{2\sqrt g(\sqrt\theta+\lambda)}\bigl(\ell+\tfrac{h^2}2\bigr)\Bigr).
\end{equation}
The function $\ell,h\mapsto f_{\text{ext}}(x,\ell,h)\chi(|h|/M)$ that effectively appears in \eqref{E:4.28} is  compactly supported in both variables; Lemma~\ref{lemma-add-field} then shows that, along subsequences where~$\zeta^D_N$ converges in law to some~$\zeta^D$, the expression in \eqref{E:4.28} converges to $\langle\zeta^D, f^{\ast\mathfrak g}_M \rangle$ where~$f^{\ast\mathfrak g}_M$ is defined by \eqref{E:4.24} with~$\mathfrak g(\textd h)$ replaced by $\chi(|h|/M)\mathfrak g(\textd h)$. 
From the known convergence of~$\wh\eta_N^D$ (see \eqref{E:1.19}) we thus conclude 
\begin{equation}
\langle\zeta^D,f^{\ast\mathfrak g}_M\rangle +O(\texte^{-\beta M^2})\,\laweq\, \cspecial(\lambda)\int Z^D_\lambda(\textd x)\otimes\texte^{-\alpha\lambda h}\textd h \,\,f(x,h),
\end{equation}
where $O(\texte^{-\beta M^2})$ is a random quantity with  $L^1$-norm  at most a constant times~$\texte^{-\beta M^2}$. Taking~$M\to\infty$ via  the  Monotone Convergence Theorem now gives \eqref{E:4.25a}.
\end{proofsect}

Working towards the proof of Theorem~\ref{thm-thick}, a key  remaining point to show is that the class of~$f^{\ast\mathfrak g}$ arising from functions~$f$ for which the integral on the right of \eqref{E:4.25a} converges absolutely is sufficiently rich so that \eqref{E:4.25a} determines the measure~$\zeta^D$ uniquely. For this we note that, by an application of the Dominated and Monotone Convergence Theorems, \eqref{E:4.25a} extends from~$C_\cc(D\times\R)$ to the class of functions~$(x,h)\mapsto 1_A(x)f(h)$, where $A\subseteq   D$ is open with~$\overline A\subseteq   D$ and $f\in C_\cc^\infty(\R)$ with~$f\ge0$. The transformation \eqref{E:4.24} only affects the second variable on which it takes the form $f\mapsto (f\ast\frake)\circ\fraks$, where the convolution is with the function
\begin{equation}
\label{E:5.33}
 \frake(z):=\sqrt{\frac\beta\pi}\,\frac{\texte^{\beta z}}{\sqrt{-z}}1_{(-\infty,0)}(z)\quad\text{for}\quad\beta:=\alpha\bigl(\sqrt\theta+\lambda\bigr)
\end{equation}
and where $h\mapsto\fraks(h)$ is the scaling map
\begin{equation}
\label{E:5.34}
\fraks(h):=\frac{h}{2\sqrt g(\sqrt\theta+\lambda)}.
\end{equation}
As it turns out, it then suffices to observe:

\begin{lemma}
\label{lemma-unique}
 Denote  $\mu_\lambda(\textd h):=\texte^{-\alpha\lambda h}\textd h$  and let~$\frake(\cdot)$ be as in \eqref{E:5.33} with~$\beta>\alpha\lambda$. Then  there is at most one Radon measure~$\nu$ on~$\R$ such that for all~$f\in C_\cc^\infty(\R)$ with $f\ge0$,
\begin{equation}
\label{E:5.35iw}
\bigl\langle\nu,f\ast\frake\bigr\rangle=\langle\mu_\lambda,f\rangle.
\end{equation}
\end{lemma}

\begin{proofsect}{Proof}
Writing \eqref{E:5.35iw} explicitly using integrals and using the fact that the class of all $f\in C_\cc^\infty(\R)$ with $f\ge0$ separates Radon measures on~$\R$ shows
\begin{equation}
\int_{\R}\nu(\textd s)\frake(s-h) = \texte^{-\alpha\lambda h},\quad h\in\R.
\end{equation}
Abbreviating $\nu_{\lambda}(\textd h):=\texte^{\alpha\lambda h}\nu(\textd h)$ and $\frake_{\lambda}(h):=\texte^{-\alpha\lambda h}\frake(h)$, this can be recast as
\begin{equation}
\int_{\R}\nu_{\lambda}(\textd s)\frake_{\lambda}( s-h ) = 1,\quad h\in\R.
\end{equation}
Integrating this against suitable test functions with respect to the Lebesgue measure and applying 
 the  Dominated Convergence Theorem, we conclude
\begin{equation}
\bigl\langle\nu_\lambda,f\ast\frake_\lambda\bigr\rangle=\langle\leb,f\rangle,\quad f\in \CalS(\R),
\end{equation}
where~$\CalS(\R)$ is the Schwartz class of functions on~$\R$. Note that this identity entails that the integral on the left-hand side converges absolutely.

Since $\CalS(\R)$ separates Radon measures on~$\R$, to conclude the statement from \eqref{E:5.35iw} it suffices to prove that, for $\theta>0$,
\begin{equation}
\label{E:5.39iw}
f\mapsto f\ast\frake_\lambda\text{ is a bijection of }\CalS(\R) \text{ onto itself}.
\end{equation}
The Fourier transform maps~$\CalS(\R)$ bijectively onto itself and so we may as well prove \eqref{E:5.39iw} in the Fourier picture. For this we note that, as $\theta>0$ we have~$\tilde\beta:=\beta-\alpha\lambda>0$ and so $z\mapsto\frake_\lambda(z)$ decays exponentially as~$z\to-\infty$. In particular, $\frake_\lambda$ is integrable and so in the Fourier transform,~$f\mapsto f\ast\frake_\lambda$ is reduced to the multiplication by
\begin{equation}
\label{E:5.40}
 \wh\frake_\lambda(k):=\int_\R\textd z\,\frake_\lambda(z)\texte^{2\pi\texti kz}
= \sqrt{\frac\beta\pi}\int_{(0,\infty)}\textd x\,\frac1{\sqrt{x}}\,\texte^{-\tilde\beta(1+2\pi\texti k/\tilde\beta) x}.
\end{equation}
 Hereby we readily check that $k\mapsto\wh\frake_\lambda(k)$ is~$C^\infty(\R)$ with bounded derivatives  which  implies that $f\mapsto  \frake_{\lambda} \ast f$ maps~$\CalS(\R)$ into~$\CalS(\R)$.  Using the substitution $x=y^2$ and computing the complex-Gaussian integral we find that
\begin{equation}
\bigl|\wh\frake_\lambda(k)\bigr|=\sqrt{\frac\beta{\tilde\beta}}\,\frac1{|1+2\pi\texti k/\tilde\beta|^{1/2}}.
\end{equation}
 As~$|\wh\frake_\lambda(k)|>0$ for all $k\in\R$, the map $f\mapsto  \frake_{\lambda} \ast f$ is injective; the fact that $|\wh\frake_\lambda(k)|^{-1}$ is bounded by a power of~$|k|$ then shows that it is also onto.  Hence \eqref{E:5.39iw} follows. 
\end{proofsect}

We are now ready to give:

\begin{proofsect}{Proof of Theorem~\ref{thm-thick}}
Consider a subsequential limit~$\zeta^D$, pick~$f\in C_\cc(\R)$ with~$f\ge0$ and let~$A\subseteq   D$ be open with~$\overline A\subseteq   D$. Using the notation \twoeqref{E:5.33}{E:5.34} we then have 
\begin{equation}
\bigl\langle\zeta^D,(1_A\otimes f)^{\ast\mathfrak g}\bigr\rangle=\bigl\langle\zeta_A^D,(f\ast\frake)\circ\fraks\bigr\rangle=\bigl\langle\zeta_A^D\circ\fraks^{-1},f\ast\frake\bigr\rangle,
\end{equation}
where~$\zeta_A^D$ is a Borel measure on~$\R$ defined by $\zeta_A^D(B):=\zeta^D(A\times B)$. Writing $\mu_\lambda(\textd h):=\texte^{-\alpha\lambda h}\textd h$, the identity \eqref{E:4.25a} then translates into
\begin{equation}
\label{E:5.43nw}
\bigl\langle\zeta_A^D\circ\fraks^{-1},f\ast\frake\bigr\rangle \laweq \cspecial(\lambda)Z^D_\lambda(A)\langle\mu_\lambda,f\rangle,
\end{equation}
where the equality in law holds simultaneously for all~$A$ and~$f$ as above. 

To infer the product form of~$\zeta^D$ from \eqref{E:5.43nw}, define (for a given~$A$ and a given realization of~$\zeta^D$) a Borel measure on~$\R$ by
\begin{equation}
\label{E:5.44nw}
\nu:=\Bigl[\alpha\lambda \bigl\langle\zeta_A^D\circ\fraks^{-1},1_{[0,\infty)}\ast\frake\bigr\rangle\Bigr]^{-1}\zeta^D_A,
\end{equation}
where the conditions on~$A$ imply $Z^D_\lambda(A)>0$ a.s.\ and so, by \eqref{E:5.43nw}, the quantity in the square bracket is strictly positive a.s. By \eqref{E:5.43nw} we have $\langle\nu\circ\fraks^{-1},f\ast\frake\rangle=\langle\mu_\lambda,f\rangle$ for all~$f\in C_\cc(\R)$ and so, by Lemma~\ref{lemma-unique}, $\nu\circ\fraks^{-1}$, and thus also~$\nu$, is determined uniquely. In particular, $\nu$ is the same for all~$A$ as above and for a.e.~realization of~$\zeta^D$. Using \eqref{E:5.43nw} in \eqref{E:5.44nw} then shows $\zeta^D_A(\textd h)\laweq\cspecial(\lambda)Z^D_\lambda(A)\nu(\textd h)$. As this holds simultaneously for all~$A$ as above, Remark~\ref{rem-restrict} permits us to conclude
\begin{equation}
\zeta^D\laweq \cspecial(\lambda)Z^D_\lambda\,\otimes\,\nu,
\end{equation}
where~$\nu$ is a uniquely-determined deterministic Radon measure on~$\R$. 

It remains to derive the explicit form of~$\nu$ which, thanks to its uniqueness, we can do by plugging the desired expression on the left-hand side of \eqref{E:4.25a} and checking for equality. Abbreviate $\tilde\alpha:=\alpha(\theta,\lambda)$ and note that
\begin{equation}
\tilde\alpha = \frac1{2\sqrt g(\sqrt\theta+\lambda)}\alpha\lambda.
\end{equation}
Pick~$f\in C_\cc(D\times\R)$ and perform  the following calculation where, in the last step, we invoke the substitution $r:=\frac1{2\sqrt g(\sqrt\theta+\lambda)}(\ell+\frac{h^2}2)$ and separate integrals using Fubini-Tonelli:
\begin{equation}
\label{E:4.38}
\begin{aligned}
\int_{D\times\R}&Z^D_\lambda(\textd x)\otimes\texte^{-\tilde \alpha\ell}\textd\ell\,f^{\ast\mathfrak g}(x,\ell)
\\
&=\int_{D\times\R\times\R}Z^D_\lambda(\textd x)\otimes\texte^{-\tilde \alpha\ell}\textd\ell\otimes\mathfrak g(\textd h)\,
f\Bigl(x,\tfrac1{2\sqrt g(\sqrt\theta+\lambda)}\bigl(\ell+\tfrac{h^2}2\bigr)\Bigr)
\\
&=\int_{D\times\R\times\R}Z^D_\lambda(\textd x)\otimes\texte^{-\tilde\alpha(\ell+\frac{h^2}2)}\textd\ell\otimes\texte^{\,\tilde\alpha\frac{h^2}2}\mathfrak g(\textd h)\,
f\Bigl(x,\tfrac1{2\sqrt g(\sqrt\theta+\lambda)}\bigl(\ell+\tfrac{h^2}2\bigr)\Bigr)
\\
&=2\sqrt g(\sqrt\theta+\lambda)\Bigl(\int_\R\mathfrak g(\textd h)\texte^{\,\tilde\alpha\frac{h^2}2}\Bigr)
\int_{D\times\R}Z^D_\lambda(\textd x)\otimes\texte^{-\alpha \lambda r}\textd r\,f(x,r)\,.
\end{aligned}
\end{equation}
As $\tilde\alpha<1/g$, the first integral on the last line converges to the root of  $(1-\tilde\alpha g)^{-1}=\frac{\sqrt\theta+\lambda}{\sqrt\theta}$ while \eqref{E:4.25a} equates the second integral to $\cspecial(\lambda)^{-1}\langle\zeta^D,f^{\ast\mathfrak g}\rangle$ in law. This implies 
\begin{equation}
 \zeta^D \,\laweq\, \frac{\theta^{1/4}}{2\sqrt{g}\,(\sqrt\theta+\lambda)^{3/2}}\,
\cspecial(\lambda)Z^D_\lambda(\textd x)\otimes\texte^{-\tilde \alpha\ell}\textd\ell. 
\end{equation}
In particular, all weakly converging subsequences of~$\{\zeta^D_N\colon N\ge1\}$ converge to this~$\zeta^D$, thus proving the desired claim.
\end{proofsect}

\section{Thin points}
\label{sec5}\noindent\nopagebreak
 Our next task is  the convergence of point measures  $\zeta^D_N$  associated with $\lambda$-thin  points.  The argument proceeds very much along the same sequence of lemmas as for the $\lambda$-thick points and so we  will concentrate on the steps where a different reasoning is needed.  Throughout we assume that~$t_N$ and~$a_N$ are sequences satisfying \eqref{E:1.22} with some~$\theta>0$ and some~$\lambda\in(0, 1 \wedge \sqrt\theta)$.  The auxiliary centering sequence~$\wh a_N$ is now defined by 
\begin{equation}
\label{E:5.1}
\wh a_N:=\sqrt{2t_N}-\sqrt{2a_N}
\end{equation}
which ensures that we still have $\wh a_N\sim 2\lambda\sqrt g\log N$ as~$N\to\infty$.  Appealing to the coupling of~$L^{D_N}_{t_N}$ and~$h^{D_N}$ to~$\wt h^{D_N}$ via  the  Dynkin isomorphism, we use~$\wh\eta_N^D$ to denote the point process associated with~$\wt h^{D_N}$ and the centering sequence~$-\wh a_N$.  

The proof again opens up by proving suitable tightness and joint-convergence statements. We start with an analogue of Lemma~\ref{lemma-no-conspire}: 

\begin{lemma}
\label{lemma-no-conspire2}
Let $0<\beta<\frac1{2g}\frac{ \lambda }{\sqrt\theta-\lambda}$. Then for each~$b\in\R$ there is~$c_5(b)\in(0,\infty)$  and, for each $M\ge0$, there is~$N'=N'(b,M)$ such that for  all~$N\ge N'$ and all~$x\in D_N$,
\begin{equation}
P^\varrho\otimes\BbbP\biggl(L^{D_N}_{t_N}(x)+\frac{(h^{D_N}_x)^2}2 \le  a_N+b\log N,\,\frac{|h^{D_N}_x|}{\sqrt{\log N}}\ge M\biggr)\le  c_5(b)\frac{W_N}{N^2}\texte^{-\beta M^2}.
\end{equation}
\end{lemma}

\begin{proofsect}{Proof}
Let us again for simplicity just deal with the case~$b=0$. Pick~$0<\delta<\sqrt\theta-\lambda$. Then the probability in question is bounded by
\begin{multline}
\qquad
P^\varrho\Bigl(L^{D_N}_{t_N}(x)\le 2g(\sqrt\theta-\lambda-\delta)^2(\log N)^2\Bigr)
\\+P^\varrho\Bigl(2g(\sqrt\theta-\lambda-\delta)^2(\log N)^2\le L^{D_N}_{t_N}(x)\le a_N-\frac12M^2\log N\Bigr).
\qquad
\end{multline}
Invoking the calculation in \eqref{E:3.18}, the first term is at most order~$N^{-2(\lambda+\delta)^2+o(1)}$ which is $o(W_N/N^2)$. The second term is now bounded using Lemma~\ref{lemma-lower} and the fact that, by the uniform bound $G^{D_N}(x,x) \le g \log N + c$ with~$c$ independent of~$N$, we have 
\begin{equation}
\min_{x\in D_N}
\frac{\log N\bigl(\sqrt{2t_N}-\sqrt{2a_N}\bigr)}{G^{D_N}(x,x)\sqrt{2a_N}}\ge\frac1g\,\frac{\lambda}{\sqrt\theta-\lambda}+o(1)
\end{equation}
in the limit~$N\to\infty$. 
Indeed, this  shows that the last exponential in \eqref{E:4.7a2} for the choice 
$b:=-\frac12M^2\log N$ is less than $\texte^{-\beta M^2}$ once~$N$ is sufficiently large.
\end{proofsect}

Next we will give an analogue of Lemma~\ref{lemma-add-field} which we restate \textit{verbatim}, albeit with a somewhat different proof: 

\begin{lemma}
\label{lemma-add-field2}
Suppose $\{N_k\}$ is a subsequence along which~$\zeta^D_N$ converges in law to~$\zeta^D$. Then
\begin{equation}
\frac1{W_N}\sum_{x\in D_N}\delta_{x/N}\otimes\delta_{(L^{D_N}_{t_N}(x)-a_N)/\log N}\otimes\delta_{h_x^{D_N}/\sqrt{\log N}}\,\,\underset{\begin{subarray}{c}
N=N_k\\k\to\infty
\end{subarray}}{\overset{\text{\rm law}}\longrightarrow}\,\,
\zeta^D\otimes\mathfrak g,
\end{equation}
where~$\mathfrak g$ is the law of~$\NN(0,g)$.
\end{lemma}

\begin{proofsect}{Proof}
Let~$\zeta^{D,\text{ext}}_N$ denote the measure on the left and let~$f\in C_\cc(D\times\R\times\R)$  obey~$f\ge0$.  As for Lemma~\ref{lemma-add-field}, the argument hinges on proving
\begin{equation}
\text{\rm Var}_{\BbbP}\bigl(\langle\zeta^{D,\text{ext}}_N,f\rangle\bigr)\,\underset{N\to\infty}\longrightarrow\,0,\quad \text{in $P^\varrho$-probability},
\end{equation}
where~$\text{\rm Var}_{\BbbP}$ denotes the variance with respect to the law of~$h^{D_N}$, conditional on~$L^{D_N}_{t_N}$. 
Invoking \eqref{E:4.11}, we treat the sum over the pairs~$|x-y|\ge\epsilon N$ via  the argument following \twoeqref{E:5.12i}{E:5.13i}. 
The key difference is that we no longer have the domination of~$\zeta^D_N$ by a DGFF process  in this case  and so we have to control the sum over the pairs $x,y\in D_N$ with~$|x-y|\le\epsilon N$ differently. 

Since~$f$ is non-negative and compactly supported, we in fact just need to show that, for any~$M>0$, the $L^1(\BbbP)$-norm of 
\begin{multline}
\label{E:5.3}
\quad\frac1{W_N^2}\sum_{\begin{subarray}{c}
x,y\in D_N\\|x-y|\le\epsilon N
\end{subarray}}
1_{\{L^{D_N}_{t_N}(x)\le a_N+M^2\log N\}}1_{\{|h^{D_N}_x|\le M\sqrt{\log N}\}}
\\
\times1_{\{L^{D_N}_{t_N}(y)\le a_N+M^2\log N\}}1_{\{|h^{D_N}_y|\le M\sqrt{\log N}\}}
\quad
\end{multline}
vanishes in $P^\varrho$-probability in the limit as~$N\to\infty$ and~$\epsilon\downarrow0$. To this end  we note that, dropping the indicators involving the DGFF, \eqref{E:5.3} is bounded by $[\zeta^D_N(D\times(-\infty,M^2])]^2$ which by Corollary~\ref{cor-tightness-lower} is  bounded in probability as~$N\to\infty$.  Therefore, it suffices to prove that \eqref{E:5.3} vanishes in the stated limits in $P^\varrho\otimes\BbbP$-probability. 

To this end pick~$b>\frac{M^2}{\sqrt{g}(\sqrt{\theta}-\lambda)}$ and note that, as soon as~$N$ is sufficiently large, the asymptotic forms of~$a_N$ along with the Dynkin isomorphism yield 
\begin{multline}
\qquad
1_{\{L^{D_N}_{t_N}(x)\le a_N+M^2\log N\}}1_{\{|h^{D_N}_x|\le M\sqrt{\log N}\}}
\\
\le 1_{\{(\wt h^{D_N}_x+\sqrt{2t_N})^2\le 2a_N+ 3 M^2\log N\}}
\le 1_{\{\wt h^{D_N}_x\le -\wh a_N+b\}}.
\qquad
\end{multline}
 It follows that \eqref{E:5.3} is bounded by 
\begin{equation}
\label{E:6.9uo}
\wh\eta_N^D\otimes\wh\eta_N^D\Bigl(\bigl\{(x,y)\colon|x-y|\le\epsilon\bigr\}\times(-\infty,b]^2\Bigr)
\end{equation}
 whose~$N\to\infty$ and~$\epsilon\downarrow0$ limits are now handled as before.
\end{proofsect}

Our next task is a derivation of a convolution identity that will, as for the thick points, ultimately characterize the limit measure uniquely: 

\begin{lemma}
\label{lemma-conv2}
Given~$f\in C_\cc(D\times\R)$ with~$f\ge0$, let (abusing our earlier notation) 
\begin{equation}
f^{\ast\mathfrak g}(x,\ell):=\int\mathfrak g(\textd h)f\Bigl(x,\tfrac1{2\sqrt g(\sqrt\theta-\lambda)}\bigl(\ell+\tfrac{h^2}2\bigr)\Bigr).
\end{equation}
Then for every subsequential weak limit~$\zeta^D$ of~$\zeta^D_N$, simultaneously for all~$f$ as above, 
\begin{equation}
\label{E:4.25b}
\langle\zeta^D, f^{\ast\mathfrak g}\rangle\laweq \cspecial(\lambda)\int Z^D_\lambda(\textd x)\otimes\texte^{\alpha\lambda h}\textd h \,\,f(x,h),
\end{equation}
where~$\alpha:=2/\sqrt g$ and~$\cspecial(\lambda)$ is as in \eqref{E:1.19}.
\end{lemma}

\begin{proofsect}{Proof}
Pick~$f$ as above and let~$\chi$ be the function as in the proof of Lemma~\ref{lemma-conv}. The fact that~$f$ has compact support gives
\begin{equation}
\frac1{W_N}\sum_{x\in D_N}f\bigl(\ffrac xN, \wt h^{D_N}_x+\wh a_N\bigr)
=\frac1{W_N}\sum_{x\in D_N} f\Bigl(\ffrac xN,-\sqrt{2a_N}+\sqrt{2L^{D_N}_{t_N}(x)+(h^{D_N}_x)^2}\Bigr)
\end{equation}
and Lemma~\ref{lemma-no-conspire2} then bounds this by~$O(\texte^{-\beta M^2})$ 
plus
\begin{equation}
\frac1{W_N}\sum_{x\in D_N} f\Bigl(\ffrac xN,-\sqrt{2a_N}+\sqrt{2L^{D_N}_{t_N}(x)+(h^{D_N}_x)^2}\Bigr)\chi\Bigl(\frac{|h_x^{D_N}|}{M\sqrt{\log N}}\Bigr).
\end{equation}
The truncation of the field now forces $L^{D_N}_{t_N}-a_N$ to be  of  order $\log N$.
Expanding the square root and using the uniform continuity with the help of Corollary~\ref{cor-tightness-lower} rewrites this~as
\begin{equation}
\label{E:5.28}
\frac1{W_N}\sum_{x\in D_N}\tilde f_{\text{ext}}\Bigl(\ffrac xN,\tfrac{L^{D_N}_{t_N}(x)-a_N}{\log N}, \tfrac{h_x^{D_N}}{\sqrt{\log N}}\Bigr)\chi\Bigl(\frac{|h_x^{D_N}|}{M\sqrt{\log N}}\Bigr)+o(1),
\end{equation}
where
\begin{equation}
\tilde f_{\text{ext}}(x,\ell,h):=f\Bigl(x,\tfrac1{2\sqrt g(\sqrt\theta-\lambda)}\bigl(\ell+\tfrac{h^2}2\bigr)\Bigr).
\end{equation}
The rest of the proof now proceeds as before. (The exponential on the right-hand side of \eqref{E:4.25b} does not get a negative sign because~$\wh\eta_N^D$ is centered along negative sequence of order~$\log N$.)
\end{proofsect}

Using  the  Dominated and Monotone Convergence  Theorems,  we now readily extend \eqref{E:4.25b} to functions of the form $1_A\otimes f$ where~$A\subseteq   D$ is open with~$\overline A\subseteq   D$ and~$f\in C_\cc(\R)$ obeys~$f\ge0$. For such~$f$ we then get
\begin{equation}
(1_A\otimes f)^{\ast\mathfrak g}=1_A\otimes (f\ast\frake')\circ\fraks'
\end{equation}
where~$\frake'$ is given by the same formula as~$\frake$ in \eqref{E:5.33} but with~$\beta$ replaced by
\begin{equation}
\beta':=\alpha\bigl(\sqrt\theta-\lambda\bigr)
\end{equation}
and~$\fraks'(h):=h/(2\sqrt g(\sqrt\theta-\lambda))$. We then state:

\begin{lemma}
\label{lemma-unique2}
 Denote  $\mu_\lambda'(\textd h):=\texte^{\alpha\lambda h}\textd h$  and let~$\frake'$ be as above with~$\beta'>-\alpha\lambda$. Then  there is at most one Radon measure~$\nu$ on~$\R$ such that for all~$f\in C_\cc^\infty(\R)$ with $f\ge0$,
\begin{equation}
\label{E:5.35iw2}
\bigl\langle\nu,f\ast\frake'\bigr\rangle=\langle\mu_\lambda',f\rangle.
\end{equation}
\end{lemma}

\begin{proofsect}{Proof}
As in the proof of Lemma~\ref{lemma-unique}, we recast \eqref{E:5.35iw2} as
\begin{equation}
\bigl\langle\nu_\lambda,f\ast\frake'_\lambda\bigr\rangle=\langle\leb,f\rangle
\end{equation}
where $\nu_\lambda(\textd h)=\texte^{-\alpha\lambda h} \nu(\textd h)$ and $\frake_\lambda'(h)=\texte^{\alpha\lambda h}\frake'(h)$. Since $\tilde\beta':=\beta'+\alpha\lambda >0$, we again get that~$\frake_\lambda'$ is integrable. Replacing~$\tilde\beta$ by~$\tilde\beta'$, the rest of the argument is then identical to that in the proof of Lemma~\ref{lemma-unique}.
\end{proofsect}

We are now ready to give:

\begin{proofsect}{Proof of Theorem~\ref{thm-thin}}
The argument proving that \eqref{E:4.25b} determines $\zeta^D$ uniquely is the same as for the thick points so we just need to perform the analogue of the calculation in \eqref{E:4.38}. Denoting, for the duration of this proof, 
\begin{equation}
\wh\alpha:=\frac1{2\sqrt g\,(\sqrt\theta-\lambda)}\alpha\lambda,
\end{equation}
we get
\begin{equation}
\label{E:4.38a}
\begin{aligned}
\int_{D\times\R}&Z^D_\lambda(\textd x)\otimes\texte^{\wh \alpha\ell}\textd\ell\,\,  f^{\ast\mathfrak g} (x,\ell)
\\
&=\int_{D\times\R\times\R}Z^D_\lambda(\textd x)\otimes\texte^{\wh \alpha\ell}\textd\ell\otimes\mathfrak g(\textd h)\,
f\Bigl(x,\tfrac1{2\sqrt g(\sqrt\theta-\lambda)}\bigl(\ell+\tfrac{h^2}2\bigr)\Bigr)
\\
&=\int_{D\times\R\times\R}Z^D_\lambda(\textd x)\otimes\texte^{\wh\alpha(\ell+\frac{h^2}2)}\textd\ell\otimes\texte^{\,-\wh\alpha\frac{h^2}2}\mathfrak g(\textd h)\,
f\Bigl(x,\tfrac1{2\sqrt g(\sqrt\theta-\lambda)}\bigl(\ell+\tfrac{h^2}2\bigr)\Bigr)
\\
&=2\sqrt g(\sqrt\theta-\lambda)\Bigl(\int_\R\mathfrak g(\textd h)\texte^{\,-\wh\alpha\frac{h^2}2}\Bigr)
\int_{D\times\R}Z^D_\lambda(\textd x)\otimes\texte^{\alpha \lambda r}\textd r\,f(x,r)\,.
\end{aligned}
\end{equation}
The Gaussian integral on the last line equals the root of $\frac{\sqrt\theta-\lambda}{\sqrt\theta}$. It follows that~$\zeta^D_N$ converges in law to the measure
\begin{equation}
\frac{\theta^{1/4}}{2\sqrt{g}\,(\sqrt\theta-\lambda)^{3/2}}\,
\cspecial(\lambda)Z^D_\lambda(\textd x)\otimes\texte^{\wh \alpha\ell}\textd\ell.
\end{equation}
This is the desired claim.
\end{proofsect}

\section{Light and avoided points}
\label{sec6}\noindent\nopagebreak
In this section we will deal with the point measures $\vartheta^D_N$ and~$\kappa^D_N$ associated with the light and avoided points, respectively. The argument follows the blueprint of the proof for the~$\lambda$-thick and~$\lambda$-thin points although important  changes  arise due to a different scaling of~$\wh W_N$ with~$N$ compared to~$W_N$.  As before, a key point of the argument is the extension of the convergence by adding information about an independent DGFF. The difference now is that this field comes without any normalization: 

\begin{lemma}
\label{lemma-7.1}
Suppose $\{N_k\}$ is a subsequence along which~$\vartheta^D_N$ converges in law to~$\vartheta^D$. Then
\begin{equation}
\label{E:7.1}
\frac{\sqrt{\log N}}{\wh W_N }\sum_{x\in D_N}\delta_{x/N}\otimes\delta_{L^{D_N}_{t_N}(x)}\otimes\delta_{h^{D_N}_x}
\,\,\underset{\begin{subarray}{c}
N=N_k\\k\to\infty
\end{subarray}}{\overset{\text{\rm law}}\longrightarrow}\,\,
\vartheta^D\otimes\frac1{\sqrt{2\pi g}}\,\leb.
\end{equation}
\end{lemma}

\begin{proofsect}{Proof}
Let $\vartheta^{D,\text{ext}}_N$ denote the measure on the left and pick~$f\in C_\cc(D\times[0,\infty)\times\R)$.  Suppose that~$f(x,\ell,h)=0$ unless~$x\in A$, where $A$ is an open set with~$\overline A\subseteq   D$, and unless $\ell,h^2\le M$ for some~$M>0$. Noting that the probability density of~$h^{D_N}_x$ is~$(1+o(1))(2\pi g\log N)^{-1/2}$ with~$o(1)\to0$ as~$N\to\infty$ uniformly over any compact interval shows, with the help of the tightness of~$\{\vartheta^D_N\colon N\ge1\}$ proved in Corollary~\ref{cor-light}, that
\begin{equation}
\label{E:6.2u}
\E\bigl\langle\vartheta^{D,\text{ext}}_N,f\bigr\rangle=o(1)+\frac1{\sqrt{2\pi g}}\bigl\langle\vartheta^D_N\otimes\leb,f\bigr\rangle,
\end{equation}
where~$o(1)\to0$ in $P^\varrho$-probability as~$N\to\infty$. The claim thus reduces to proving concentration of $\langle\vartheta^{D,\text{ext}}_N,f\rangle$ around the (conditional) expectation with respect to~$h^{D_N}$.

Due to the additional $\sqrt{\log N}$ factor in the normalization $\vartheta^{D,\text{ext}}_N$, the domination arguments  for the conditional second moment of $\langle\vartheta^{D,\text{ext}}_N,f\rangle$ of the kind \twoeqref{E:5.12i}{E:4.22uo} for the thick points and \twoeqref{E:5.3}{E:6.9uo} for the thin points seem to fail, so we will instead work with the Laplace transform of $\langle\vartheta^{D,\text{ext}}_N,f\rangle$.  This is motivated by noting that, for~$f\ge0$,  the conditional Jensen inequality and \eqref{E:6.2u} yield
\begin{equation}
\label{E:7.3o}
E^\varrho\otimes\E\bigl(\texte^{-\langle\vartheta^{D,\text{ext}}_N,f\rangle}\bigr)
\ge\texte^{o(1)}E^\varrho\Bigl(\texte^{-(2\pi g)^{-1/2}\langle\vartheta^{D}_N\otimes\leb,f\rangle}\Bigr).
\end{equation}
It thus suffices to derive the opposite inequality which  will require a somewhat technical argument. A key point is to restrict the measure~$\vartheta^{D,\text{ext}}_N$ by a suitable truncation. 

We start with the definition of a truncation event.  Writing temporarily $L(x)$, resp., $h_x$ instead of~$L^{D_N}_{t_N}(x)$, resp.,~$h^{D_N}_x$, given any~$\epsilon,\delta>0$, let
\begin{equation}
F_{N,M,\epsilon,\delta}(x):=\biggl\{\sum_{\begin{subarray}{c}
y\in D_N\\|x-y|<\epsilon N
\end{subarray}}
1_{\{L(y)+\frac12h_y^2\le 2M\}}\le\delta  \frac{\wh W_N}{\sqrt{\log N}}  \biggr\}.
\end{equation}
 We claim that, with probability tending to one as~$N\to\infty$ and~$\epsilon\downarrow0$ (for any~$\delta>0$ fixed), the event $F_{N,M,\epsilon,\delta}(x)$ will not occur for any~$x\in D_N$. For this  let $B_r(x):=\{y\in\R^2\colon |y-x|<r\}$ and let $x_1,\dots,x_m\in D$ be such that $\{B_\epsilon(x_i)\colon i=1,\dots,m\}$ cover~$D$. Writing~$\wh\eta^D_N$ for the DGFF measure associated with the field~$\wt h^{D_N}$ and centering sequence $\{\sqrt{2t_N}\}_{N\ge1}$ and noting that the normalization factor $\wh W_N/\sqrt{\log N}$ in \eqref{E:7.1} then coincides with~$K_N$ (for the centering sequence $\sqrt{2t_N}$), the coupling from Theorem~\ref{thm-Dynkin} yields
\begin{equation}
\bigcup_{x\in D_N}F_{N,M,\epsilon,\delta}(x)^\cc\subseteq  \bigcup_{i=1}^m \Bigl\{\wh\eta^D_N\bigl(B_{2\epsilon}(x_i)\times[-2\sqrt{M},2\sqrt{M}]\bigr)>\delta\Bigr\}.
\end{equation}
Since $\wh\eta^D_N$ is known to converge to a measure with no-atoms,  the probability of the event on the right-hand side tends to zero as~$N\to\infty$ and~$\epsilon\downarrow0$ for any~$\delta>0$, as claimed.  

Introduce the truncated measure 
\begin{equation}
\vartheta^{D,\text{ext}}_{N,M,\epsilon,\delta}:=\frac1{ K_N}
\sum_{x\in D_N}1_{F_{N,M,\epsilon,\delta}(x)}\,\delta_{x/N}\otimes\delta_{L(x)}\otimes\delta_{h_x},
\end{equation}
 where we write $K_N$ for~$\wh W_N/\sqrt{\log N}$. 
We then get
\begin{equation}
\label{E:7.7nw}
\lim_{\epsilon\downarrow0}\,\limsup_{N\to\infty}\,P\Bigl(\langle\vartheta^{D,\text{ext}}_{N,M,\epsilon,\delta},f\rangle\ne\langle \vartheta^{D,\text{ext}}_N,f\rangle\Bigr)=0
\end{equation}
for any~$\delta>0$ (and any~$f$ and~$M$ as above).
Next we will invoke the fact that, for each~$M>0$ and each~$A\subseteq  D$ open with~$\overline A\subseteq D$,
\begin{equation}
\label{E:7.8nw}
\bigl\{\wh\eta^D_N(A\times[-M,M])\colon N\ge1\bigr\}\text{ is uniformly integrable},
\end{equation}
which follows from the convergence in the mean and control of moments implied by~\cite[Lemmas~4.1 and 4.2]{BL4}. Theorem~\ref{thm-Dynkin} then extends \eqref{E:7.8nw} to the uniform integrability of $\{\langle\vartheta^{D,\text{ext}}_N,f\rangle\colon N\ge1\}$. Using \eqref{E:7.7nw} we then get
\begin{multline}
\label{E:7.9nwt}
\quad\lim_{\epsilon\downarrow0}\,\limsup_{N\to\infty}\biggl|E^\varrho\otimes\E\Bigl(\langle\vartheta^{D,\text{ext}}_N,f\rangle\,\texte^{-s\langle\vartheta^{D,\text{ext}}_N,f\rangle}\Bigr)
\\
-E^\varrho\otimes\E\Bigl(\langle\vartheta^{D,\text{ext}}_N,f\rangle\,\texte^{-s\langle\vartheta^{D,\text{ext}}_{N,M,2\epsilon,\delta},f\rangle}\Bigr)\biggr|=0
\quad
\end{multline}
uniformly in~$s\in[0,1]$.  (We write~$2\epsilon$ for reasons to be clear in a moment.) 
As a consequence, we may thus focus on the second expectation from now on.

We first use the explicit form of the measure~$\vartheta^{D,\text{ext}}_N$ and, noting that~$f\ge0$, apply the conditional Jensen inequality as
\begin{equation}
\label{E:7.4uo}
\begin{aligned}
E^\varrho\otimes\E\Bigl(\langle\vartheta^{D,\text{ext}}_N,&f\rangle\,\texte^{-s\langle\vartheta^{D,\text{ext}}_{N,M,2\epsilon,\delta},f\rangle}\Bigr)
\\
&=\frac1{ K_N}\sum_{\begin{subarray}{c}
x\in D_N\\x/N\in A
\end{subarray}}
E^\varrho\otimes\E\Bigl( f\bigl(x/N,L(x),h_x\bigr)
\texte^{-s\langle\vartheta^{D,\text{ext}}_{N,M,2\epsilon,\delta},f\rangle}\Bigr)
\\
&\ge
\frac1{ K_N}\sum_{\begin{subarray}{c}
x\in D_N\\x/N\in A
\end{subarray}}
E^\varrho\otimes\E\Bigl( f\bigl(x/N,L(x),h_x\bigr)
\texte^{-s\E(\langle\vartheta^{D,\text{ext}}_{N,M,2\epsilon,\delta},f\rangle|\sigma(h_x))}\Bigr).
\end{aligned}
\end{equation}
Reflecting on the positivity and support restrictions for~$f$, the conditional expectation in the exponent is dominated via
\begin{multline}
\label{E:7.5iue}
\quad
\E\bigl(\langle\vartheta^{D,\text{ext}}_{N,M,2\epsilon,\delta},f\rangle\,\big|\,\sigma(h_x)\bigr)
\le
\frac1{ K_N}\sum_{\begin{subarray}{c}
y\in D_N\\|x-y|\ge\epsilon N
\end{subarray}}\,\E\Bigl(f\bigl(y/N,L(y),h_y\bigr)\,\Big|\,\sigma(h_x)\Bigr)
\\+
\frac{\Vert f\Vert_{\infty}}{ K_N}\,\E\biggl(\,\sum_{\begin{subarray}{c}
y\in D_N\\|x-y|<\epsilon N
\end{subarray}}1_{\{L(y)+\frac12 h_y^2\le 2M\}}1_{F_{N,M,2\epsilon,\delta}(y)}\,\Big|\,\sigma(h_x)\Bigr).
\end{multline}
As a result of the truncation, since the ball of radius $2\epsilon N$ around any~$y$ with $|y-x|<\epsilon N$ includes the ball of radius~$\epsilon N$ around~$x$, as soon as $F_{N,M,2\epsilon,\delta}(y)$ occurs for at least one~$y$ with $|y-x|<\epsilon N$, the sum in the second term on the right 
is at most~$\delta  K_N $. This bounds the second term on the right of \eqref{E:7.5iue} by~$\delta \Vert f\Vert_{\infty}$ pointwise. 

We have reduced estimating the conditional expectation to a bound on the first term on the right of \eqref{E:7.5iue}. Denoting, for any~$r>0$,
\begin{equation}
\label{E:7.12uiq}
\text{osc}_{f,M}(r):=\sup_{z\in D}\,\,\sup_{\ell\le M}\,\,\sup_{\begin{subarray}{c}
h,h'\in[-\sqrt{M},\sqrt{M}]\\|h-h'|\le r
\end{subarray}}
\bigl|f(z,\ell,h)-f(z,\ell,h')\bigr|,
\end{equation}
the decomposition of~$h_y=\frakb_{D_N, x}(y)h_x+\wh h_y$ from \eqref{E:5.12i}, 
where~$\wh h_y$ is the DGFF in~$D_N\smallsetminus\{x\}$ independent of~$h_x$, along with the support restrictions on~$f$ show that, on $\{h_x^2\le M\}$,
\begin{multline}
\label{E:7.13nw}
\biggl|\E\Bigl(f\bigl(y/N,L(y),h_y\bigr)\,\Big|\,\sigma(h_x)\Bigr)-\E\Bigl(f\bigl(y/N,L(y),h_y\bigr)\Bigr)\biggr|
\\
\le 
\biggl[\text{osc}_{f,M}\Bigl(\frakb_{D_N, x}(y)\bigl[\sqrt M+(\log N)^{3/4})\bigr]\Bigr)
\,\BbbP\Bigl(|\wh h_y|\le \sqrt M+\frakb_{D_N, x}(y)(\log N)^{3/4}\Bigr)
\\+\Vert f\Vert_{\infty}\,\BbbP\bigl(|h_x|>(\log N)^{3/4}\bigr)\biggr]\,1_{\{L(y)\le M\}},
\quad
\end{multline}
where the two terms in the large square bracket arise by splitting the second expectation in the absolute value on the left according to whether the (implicit) absolute value of the DGFF at~$x$ is less than or in excess of~$(\log N)^{3/4}$. Next observe that, since~$|y-x|\ge\epsilon N$, the bound \eqref{E:5.13i} applies.  Using that  $\text{osc}_{f,M}(r)\to0$ as~$r\downarrow0$ by the uniform continuity of~$f$ and that $\wh h_y$ and~$h_x$ have variance of order~$\log N$, the right-hand side of \eqref{E:7.13nw} is at most $o((\log N)^{-1/2})1_{\{L(y)\le M\}}$ uniformly in~$y$. 

Invoking $o((\log N)^{-1/2})/ K_N =o(1/\wh W_N)$ we conclude that, for a non-random~$o(1)$ that obeys~$o(1)\to0$ as~$N\to\infty$ followed by~$\epsilon\downarrow0$, uniformly on $\{h_x^2\le M\}$,
\begin{equation}
\label{E:7.14uiq}
\E\bigl(\langle\vartheta^{D,\text{ext}}_{N,M,2\epsilon,\delta},f\rangle\,\big|\,\sigma(h_x)\bigr)\le\delta
\Vert f\Vert_{\infty}+\E\bigl(\langle\vartheta^{D,\text{ext}}_N,f\rangle\bigr)+o(1)\vartheta^D_N\bigl(\overline D\times[0,M]\bigr).
\end{equation}
Plugging this in \twoeqref{E:7.9nwt}{E:7.4uo}, invoking \eqref{E:6.2u} along with the tightness of $\{\vartheta^D_N\colon N\ge1\}$ and the uniform integrability of $\{\langle\vartheta^{D,\text{ext}}_N,f\rangle\colon N\ge1\}$ implied by \eqref{E:7.8nw} and, finally, taking~$\delta\downarrow0$ after $N\to\infty$  (and, now implicit,~$\epsilon\downarrow0$)  shows
\begin{multline}
E^\varrho\otimes\E\Bigl(\langle\vartheta^{D,\text{ext}}_N,f\rangle\,\texte^{-s\langle\vartheta^{D,\text{ext}}_N,f\rangle}\Bigr)
\\
\ge o(1)+(2\pi g)^{-1/2}
E^\varrho\Bigl(\langle\vartheta^{D}_N\otimes\leb,f\rangle\,\texte^{-s(2\pi g)^{-1/2}\langle\vartheta^{D}_N\otimes\leb,f\rangle}\Bigr),
\end{multline}
where~$o(1)\to0$  as~$N\to\infty$ uniformly in~$s\in[0,1]$.
Integrating both sides over~$s\in[0,1]$ with respect to the Lebesgue measure then gives
\begin{equation}
E^\varrho\otimes\E\bigl(\texte^{-\langle\vartheta^{D,\text{ext}}_N,f\rangle}\bigr)
\le o(1)+E^\varrho\Bigl(\texte^{-(2\pi g)^{-1/2}\langle\vartheta^{D}_N\otimes\leb,f\rangle}\Bigr).
\end{equation}
This, in combination with \eqref{E:7.3o}, proves the desired claim. 
\end{proofsect}

Next we prove an analogue of Lemma~\ref{lemma-conv}:

\begin{lemma}
Given~$f\in C_\cc(D\times[0,\infty))$ with~$f\ge0$ denote
\begin{equation}
f^{\ast\leb}(x,\ell):=\frac1{\sqrt{2\pi g}}\int_{\R}\textd h\,\,f\bigl(x,\ell+\tfrac{h^2}2\bigr).
\end{equation}
Then for every weak subsequential limit~$\vartheta^D$ of~$\vartheta^D_N$,
\begin{equation}
\label{E:6.3}
\bigl\langle\vartheta^D,f^{\ast\leb}\bigr\rangle\laweq\cspecial(\sqrt\theta)\int Z^D_{\sqrt\theta}(\textd x)\otimes\texte^{\alpha\sqrt\theta\,h}\textd h \,\,f\bigl(x,\tfrac{h^2}2\bigr)
\end{equation}
simultaneously for all~$f$ as above. 
\end{lemma}

\begin{proofsect}{Proof}
Pick~$f\in C_\cc(D\times[0,\infty))$ with~$f\ge0$ and set $f^{\text{ext}}(x,\ell,h):=f(x,\ell+\tfrac12h^2)$. Then
\begin{equation}
\begin{aligned}
\sum_{x\in D_N}f\Bigl(\ffrac xN,\frac12(\wt h^{D_N}_x+\sqrt{2t_N}\bigr)^2\Bigr)
&=\sum_{x\in D_N}f\Bigl(\ffrac xN,L^{D_N}_{t_N}(x)+\frac12(h^{D_N}_x)^2\Bigr)
\\
&=\sum_{x\in D_N}f^{\text{ext}}\Bigl(\ffrac xN, L^{D_N}_{t_N}(x),h^{D_N}_x\Bigr).
\end{aligned}
\end{equation}
Since~$f^{\text{ext}}$ is compactly supported in all variables, Lemma~\ref{lemma-7.1} tells us that, after multiplying by $\sqrt{\log N}\,/\wh W_N $ and specializing~$N$ to the  subsequence along which~$\vartheta^D_N$ tends in law to~$\vartheta^D$, the right-hand side tends to $\bigl\langle\vartheta^D,f^{\ast\leb}\bigr\rangle$.  By \eqref{E:1.19} and the fact that~$\sqrt{2t_N}\sim2\sqrt g\sqrt\theta\log N$, the left-hand side tends to the measure on the right of \eqref{E:6.3}.
\end{proofsect}

With these in hand we are ready to prove convergence of~$\vartheta^D_N$'s:

\begin{proofsect}{Proof of Theorem~\ref{thm-light}}
Pick~$A\subseteq   D$ open with $\overline A\subseteq   D$. Taking a sequence of compactly supported functions converging  upward to $f(x,h):=1_A(x)\texte^{-s h}1_{[0,\infty)}(h)$, where~$s>0$, and denoting
\begin{equation}
\tilde\mu_A(B):=\vartheta^D(A\times B),
\end{equation}
the Tonelli and Monotone Convergence Theorems yield
\begin{equation}
\label{E:7.18iw}
\int_0^\infty\tilde\mu_A(\textd\ell)\texte^{-\ell s} \laweq \sqrt{2\pi g}\,\cspecial(\sqrt\theta)Z^D_{\sqrt\theta}(A)\,\texte^{\frac{\alpha^2\theta}{2s}},\quad s>0.
\end{equation}
Note that $s\mapsto \texte^{\frac{\alpha^2\theta}{2s}}$ is the Laplace transform of the measure in \eqref{E:1.33}. Since the Laplace transform determines Borel measures on $[0,\infty)$ uniquely, the claim follows by the fact that the right-hand side is a Borel measure in~$A$ which is determined by its values on~$A$ open with the closure in~$D$. 
\end{proofsect}

In order to extend Theorem~\ref{thm-light} to the control of the measure $\kappa^D_N$ associated with the avoided points, we need the following estimate:

\begin{lemma}
\label{lemma-6.3}
Let~$A\subseteq   D$ be open with~$\overline A\subseteq   D$. Then
\begin{equation}
\lim_{\epsilon\downarrow0}\,\limsup_{N\to\infty}\,\frac{N^2}{\wh W_N }\max_{\begin{subarray}{c}
x\in D_N\\ x/N\in A
\end{subarray}}
P^\varrho\Bigl(0<L^{D_N}_{t_N}(x)\le\epsilon\Bigr)=0.
\end{equation}
\end{lemma}

\begin{proofsect}{Proof}
First note that, using Dynkin's isomorphism, we get
\begin{equation}
\label{E:monster}
\begin{aligned}
P^\varrho\Bigl(0&<L^{D_N}_{t_N}(x)\le\epsilon\Bigr)P\Bigl(\frac12 (h^{D_N}_x)^2\le\epsilon\Bigr)
\\
&\le P^\varrho\otimes\BbbP\Bigl(L^{D_N}_{t_N}(x)+\frac12 (h^{D_N}_x)^2\le 2\epsilon,\,L^{D_N}_{t_N}(x)>0\Bigr)
\\
&=\BbbP\Bigl(\frac12(\wt h^{D_N}_x-\sqrt{2t_N})^2\le 2\epsilon\Bigr)-\BbbP\Bigl(\frac12(h^{D_N}_x)^2\le 2\epsilon\Bigr)P^\varrho\bigl(L^{D_N}_{t_N}(x)=0\bigr).
\end{aligned}
\end{equation}
The fact that $|G^{D_N}(x,x)- g\log N|$ is bounded uniformly for all~$x\in D_N$ with~$x/N\in A$ then shows
\begin{equation}
\BbbP\Bigl(\frac12(\wt h^{D_N}_x-\sqrt{2t_N})^2\le \epsilon\Bigr)
=\bigl(2+o(1)\bigr)\sqrt{2\epsilon}\,\frac1{\sqrt{2\pi G^{D_N}(x,x)}}\,
\,\texte^{-\frac{t_N}{G^{D_N}(x,x)}}
\end{equation}
while
\begin{equation}
P\Bigl(\frac12 (h^{D_N}_x)^2\le\epsilon\Bigr)= 
\bigl(2+o(1)\bigr)\sqrt{2\epsilon}\,\frac1{\sqrt{2\pi G^{D_N}(x,x)}},
\end{equation}
where $o(1)\to0$ as~$\epsilon\downarrow0$ uniformly in~$N\ge1$ and~$x$ as above. In light of \eqref{E:L-vanish}, the right-hand side of \eqref{E:monster} divided by $\wh W_N /N^2$-times the DGFF probability on the extreme left tends to zero as~$N\to\infty$ and~$\epsilon\downarrow0$. 
\end{proofsect}

We are ready to give:

\begin{proofsect}{Proof of Theorem~\ref{thm-avoid}}
Take~$f_n\in C_\cc([0,\infty))$ such that~$f_n(h):=(1-nh)\vee0$ and pick~$A\subseteq   D$ open with $\overline A\subseteq   D$. Then
\begin{equation}
 E^\varrho  \bigl|\langle\kappa^D_N,1_A\rangle-\langle\vartheta^D_N,1_A\otimes f_n\rangle\bigr|
\le\frac2{\wh W_N }\sum_{\begin{subarray}{c}
x\in D_N\\ x/N\in A
\end{subarray}}P^\varrho\bigl(0<L^{D_N}_{t_N}(x)\le1/n\Bigr)\,.
\end{equation}
By Lemma~\ref{lemma-6.3}, the sum on the right-hand side tends to zero in the limits~$N\to\infty$ followed by~$n\to\infty$. Theorem~\ref{thm-light} in turn shows that
\begin{equation}
\langle\vartheta^D_N,1_A\otimes f_n\rangle\,\,\underset{N\to\infty}\Lawarrow\,\,
\sqrt{2\pi g}\,\cspecial(\sqrt\theta)Z^D_{\sqrt\theta}(A)\Bigl[1+\int_{(0,1/n]} \mu(\textd h) f_n(h)\Bigr],
\end{equation}
where~$\mu$ is the measure in \eqref{E:1.33}. The claim follows by noting that the integral on the right tends to zero as~$n\to\infty$.
\end{proofsect}

\section{Local structure}
\label{sec8}\noindent
In this section we deal with local structures of the exceptional level sets associated with the local time~$L^{D_N}_{t_N}$. Throughout we again rely on the coupling of $L_{t_N}^{D_N}$ and an independent DGFF $h^{D_N}$ 
to another DGFF~$\wt h^{D_N}$ via the Dynkin isomorphism (Theorem~\ref{thm-Dynkin}). We start with the thick points.

\subsection{Local structure of thick points}
Let $a_N$ and $t_N$ satisfy \eqref{E:1.20} with some $\theta > 0$ and some $\lambda \in (0,1)$ and recall the notation~$\wh \zeta^D_N$ for the extended point measures from \eqref{E:zeta-local} that describe the $\lambda$-thick points along with their local structure. Let $\wh a_N$ be the sequence given by (\ref{E:4.1}).
We will compare~$\wh \zeta_N^D$ to the point measures
\begin{equation}
\label{E:full-etaDGFF}
\wh \eta^D_N :=
\frac{1}{W_N} \sum_{x \in D_N} \delta_{x/N} \otimes 
\delta_{\wt h_x^{D_N} - \wh a_N} 
\otimes \delta_{\left\{\frac{2\sqrt{2 a_N}+2(\wt h_x^{D_N} - \wh a_N) + (\wt h_{x+z}^{D_N} - \wt h_x^{D_N})}{2\log N} (\wt h_x^{D_N} - \wt h_{x+z}^{D_N})~:~z \in \mathbb{Z}^2 \right\}}
\end{equation} 
associated with the DGFF~$\wt h^{D_N}$. For that we need:

\begin{lemma}[Gradients of squared DGFF]
\label{lemma-gradient-DGFF}
For all $b \in\R$, all $M \ge1$  and all $r > 0$, 
\begin{multline}
\label{E:9.3}
\lim_{N \to \infty} \frac{1}{W_N}\sum_{x \in D_N} P^{\varrho} 
\bigl(L_{t_N}^{D_N} (x) \geq a_N + b \log N \bigr) \\
\times \BbbP \biggl(\,\bigcup_{ z \in \Lambda_r (0)}\Bigl\{
\bigl|(h_x^{D_N})^2 - (h_{x+z}^{D_N})^2 \bigr| > (\log N)^{3/4},
\,|h_x^{D_N}| \leq M \sqrt{\log N}\Bigr\}\biggr)
= 0,
\end{multline}
where $\Lambda_r (x) := \{z \in \Bbb{Z}^2 \colon |z-x| \leq r \}$.
\end{lemma}

\begin{proofsect}{Proof}
When $|h_x^{D_N}| \leq M \sqrt{\log N}$,
we have
\begin{equation}
\bigl|(h_x^{D_N})^2 - (h_{x+z}^{D_N})^2 \bigr| \leq
 \bigl|h_x^{D_N} - h_{x+z}^{D_N} \bigr|^2 + 2M \sqrt{\log N} \bigl|h_x^{D_N} - h_{x+z}^{D_N} \bigr|.
\end{equation}
Thus, for $M\ge1$, the term corresponding to $x \in D_N$ 
on the left-hand side of (\ref{E:9.3}) is bounded from above by
\begin{equation}
\label{E:8.4nw}
\sum_{z \in \Lambda_r (0)}
P^{\varrho}
\bigl(L_{t_N}^{D_N} (x) \geq a_N + b \log N \bigr) 
\mathbb{P} \Bigl(\bigl|h_x^{D_N} - h_{x+z}^{D_N} \bigr| > (4M)^{-1}(\log N)^{1/4}  \Bigr).
\end{equation}
For $\epsilon > 0$, abbreviate $D_N^{\epsilon} := \{x \in D_N \,\colon\,\dinfty(x, D_N^{\text{c}}) > \epsilon N \}$. 
Then for any $x \in D_N^{\epsilon}$ and $z \in \Lambda_r (0)$,
$\text{Var}_{\BbbP} (h_x^{D_N} - h_{x+z}^{D_N})$ is equal to
\begin{multline}
\qquad
G^{D_N} (x, x) + G^{D_N} (x+z,x+z) - 2G^{D_N} (x,x+z) 
\\= g\log N + g \log N - 2 g \log (N/(1+|z|)) + O(1) 
\\\leq 2g \log (1+r) + O(1).
\qquad
\end{multline}
The standard Gaussian tail estimate bounds \eqref{E:8.4nw} by $o(1)P^\varrho(L_{t_N}^{D_N} (x) \geq a_N + b \log N)$ with $o(1)\to0$ uniformly in~$x\in D_N^\epsilon$. Lemma~\ref{lemma-upper} subsequently shows that the sum over $x \in D_N^{\epsilon}$ on the left-hand side of \eqref{E:9.3} is $o(1)$ as $N \to \infty$.
The sum over $x \in D_N \smallsetminus D_N^{\epsilon}$ is bounded from above by
$E^\varrho(\zeta_N^D (D \smallsetminus D^{\epsilon} \times [b, \infty) ))$ which tends to $0$ 
as $N \to \infty$ followed by $\epsilon \downarrow 0$ by Corollary~\ref{cor-tightness-upper}. 
\end{proofsect}

We are ready to give:

\begin{proofsect}{Proof of Theorem~\ref{thm-thick-local}}
Pick any $f = f(x, \ell, \phi) \in C_{\text{c}} (D \times \R \times \R^{\mathbb{Z}^2})$
which depends only on a finite number of coordinates of $\phi$, say, those in $\Lambda_r (0)$
for some $r > 0$. The following identity is key for the entire proof
\begin{multline}
\qquad
\left\{\sqrt{2 a_N}+(\wt h_x^{D_N} - \wh a_N) + \frac{1}{2}(\wt h_{x+z}^{D_N} - \wt h_x^{D_N}) \right\}
(\wt h_x^{D_N} - \wt h_{x+z}^{D_N})
\\
= L_{t_N}^{D_N} (x) - L_{t_N}^{D_N} (x+z) + \frac{1}{2} (h_x^{D_N})^2 - \frac{1}{2} (h_{x+z}^{D_N})^2.
\qquad
\end{multline} 
Indeed,  writing $\nabla_z s (x) := s (x) - s (x+z)$ for a version of the discrete gradient of $s\colon \Bbb{Z}^2 \to \R$,  we then get
\begin{multline}
\qquad
\label{E:9.10}
\langle \wh \eta_N^D, f \rangle
= o(1)+\frac{1}{W_N}
\sum_{x \in D_N} f \Biggl(\frac{x}{N}, \sqrt{2 L_{t_N}^{D_N} (x) + (h_x^{D_N})^2} - \sqrt{2a_N}, 
\\
\biggl\{\frac{\nabla_z L_{t_N}^{D_N} (x)}{\log N}
+ \frac{\nabla_z (h^{D_N})^2 (x)}{2 \log N}\colon z \in \mathbb{Z}^2 \biggr\} \Biggr),
\qquad
\end{multline}
 where $o(1)$ stands for the analogue of the second term on the extreme right of \eqref{E:4.26}; this term tends to zero in probability as~$N\to\infty$ by exactly the same argument. 

In order to control  the  gradients of the DGFF squared that appear on the right-hand side of \eqref{E:9.10}, set
\begin{equation}
\label{E:gradient-squared-DGFF}
G_{N, r} (x) := \bigcap_{z \in \Lambda_r (0)}\Bigl\{\bigl|\nabla_z (h^{D_N})^2 (x)\bigr| \leq (\log N)^{3/4}
\Bigr\}
\end{equation}
and let, as before, $\chi \colon [0, \infty) \to [0,1]$ be a non-increasing, continuous function
with $\chi(x) = 1$ for $0 \leq x \leq 1$ and $\chi(x) = 0$ for $x \geq 2$.
By Lemmas \ref{lemma-no-conspire} and \ref{lemma-gradient-DGFF},
we may truncate \eqref{E:9.10} by introducing $1_{G_{N, r} (x)}$ and $\chi (M^{-1}|h_x^{D_N}|/ \sqrt{\log N})$
for $M > 0$ under the sum and write $\langle \wh \eta_N^D, f \rangle$ as a random quantity whose $L^1$ norm is at most a constant times  $\|f\|_{\infty} \texte^{- \beta M^2}$
uniformly in~$N$ plus the quantity
\begin{multline}
\label{E:9.12}
\quad
\frac{1}{W_N}
\sum_{x \in D_N} 1_{G_{N,r} (x)} f \Biggl(\frac xN, \sqrt{2 L_{t_N}^{D_N} (x) + (h_x^{D_N})^2} - \sqrt{2a_N}, \\
\biggl\{\frac{\nabla_z L_{t_N}^{D_N} (x)}{\log N}
+ \frac{\nabla_z (h^{D_N})^2 (x)}{2 \log N}\colon z \in \mathbb{Z}^2 \biggr\} \Biggr)
 \chi \biggl(\frac{|h_x^{D_N}|}{M \sqrt{\log N}} \biggr).
 \quad
\end{multline}
Using the uniform continuity of $f$
and Corollary~\ref{cor-tightness-upper} and
Lemma~\ref{lemma-gradient-DGFF},
we rewrite (\ref{E:9.12}) by a random quantity which tends to $0$
as $N \to \infty$ in probability plus the quantity
\begin{equation}
\label{E:9.13}
\frac{1}{W_N} \sum_{x \in D_N} f_{\text{ext}} \Biggl(\frac xN, \frac{L_{t_N}^{D_N} (x) - a_N}{\log N},
\Bigl\{\frac{\nabla_z L_{t_N}^{D_N} (x)}{\log N} \Bigr\}_{z \in \mathbb{Z}^2},
\frac{h_x^{D_N}}{\sqrt{\log N}} \Biggr) \chi \biggl(\frac{|h_x^{D_N}|}{M \sqrt{\log N}} \biggr),
\end{equation}
where we introduced
\begin{equation}
f_{\text{ext}} (x, \ell, \phi, h) := f 
\biggl(x, \frac{1}{2 \sqrt{g} (\sqrt{\theta} + \lambda)}
\bigl(\ell + \tfrac{1}{2} h^2 \bigr), \phi \biggr).
\end{equation}
Note that Corollary~\ref{cor-tightness-upper} implies that $\{\wh \zeta_N^D\colon N\ge1\}$ is tight. 
Let $\wh \zeta^D$ be a subsequential weak limit of $\wh \zeta_N^D$ along the subsequence $\{N_k\}$.
By the same argument as in the proof of Lemma~\ref{lemma-add-field},
as $k \to \infty$ followed by $M \to \infty$, 
$\langle \wh \eta_{N_k}^D, f \rangle$ converges in law to
\begin{equation}
\label{E:8.13nw}
\int \wh \zeta^D (\textd x\,\textd\ell\,\textd\phi) \otimes \mathfrak{g} (\textd h)
f_{\text{ext}}(x, \ell, \phi, h).
\end{equation}
On the other hand, noting that~$\sqrt{2a_N}/\log N\to 2\sqrt g(\sqrt\theta+\lambda)$, \cite[Theorem 2.1]{BL4} shows that $\langle \wh \eta_N^D, f \rangle$ converges, as $N \to \infty$, in law to
\begin{equation}
  \int \cspecial(\lambda)Z_{\lambda}^D (\textd x) \otimes \texte^{- \alpha \lambda h} \textd h \otimes 
\nu_{\theta, \lambda} (\textd \phi) f(x, h, \phi).
\end{equation}
The arguments in the proof of Theorem~\ref{thm-thick} show that the class of functions $f_{\text{ext}}$ arising from~$f\in C_\cc(D\times\R\times\R^{\Z^2})$ above determines the measure~$\wh\zeta^D$ uniquely from \eqref{E:8.13nw}; the calculation \eqref{E:4.38} then gives
\begin{equation}
\wh \zeta^D \,\laweq\,
\frac{\theta^{1/4}}{2\sqrt{g} (\sqrt{\theta} + \lambda)^{3/2}} \cspecial(\lambda) 
Z_{\lambda}^D (\textd x) \otimes \texte^{- \alpha (\theta, \lambda) \ell} \textd\ell \otimes 
\nu_{\theta, \lambda} (\textd\phi).
\end{equation}
This is the desired claim.
\end{proofsect}

\subsection{Local structure of thin points}
We move to the proof of the convergence of point measures $\wh \zeta_N^D$ 
associated with $\lambda$-thin points. The proof follows very much the same steps as for the thick points so we stay quite brief. 
Assume that~$a_N$ and~$t_N$ satisfy \eqref{E:1.22} with some $\theta > 0$ and some $\lambda \in (0, 1 \wedge \sqrt{\theta} )$.
As a counterpart to Lemma~\ref{lemma-gradient-DGFF}, we need the following:

\begin{lemma}[Gradients of squared DGFF]
\label{lemma-control-gradient-DGFF2}
For all $b > 0$, all $M \ge1$ and all $r > 0$,
\begin{multline}
\lim_{N \to \infty}
\frac{1}{W_N}\sum_{x \in D_N} 
P^{\varrho} \Bigl(a_N - b \log N \leq L_{t_N}^{D_N} (x) \leq a_N + b \log N \Bigr) \\
\times 
\BbbP \biggl(\,\bigcup_{ z \in \Lambda_r (0)}\Bigl\{
\bigl|(h_x^{D_N})^2 - (h_{x+z}^{D_N})^2 \bigr| > (\log N)^{3/4},
\,|h_x^{D_N}| \leq M \sqrt{\log N}\Bigr\}\biggr)
= 0.
\end{multline}
\end{lemma}

\begin{proofsect}{Proof}
The proof is the same as that of Lemma~\ref{lemma-gradient-DGFF}
except that we use Lemma~\ref{lemma-lower} 
and Corollary~\ref{cor-tightness-lower}
instead of Lemma~\ref{lemma-upper} and Corollary~\ref{cor-tightness-upper}, respectively.
\end{proofsect}

We are again ready to start:

\begin{proofsect}{Proof of Theorem~\ref{thm-thin-local}}
Set 
\begin{equation}
\wh a_N := \sqrt{2t_N} - \sqrt{2a_N}
\end{equation}
and pick any
$f = f(x, \ell, \phi) \in C_{\text{c}} (D \times \R \times \R^{\mathbb{Z}^2})$
that depends only on a finite number of coordinates of~$\phi$.
Let $\wh \eta_N^D$ be the point process obtained from (\ref{E:full-etaDGFF})
by replacing~$\wh a_N$ by~$- \wh a_N$. Using the calculation
\begin{multline}
\qquad
\left\{\sqrt{2 a_N}+(\wt h_x^{D_N} + \wh a_N) + \frac{1}{2}(\wt h_{x+z}^{D_N} - \wt h_x^{D_N}) \right\}
(\wt h_x^{D_N} - \wt h_{x+z}^{D_N}) 
\\
= L_{t_N}^{D_N} (x) - L_{t_N}^{D_N} (x+z) + \frac{1}{2} (h_x^{D_N})^2 - \frac{1}{2} (h_{x+z}^{D_N})^2
\qquad
\end{multline}
we then again have \eqref{E:9.10} for $\langle\wh\eta^D_N,f\rangle$.
Using Corollary~\ref{cor-tightness-lower} and
Lemmas \ref{lemma-no-conspire2} and \ref{lemma-control-gradient-DGFF2},
we rewrite \eqref{E:9.10} as a random quantity whose $L^1$ norm 
is at most a constant times $\|f\|_{\infty} \texte^{- \beta M^2}$
uniformly in~$N$ plus \eqref{E:9.13},
where, in this case,
\begin{equation}
f_{\text{ext}} (x, \ell, \phi, h) := f \biggl(x, \frac{1}{2 \sqrt{g} (\sqrt{\theta} - \lambda)}
\bigl(\ell + \tfrac{1}{2} h^2 \bigr), \phi \biggr).
\end{equation}
Note that Corollary~\ref{cor-tightness-lower} implies the tightness of $\{\wh \zeta_N^D\colon N\ge1\}$.
Let $\wh \zeta^D$ be any subsequential weak limit of $\wh \zeta_N^D$ along the subsequence~$\{N_k\}$. 
By the same argument as in the proof of Lemma~\ref{lemma-add-field2},
as $k \to \infty$ and $M \to \infty$, $\langle \wh \eta_{N_k}^D, f \rangle$ tends in law to
\begin{equation}
\int \wh \zeta^D (\textd x\,\textd\ell\,\textd\phi) \otimes \mathfrak{g} (\textd h)
f_{\text{ext}}(x, \ell, \phi, h).
\end{equation}
On the other hand, by \cite[Theorem~2.1]{BL4}, as $N \to \infty$,
$\langle \wh \eta_N^D, f \rangle$ converges in law to
\begin{equation}
 \int  \cspecial(\lambda)Z_{\lambda}^D (\textd x) \otimes \texte^{\alpha \lambda h} \textd h \otimes 
\wt \nu_{\theta, \lambda} (\textd\phi) f(x, h, \phi).
\end{equation}
The arguments in the proof of Theorem~\ref{thm-thin} and the calculation \eqref{E:4.38a} then show
\begin{equation}
\wh \zeta^D \,\laweq\,
\frac{\theta^{1/4}}{2\sqrt{g} (\sqrt{\theta} - \lambda)^{3/2}} \,\cspecial(\lambda) \,
Z_{\lambda}^D (\textd x) \otimes \texte^{\wt \alpha (\theta, \lambda) \ell} \textd\ell \otimes 
\wt \nu_{\theta, \lambda} (\textd\phi).
\end{equation}
This is the desired claim.
\end{proofsect}

\subsection{Local structure of avoided points}
In this section we will prove the convergence of the point measures 
associated with the local structure of the avoided points. The proof will make use of the Pinned Isomorphism Theorem~(see Theorem~\ref{thm-pinned-iso}) but that so only at the very end. Most of the argument consists of careful manipulations with the doubly extended measure
\begin{equation}
\wh \kappa_N^{D, \text{ext}} :=
\frac{\sqrt{\log N}}{\wh W_N}
\sum_{x \in D_N}
1_{\{L_{t_N}^{D_N} (x) = 0 \}} \delta_{x/N} \otimes \delta_{\{L_{t_N}^{D_N} (x+z) \,:\,z \in \Bbb{Z}^2 \}}
\otimes \delta_{h_x^{D_N}} \otimes \delta_{\{\wh h_{x+z}^{D_N \smallsetminus \{x\}}\,:\,z \in \Bbb{Z}^2 \}},
\end{equation}
where, for $\frakb_{D_N,x}$ as in \eqref{E:5.12i},
\begin{equation}
\label{E:8.24uiw}
\wh h^{D_N \smallsetminus \{x \}}_z :=h^{D_N}_z-h^{D_N}_x\frakb_{D_N,x}(z),\quad z\in\Z^2.
\end{equation}
By \eqref{E:5.12i}, $\wh h^{D_N \smallsetminus \{x \}}$ is the field $h^{D_N}$ conditioned on $h^{D_N}_x=0$. In particular, 
\begin{equation}
\label{E:8.25uiw}
\wh h^{D_N \smallsetminus \{x \}}\independent h^{D_N}_x.
\end{equation}
Corollary~\ref{cor-light} implies that $\{\wh \kappa_N^{D, \text{ext}} \colon N\ge1\}$ is tight with respect to vague convergence of measures on the product space $D\times[0,\infty)^{\Z^2}\times\R\times\R^{\Z^2}$.  As before, a key ingredient we need is factorization of the subsequential limits: 

\begin{lemma}
\label{lemma-avoid-local-add}
Suppose $\{N_k\}$ is a subsequence along which $\wh \kappa_N^D$
converges in law to $\wh \kappa^D$. 
Then
\begin{equation}
\wh \kappa_N^{D, \text{\rm ext}} \,\,\,\underset{\begin{subarray}{c} N = N_k \\ k\to\infty \end{subarray}}\Lawarrow\,\,\,
 \frac{1}{\sqrt{2\pi g}} \,\wh \kappa^D \otimes \leb \otimes \nu^0.
\end{equation}
\end{lemma}

\begin{proofsect}{Proof}
Let $f = f(x, \ell, h, \phi)\colon D\times[0,\infty)^{\Z^2}\times\R\times\R^{\Z^2} \to \R$
be a continuous, compactly-supported function
that depends only on a finite number of coordinates of $\ell$ and $\phi$, say, 
those in $\Lambda_{r_0} (0)$ for some~$r_0>0$.
Suppose in addition that $f(x, \ell, h, \phi) = 0$
unless $x \in A$ for some open $A \subseteq   D$ with $\ol A \subseteq   D$, and unless
$|h|^2 \leq M$ and $\ell_z, |\phi_z| \leq M$ for all $z \in \Lambda_{r_0} (0)$ for some $M > 0$.

Noting that only the second pair of the variables of~$\wh \kappa_N^{D, \text{ext}}$ is affected by expectation~$\E$ with respect to the law of~$h^{D_N}$, we now claim
\begin{equation}
\label{E:(8.26)}
\Bbb{E} \langle \wh \kappa_N^{D, \text{ext}}, f \rangle
= \frac{1}{\sqrt{2\pi g}}\,\langle \wh \kappa_N^D \otimes \leb \otimes \nu^0, f \rangle + o(1),
\end{equation}
where $o(1) \to 0$ in $P^{\varrho}$-probability as $N \to \infty$ and where~$\nu^0$ the law of the pinned DGFF. As in the proof of Lemma~\ref{lemma-7.1},  \eqref{E:(8.26)}  follows by noting that the probability density of~$h_x^{D_N}$ multiplied by~$\sqrt{\log N}$ tends to $(2\pi g)^{-1/2}$ uniformly over any compact interval and by the fact $(\wh h_{x+z}^{D_N \smallsetminus \{x\}})_{z \in \Lambda_{r_0} (0)}$ tends in law to $(\phi_z)_{z \in \Lambda_{r_0}(0)}$ (which can be gleaned from the representation of the Green function by the potential kernel, see \cite[Lemma~B.3]{BL3}, and the asymptotic expression for the potential kernel, see \cite[Lemma~B.4]{BL3}). These two convergences may be applied jointly in light of the independence \eqref{E:8.25uiw} and the Bounded Convergence Theorem enabled by the tightness of~$\{\wh\kappa^D_N\colon N\ge1\}$.

In order to convert  the   convergence in the mean to the convergence in law,  we proceed as in the proof of Lemma~\ref{lemma-7.1}.
 Let us abbreviate~$L^{D_N}_{t_N}(x)$, $h^{D_N}_x$ and $\wh h_{x+z}^{D_N \smallsetminus \{x \}}$
by $L(x)$,~$h_x$ and~$\phi_z^{(x)}$, respectively, for the duration of this proof.  
Recall the event $F_{N, M, \epsilon, \delta} (x)$ in the proof of Lemma~\ref{lemma-7.1}.
By the argument leading up to \eqref{E:7.9nwt},
for the truncated measure 
\begin{equation}
\wh \kappa^{D,\text{ext}}_{N,M,\epsilon,\delta}:=\frac1{ K_N }
\sum_{x\in D_N}1_{F_{N,M,\epsilon,\delta}(x)}\,1_{\{L (x) = 0\}}\,
\delta_{x/N}\otimes\delta_{\{L(x+z)\colon z \in \mathbb{Z}^2 \}}\otimes\delta_{h_x}
\otimes\delta_{\{\phi_z^{(x)} : z \in \mathbb{Z}^2 \}},
\end{equation}
 where $K_N$ abbreviates $\wh W_N/\sqrt{\log N}$   we get
\begin{multline}
\quad\lim_{\epsilon\downarrow0}\,\limsup_{N\to\infty}\biggl|E^\varrho\otimes\E\Bigl(\langle\wh \kappa^{D,\text{ext}}_N,f\rangle\,\texte^{-s\langle\wh \kappa^{D,\text{ext}}_N,f\rangle}\Bigr)
\\
-E^\varrho\otimes\E\Bigl(\langle\wh \kappa^{D,\text{ext}}_N,f\rangle\,\texte^{-s\langle\wh \kappa^{D,\text{ext}}_{N,M,2\epsilon,\delta},f\rangle}\Bigr)\biggr|=0
\quad
\end{multline}
uniformly in~$s\in[0,1]$. 
 Focusing attention on the second expectation and writing $\mathcal{G}_{r_0} (x)$ for the $\sigma$-field generated by $\{h_{x+z}\colon z \in \Lambda_{r_0} (0)\}$,   the conditional Jensen inequality shows
\begin{equation}
\label{E:(8.29)}
\begin{aligned}
&E^\varrho\otimes\E\Bigl(\langle\wh \kappa^{D,\text{ext}}_N,f\rangle\,\texte^{-s\langle\wh \kappa^{D,\text{ext}}_{N,M,2\epsilon,\delta},f\rangle}\Bigr)
\\
&\ge
\frac1{ K_N }\sum_{\begin{subarray}{c}
x\in D_N\\x/N\in A
\end{subarray}}
E^\varrho\otimes\E\Bigl(1_{\{L (x) = 0 \}} f\bigl(x/N, L(x+\cdot),h_x, \phi^{(x)}\bigr)
\texte^{-s\E(\langle\wh \kappa^{D,\text{ext}}_{N,M,2\epsilon,\delta},f\rangle|\mathcal{G}_{r_0} (x))}\Bigr).
\end{aligned}
\end{equation}
The conditional expectation in the exponent is bounded by
\begin{multline}
\label{E:(8.30)}
\quad
\E\bigl(\langle \wh \kappa^{D,\text{ext}}_{N,M,2\epsilon,\delta},f\rangle\,\big|\,\mathcal{G}_{r_0} (x)\bigr)
\\\le
\frac1{ K_N }\sum_{\begin{subarray}{c}
y\in D_N\\|x-y|\ge\epsilon N
\end{subarray}}1_{\{L (y) = 0\}}\,\E\Bigl(f\bigl(y/N,L(y+\cdot),h_y, \phi^{(y)}\bigr)\,\Big|\,\mathcal{G}_{r_0} (x)\Bigr)
\\+
\frac{\Vert f\Vert_{\infty}}{ K_N }\,\E\biggl(\,\sum_{\begin{subarray}{c}
y\in D_N\\|x-y|<\epsilon N
\end{subarray}}1_{\{L(y)+\frac12 h_y^2\le 2M\}}1_{F_{N,M,2\epsilon,\delta}(y)}\,\Big|\,\mathcal{G}_{r_0} (x)\Bigr).
\end{multline}
As in \eqref{E:7.5iue}, the second term on the right 
is bounded by~$\delta\Vert f\Vert_{\infty}$ pointwise. 

Concerning the first term on the right of \eqref{E:(8.30)}, we consider the analogue of the quantity $\text{osc}_{f,M}(r)$ in \eqref{E:7.12uiq} defined, for any~$r>0$, by
\begin{equation}
\sup_{z\in D}\,\,\sup_{\ell\in [0, M]^{\Lambda_{r_0} (0)}}\,\,\sup_{\begin{subarray}{c}
h,h'\in[-\sqrt{M},\sqrt{M}]\\|h-h'|\le r
\end{subarray}}
\sup_{\begin{subarray}{c} \phi, \phi' \in [-M, M]^{\Lambda_{r_0} (0)}\\|\phi_z - \phi'_z|\le r, \, \forall z \in \Lambda_{r_0} (0) \end{subarray}}
\bigl|f(z,\ell,h,\phi)-f(z,\ell,h',\phi')\bigr|.
\end{equation}
Consider 
the decomposition of~$h_y=\sum_{z \in \Lambda_{r_0} (0)}\frakb_z^{(x)}(y)h_{x+z}+h_y^{x, r_0}$ , where~
$\frakb_z^{(x)}(y) := P^y (H_{\Lambda_{r_0} (x)} < H_{\varrho}, X_{H_{\Lambda_{r_0} (x)}} = x+z)$ and
$h_y^{x, r_0}$ is the DGFF in~$D_N\smallsetminus \Lambda_{r_0} (x)$ independent of~$h_{x+z}$, $z \in \Lambda_{r_0} (0)$.
On the event $\{h_x^2\le M\} \cap \bigcap_{z \in \Lambda_{r_0} (0)} \{|\phi_z^{(x)}| \le M \}$,
we have
\begin{multline}
\biggl|\E\Bigl(f\bigl(y/N,L(y+\cdot),h_y,\phi^{(y)}\bigr)\,\Big|\,\mathcal{G}_{r_0}(x)\Bigr)-\E\Bigl(f\bigl(y/N,L(y+\cdot),h_y,\phi^{(y)}\bigr)\Bigr)\biggr|
\\
\le 
\text{osc}_{f,M}\Bigl(2\frakb^{(x)} (y)\bigl[2 M+(\log N)^{3/4})\bigr]\Bigr)
\,\BbbP\Bigl(|h_y^{x,r_0}|\le \sqrt M+\frakb^{(x)}(y)(\log N)^{3/4}\Bigr)
\\+\Vert f\Vert_{\infty}\,\sum_{z \in \Lambda_{r_0}(0)}\BbbP\bigl(|h_{x+z}|>(\log N)^{3/4}\bigr),
\quad
\end{multline}
where $\frakb^{(x)} (y) := \max_{z \in \Lambda_{r_0} (0)} P^{y+z} [H_{\Lambda_{r_0} (x)} < H_{\rho}]$.
Since~$|y-x|\ge\epsilon N$, the bound \eqref{E:5.13i} dominates
$\frakb^{(x)} (y)$ by $c(\log N)^{-1}$, where $c > 0$ depends on $\epsilon$ and $r_0$.

Using these observations (as in \eqref{E:7.14uiq}), the conditional expectation on the right of~\eqref{E:(8.29)} is at most $\E(\langle\wh \kappa^{D,\text{ext}}_N,f\rangle)+\delta\Vert f\Vert_\infty +o(1)$ where~$o(1)\to0$ in probability as~$N\to\infty$. The rest of the proof of Lemma~\ref{lemma-7.1} then applies to give the desired claim. 
\end{proofsect}

We are now ready to give:

\begin{proofsect}{Proof of Theorem~\ref{thm-avoid-local}}
Consider the coupling from Theorem~\ref{thm-Dynkin} between the local time~$L^{D_N}_{t_N}$ and two copies~$h^{D_N}$ and~$\wt h^{D_N}$ of the DGFF in~$D_N$, with the former independent of~$L^{D_N}_{t_N}$. Recall the definition of $\wh h^{D_N\smallsetminus\{x\}}$ from \eqref{E:8.24uiw}, write $\phi^{(x)}_z:=\wh h^{D_N\smallsetminus\{x\}}_{x+z}$ and abbreviate $\nabla_z s(x):= s(x)-s(x+z) $. Then for each $x \in D_N$ and $z \in \Bbb{Z}^2$, we have
\begin{multline}
\label{E:9.37}
\qquad
\Bigl(\wt h_x^{D_N} + \sqrt{2t_N} - \frac{1}{2} \nabla_z \wt h^{D_N} (x) \Bigr)
\Bigl(- \nabla_z \wt h^{D_N} (x) \Bigr) \\
= - \nabla_z L_{t_N}^{D_N} (x) 
+ \frac{1}{2} \bigl(\,\phi^{(x)}_z
+ \frakb_{D_N,x} (x+z) h_x^{D_N}\bigr)^2
- \frac{1}{2} (h_x^{D_N})^2.
\qquad
\end{multline}
Let $\Phi_x (z)$ and $\Psi_x (z)$ denote the left-hand side and the right-hand side
of \eqref{E:9.37}, respectively.
Then for each 
$f \colon D \times [0, \infty) \times \R^{\Bbb{Z}^2}\to \R$,
\begin{multline}
\label{E:9.38}
\quad
\frac{\sqrt{\log N}}{\wh W_N}
\sum_{x \in D_N}
f \left(x/N,\, \frac{1}{2} (\wt h_x^{D_N} + \sqrt{2t_N})^2,\,
\{\Phi_x (z)\,\colon\,z \in \Bbb{Z}^2\} \right) 
\\
= \frac{\sqrt{\log N}}{\wh W_N}
\sum_{x \in D_N}
f \left(x/N,\, L_{t_N}^{D_N} (x) + \frac{1}{2} (h_x^{D_N})^2,\,
\{\Psi_x (z)\,\colon\,z \in \Bbb{Z}^2\} \right).
\quad
\end{multline}
Next pick $F\colon \R^{\Bbb{Z}^2}\to\R$ that is continuous and
depends only on a finite number of coordinates, say, in~$\Lambda_r (0)$,
and obeys $F(\phi) = 0$ unless $|\phi_z| \leq M$ for all $z \in \Lambda_r (0)$ for some $M > 0$. Then set $f(x, \ell, \phi) := 1_A(x) f_n (\ell) F(\phi)$,
where $A \subseteq   D$ is an open set with $\ol{A} \subseteq   D$ and $f_n\colon [0, \infty) \to [0, 1]$
are given by $f_n (\ell) := (1 - n \ell) \vee 0$. The Bounded Convergence Theorem ensures that \eqref{E:9.38} applies to these~$f$'s as well so we will now explicitly compute both sides (suitably scaled) in the \myemph{joint} distributional limit as $N\to\infty$ and~$n\to\infty$. Note that taking the limit jointly preserves pointwise equality.

Starting with the right hand side of \eqref{E:9.38}, the uniform continuity of $F$ and Corollary~\ref{cor-light}, we may rewrite it as a random quantity whose $L^1$-norm under $P^{\varrho} \otimes \BbbP$
is at most $o(1) n^{-1/2}$, with $o(1) \to 0$ as $n \to \infty$, plus the quantity
\begin{equation}
\label{E:9.39}
\frac{\sqrt{\log N}}{\wh W_N}
\sum_{x\in A_N}
f_n \biggl(L_{t_N}^{D_N} (x) + \frac{1}{2} (h_x^{D_N})^2 \biggr)
F \biggl(\Bigl\{L_{t_N}^{D_N} (x+z) + \frac{1}{2} (\phi^{(x)}_z)^2
\colon z \in \Bbb{Z}^2 \Bigr\} \biggr),
\end{equation}
where we denoted $A_N:=\{x\in\Z^2\colon x/N\in A\}$.
Decomposing the sum over $x$ with $L_{t_N}^{D_N} (x) = 0$
and the sum over $x$ with $L_{t_N}^{D_N} (x) > 0$
and applying Lemma~\ref{lemma-6.3} to the latter,
we rewrite (\ref{E:9.39})
as
\begin{equation}
\label{E:29}
\frac{\sqrt{\log N}}{\wh W_N} \sum_{x\in A_N}
1_{\{L_{t_N}^{D_N} (x) = 0 \}}\,
f_n \Bigl(\frac{1}{2} (h_x^{D_N})^2 \Bigr)
F \biggl(\Bigl\{L_{t_N}^{D_N} (x+z) + \frac{1}{2} (\phi^{(x)}_z)^2
\colon z \in \Bbb{Z}^2 \Bigr\} \biggr)
\end{equation}
plus a random quantity whose $L^1$-norm under $P^{\rho} \otimes \BbbP$
is at most $o(1) n^{-1/2}$ with $o(1) \to 0$ as $N \to \infty$
followed by $n \to \infty$.
Let $\wh \kappa^D$ be a (subsequential) weak limit of $\wh \kappa_N^D$ along the subsequence $\{N_k\}$. 
By Lemma~\ref{lemma-avoid-local-add}, as $k \to \infty$, \eqref{E:29} converges in law to
\begin{multline}
\qquad
\frac{1}{\sqrt{2 \pi g}} 
\int \wh \kappa^D (\textd x\,\textd\ell)\otimes \textd h \otimes \nu^0 (\textd\phi) 
1_A (x) f_n\bigl(\tfrac{h^2}2\bigr) 
F \biggl(\Bigl\{\ell_z + \frac{1}{2} \phi_z^2\colon z \in \Bbb{Z}^2 \Bigr\} \biggr) 
\\
= \frac{4}{3\sqrt{\pi g n}} \int \wh \kappa^D (\textd x\,\textd\ell)\otimes \nu^0 (\textd\phi) 
1_A (x)
F \biggl(\Bigl\{\ell_z + \frac{1}{2} \phi_z^2\colon z \in \Bbb{Z}^2 \Bigr\} \biggr) 
\qquad
\end{multline}
where we used the explicit form of~$f_n$ to perform the integral over~$h$. Multiplying this by $\frac{3}{4} \sqrt{\frac{n}{2}}$, as $n\to\infty$ this converges to
\begin{equation}
\label{E:8.43uiw}
\frac1{\sqrt{2\pi g}}\int \wh \kappa^D_A (\textd\ell)\otimes \nu^0 (\textd\phi) 
F \biggl(\Bigl\{\ell_z + \frac{1}{2} \phi_z^2\colon z \in \Bbb{Z}^2 \Bigr\} \biggr)
\end{equation}
as~$n\to\infty$ where $\wh\kappa^D_A(B):=\wh\kappa^D(A\times B)$. This is the $N\to\infty$ and~$n\to\infty$ limit of the (rescaled) right-hand side of \eqref{E:9.38}.

Concerning the left-hand side of \eqref{E:9.38}, whenever~$A$ is such that $\leb(\partial A)=0$ (which implies $Z^D_{\sqrt\theta}(\partial A)=0$ a.s.), \cite[Theorem~2.1]{BL4} yields convergence to
\begin{equation}
\label{E:9.42}
\cspecial(\sqrt{\theta}) Z_{\sqrt{\theta}}^D (A)
\int \textd h\otimes \nu_{\sqrt{\theta}} (\textd\phi)
\,\texte^{\alpha \sqrt{\theta} h} \,f_n\bigl(\tfrac{h^2}2\bigr)\, F\biggl(\Bigl\{\bigl(h-\tfrac{1}{2} \phi_z\bigr) (-\phi_z)\colon z \in \Bbb{Z}^2 \Bigr\} \biggr),
\end{equation}
where $\nu_{\sqrt{\theta}}$ is the law of $\phi+\alpha\sqrt{\theta}\, \fraka$ with~$\phi$ distributed according to~$\nu^0$.
Using that
\begin{equation}
\int \textd h\,\,\texte^{\alpha \sqrt{\theta} h} f_n\bigl(\tfrac{h^2}2\bigr) 
= \frac{4\sqrt{2}}{3\sqrt{n}} + O(n^{-3/2}),\quad n\to\infty,
\end{equation}
\eqref{E:9.42} multiplied by $\frac{3}{4} \sqrt{\frac{n}{2}}$ converges to
\begin{equation}
\label{E:8.46uiw}
\cspecial(\sqrt{\theta}) Z_{\sqrt{\theta}}^D (A)
\int \nu^0 (\textd\phi) F\Bigl(\bigl\{\tfrac{1}{2} (\phi_z+\alpha\sqrt\theta\,\fraka)^2\colon z \in \Bbb{Z}^2 \bigr\}\Bigr)
\end{equation} 
as~$n\to\infty$. This is the $N\to\infty$ and~$n\to\infty$ limit of the (rescaled) left-hand side of \eqref{E:9.38}.

We now finally have a chance to invoke the Pinned Isomorphism Theorem of~\cite{R18}. Indeed, since $2 \sqrt{2u} =\alpha\sqrt\theta$ implies $u=\pi\theta$, \eqref{E:3.6uiw} equates \eqref{E:8.46uiw} (and thus \eqref{E:8.43uiw})~with
\begin{equation}
\label{E:8.47uiw}
\cspecial(\sqrt{\theta}) Z_{\sqrt{\theta}}^D (A)
\int \nu_{\theta}^{\text{\rm RI}} (\textd\ell) \otimes \nu^0 (\textd\phi)
F\Bigl(\bigl\{\ell_z + \tfrac{1}{2} \phi_z^2\colon z \in \Bbb{Z}^2 \bigr\}\Bigr).
\end{equation} 
The Bounded Convergence Theorem extends the equality of \eqref{E:8.43uiw} and \eqref{E:8.47uiw} to~$F$ of the form $F(\ell):=\exp\{-\sum_{z\in\Lambda_r(0)}b_z\ell_z\}$ for any~$b_z\ge0$. This effectively transforms the term~$\frac12\phi_z^2$ away from both expressions and, thanks to the Cram\'er-Wold device, implies
\begin{equation}
\wh\kappa_A^D(\textd \ell)  \overset{\text{\rm law}}= \sqrt{2\pi g}\,\cspecial(\sqrt{\theta})Z^D_{\sqrt\theta}(A)\,\nu_{\theta}^{\text{\rm RI}}(\textd\ell).
\end{equation}
As this holds for all open~$A\subseteq   D$ with~$\ol A\subseteq   D$, the claim follows.
\end{proofsect}

\section*{Acknowledgments}
\nopagebreak\noindent
The first author has been supported in part by JSPS KAKENHI, Grant-in-Aid for Early-Career Scientists 18K13429. The second author has been partially supported by the NSF award DMS-1712632.  We wish to thank an anonymous referee for a number of important corrections to the initial version of this manuscript.

\bibliographystyle{abbrv}

\end{document}